\documentclass[11pt,leqno]{article}
\usepackage{amsthm,amsfonts,amssymb,amsmath,oldgerm}
\usepackage{epsfig}
\usepackage[font=footnotesize, format=hang, justification=justified]{caption}
\numberwithin{equation}{section} 

\newcommand{\Trace}{\hbox{\rm Trace}\,}

\def\eps{\varepsilon }

\def\eps{\varepsilon}

\newcommand\br{\begin{remark}}
\newcommand\er{\end{remark}}
\newcommand\bp{\begin{pmatrix}}
\newcommand\ep{\end{pmatrix}}
\newcommand\be{\begin{equation}}
\newcommand\ee{\end{equation}}
\newcommand\ba{\begin{equation}\begin{aligned}}
\newcommand\ea{\end{aligned}\end{equation}}

\newcommand{\bap}{\begin{app}}
\newcommand{\eap}{\end{app}}
\newcommand{\begs}{\begin{exams}}
\newcommand{\eegs}{\end{exams}}
\newcommand{\beg}{\begin{example}}
\newcommand{\eeg}{\end{exaplem}}
\newcommand{\bpr}{\begin{proposition}}
\newcommand{\epr}{\end{proposition}}
\newcommand{\bt}{\begin{theorem}}
\newcommand{\et}{\end{theorem}}
\newcommand{\bc}{\begin{corollary}}
\newcommand{\ec}{\end{corollary}}
\newcommand{\bl}{\begin{lemma}}
\newcommand{\el}{\end{lemma}}
\newcommand{\bd}{\begin{definition}}
\newcommand{\ed}{\end{definition}}
\newcommand{\brs}{\begin{remarks}}
\newcommand{\ers}{\end{remarks}}

\newtheorem{theo}{Theorem}[section]

\newtheorem{exams}[theo]{Examples}

\numberwithin{equation}{section}

\newcommand{\E }{{\mathcal{E}}}

\newtheorem{theorem}{Theorem}[section]
\newtheorem{proposition}[theorem]{Proposition}
\newtheorem{corollary}[theorem]{Corollary}
\newtheorem{lemma}[theorem]{Lemma}
\newtheorem{definition}[theorem]{Definition}

\newtheorem{example}[theorem]{Example}
\newtheorem{remark}[theorem]{Remark}

\title{
A numerical stability investigation of strong ZND detonations
for Majda's model}

\author{\sc \small 
Blake Barker
\thanks{Indiana University, Bloomington, IN 47405;
bhbarker@indiana.edu:
Research of B.B. was partially supported
under NSF grant no. DMS-0801745. }
~~~~
Kevin Zumbrun\thanks{Indiana University, Bloomington, IN 47405;
kzumbrun@indiana.edu:
Research of K.Z. was partially supported
under NSF grants no. DMS-0300487 and DMS-0801745.
 }}
\begin{document}

\maketitle


\begin{abstract}
We carry out a systematic numerical stability analysis of
ZND detonations of Majda's model with Arrhenius-type ignition
function, a simplified model for reacting flow, 
as heat release and activation energy are varied.
Our purpose is, first, to answer a question of Majda
whether oscillatory instabilities can occur for high activation
energies as in the full reacting Euler equations, and, second,
to test the efficiency of various versions of
a numerical eigenvalue-finding scheme
suggested by Humpherys and Zumbrun against
the standard method of Lee and Stewart.
Our results suggest that instabilities do not occur for Majda's model 
with Arrhenius-type ignition function, 
nor with a modified Arrhenius-type ignition function suggested by Lyng--Zumbrun,
even in the high-activation energy limit.
We find that the algorithm of Humpherys--Zumbrun is in the context
of Majda's model $100$-$1,000$ times
faster than the one described in the classical work of Lee and Stewart
and $1$-$10$ times faster than an optimized version
of the Lee--Stewart algorithm using an adaptive-mesh ODE solver.
\end{abstract}


\bigbreak

\section{Introduction}\label{s:majda}
In this paper, we carry out a systematic numerical stability
investigation for detonation solutions of Majda's model
with Arrhenius-type ignition function in the high-activation energy
limit, at the same time testing
and comparing
various different techniques for numerical stability analysis.
Our results should have application also to the effective design of numerical
methods for more complicated detonation models. 

Majda's model \cite{M1}, a simplified model for reacting compressible
gas dynamics, is often used as a testing ground for theory
and numerical methods designed for application to
the more complicated reacting compressible Euler equations \cite{BMR,M2}.
A question posed by Majda in \cite{M2} is whether 
this simplified model is sufficient to capture the complicated
Hopf bifurcation/pulsating instability phenomena
that occur for the full equations 
in the high-activation energy limit.\footnote{ 
``In particular, ...
there is the possibility of {\it Hopf bifurcation} to pulsating reacting
fronts with an associated exchange of stability...''-- problem 1, p.
25, \cite{M2}.
}
To explore this theoretical question is one motivation for the present work.

A second motivation comes from the numerics themselves.
In their foundational paper \cite{LS} on numerical stability analysis,
in which they introduce the standard scheme now in use,
Lee and Stewart describe the numerical determination of  
detonation stability as computationally intensive, and identify the
development of more efficient schemes as an important problem in
the physical detonation theory.\footnote{
``Finally, we point out that even though
our scheme is direct and easy to implement, complete 
investigation of the various regions of parameter space is 
computationally intensive. 
Any equivalent or more efficient numerical method for computation 
of detonation should be considered a valuable contribution 
and such approaches are needed to further explore the 
parameter regimes of instability.''-- closing paragraph, p. 130, \cite{LS}.}

We take advantage, therefore, of the simplified context of Majda's model,
as a proving ground for various different numerical schemes, in particular
testing whether recently-developed
techniques from the related problem of stability of viscous
shock waves \cite{HuZ2,HLZ,HLyZ,BHZ,Z3} can be imported
in a useful way.

Our main object is to test whether an alternative ``Evans function-type''
scheme proposed by Humpherys--Zumbrun \cite{HuZ1} can outperform the standard
scheme of Lee--Stewart, at the same time exploring optimal
implementations for both.
There is some reason to hope for improvement, since, as pointed
out in \cite{HuZ1}, the Lee--Stewart shooting method can be reformulated as
an adjoint Evans computation carried out in a backward direction,
the direction from $x=0$ toward $x=-\infty$ in which eigenfunctions
(normal modes) are required to decay.  (See also Section \ref{s:LS}.)
As pointed out in \cite{Br,BrZ,HuZ2,Z3}, there is a numerical advantage
in shooting, rather, in the direction $x=-\infty$ toward $x=0$
that eigenfunctions are expected to grow, since error modes then decay
exponentially relative to the mode being computed, with the advantage
being proportional to the spectral gap between growing and decaying modes
of the eigenvalue ODE.

Other novelties of our investigation are the systematic use of adaptive-step 
mesh both in spatial and frequency variables, the development of higher-order
high-frequency asymptotics, the derivation of rigorous if sometimes
conservative bounds on the maximum size of unstable eigenvalues,
and the introduction of a hybrid and limiting schemes designed for
the singular, square-wave limit.


\subsection{Results}
Our results are, first, that Majda's model does not appear to support
instabilities of any kind for Arrhenius or modified Arrhenius
ignition function, even in the high-activation energy limit.
(For zero activation energy, see the analytical proof of \cite{JY}.)

Second, it does appear that the optimum version found for
the Humpherys--Zumbrun scheme outperforms the optimum version found
for the Lee--Stewart scheme, by a factor of $1-10$ depending on
frequency, and on the average perhaps $2-4$.
However, this difference is dwarfed by the (scheme independent)
one obtained by using an adaptive-step
ODE solver in the $x$-evolution, which yields improvement over the
fixed-mesh solver originally used in \cite{LS} by a factor of $100-1000$.
Likewise, in our implementation, the difference between the implicit
$z$-coordinatization of \cite{Er2,LS},
in which the profile is explicitly solvable as a function of $z$, 
and the $x$-coordinatization 
in which the problem presents itself, is improvement by a uniform 
factor of $6$
(apparently due to the cost of evaluating the profile;
see Remark \ref{effic}).

The message is that the choice of numerical scheme requires a bit of care.  For,
with the wrong choices, performance can degrade by a total factor of
as much as $2400$-$24000$!
With all the right choices, on the other hand, the computation is for
reasonable values of activation energy quite numerically well-conditioned,
at least in the simplified context of Majda's model, with computation
times on the order of that seen for a scalar or $2\times 2$ viscous
conservation law, of a few seconds for an entire stability computation
for a given set of model parameters.

In the high-activation energy limit, as for any singular limit, the computational performance degrades.
In such cases, we find it necessary to factor out as much growth/decay of the solution as we possibly can, with less effective schemes not even converging 
for reasonable precisions and computation times.

\subsection{Discussion and open problems}

As observed by Majda \cite{M2}, the Majda model may be derived
from the full reacting compressible Euler equations via
weakly nonlinear geometric optics in certain limiting regimes.
However, as noted in \cite{Z1}, these can be seen to lie in the
small heat-release/small-activation energy
region for which detonations are known to be stable \cite{Z4}.
Thus, there is no physical reason to expect that parallels
should extend to the high-activation energy limit and associated
pulsating instability phenomena.

Nonetheless, the structural analogy between the equations persists,
and so the observed results of universal stability are somewhat 
striking.
It would be interesting to pursue further the difference between
the two models, both at formal and rigorous levels.
In particular, a very interesting open problem
is the analytic verification of stability for general
$\E$, as done for $\E=0$ in \cite{JY}.

A novelty here as regards detonation literature is the derivation
of rigorous bounds on the size of unstable eigenvalues.
However, these in some cases degrade rapidly in the singular, high-activation
energy limit, blowing up as $e^{O(\E)}$ as $\E\to \infty$.\footnote{
In some interesting cases, including the square-wave limit, they
do not; see Remark \ref{squareok}.}
Indeed, our convergence studies show these bounds to be much too
conservative, suggesting a sharp bound rather of $O(\E)$.
The application of semiclassical limit/turning point techniques
to obtain better bounds would be a
very interesting technical question; see Remark \ref{turnrmk}.
Another very interesting open problem would be to carry
out a complete stability analysis in the limit as $\mathcal{E}\to \infty$,
a problem intermediate to the analysis of general $\mathcal{E}$.
This together with the numerical studies carried out here for
bounded $\mathcal{E}$, would resolve the question of general
$\mathcal{E}$
by a combined numerical and analytical approach,
similarly as was done for viscous shock waves in \cite{HLZ,HLyZ,BHZ,BLZ}, 

We have carried out here a systematic confirmation/examination 
in a simple context of a number of aspects of numerical detonation 
stability analysis, which we hope will serve as a
useful reference for further developments.
In particular, a very interesting direction for future investigation
is to determine whether the gains in efficiency observed for
Majda's model carry over to the full reacting compressible
Euler equations.

\section{Preliminaries}
\subsection{Equations and assumptions}
Consider the inviscid Majda model
\ba\label{1.1}
u_t+\displaystyle\left(\frac{u^2}{2}\right)_x & =  q k \phi(u)z ,
\\
z_t & = -k\phi(u)z,\\
\ea
$u,z, \phi, q, k \in \mathbb{R}^1$, 
$q\ge 0$, $k>0$, with Arrhenius type ignition function
\begin{equation}\label{Ar}
\phi(u)=\begin{cases}
Ce^{-\mathcal{E}/T(u)}& T>0,\\
0& T\le 0.
\end{cases}
\end{equation}

Here $T(u)$ is a relation approximating the temperature/velocity relation
for a full ZND profile.  We consider here the simple cases $T(u)=u$ 
as proposed by Majda \cite{BMR,M1,M2} and $T(u)=1-(1-1.5)^2$, 
a downward quadratic relation qualitatively similar to that 
for the full ZND equations as proposed by Lyng--Zumbrun \cite{LyZ2}.

A strong detonation wave of \eqref{1.1} is a traveling-wave solution
\begin{equation}\label{1.2}
(u,z)(x,t)  =  (\overline{u},\overline{z})(x-st), \quad
\lim_{\xi\rightarrow\pm\infty}(\overline{u},\overline{z})(\xi)=(u_\pm,
z_\pm)
 \end{equation}
in the weak, or distributional, sense, smooth except at a single
shock discontinuity, without loss of generality at $x=0$,
where $u$ jumps from $u_*=\bar u(0^-)$ to $\bar
u(0^+)$ as $x$ crosses zero from left to right, with
    $z_-=0$,   $z_+=1$, $u_-> u_i >u_+$,
and
\begin{equation}\label{Lax}
u_->s>u_+,
\end{equation}
where $u_i$, defined as the minimum value of $u$ for which $T(u_i)=0$,
is the ignition temperature.

That is, a strong detonation wave consists of
a shock advancing to the right into a quiescent (i.e., nonreacting)
constant state with reactant mass fraction $z=1$, raising $u$ above
ignition level $u_i$, followed by a smooth ``reaction tail'' in
which combustion (reaction) occurs, in which $z$ decays
exponentially to value $z_-=0$ and $u$ to $1\le u_-<u_*$ as $x\to
-\infty$ \cite{Er1,M1,LyZ1,LyZ2,JLW,Z1,Z2}.

\subsection{Parametrization}
By the invariances of \eqref{1.1}, we may take without loss of generality
$u_*=2$, $u_+=0$, $s=k=1$, and $C>0$ arbitrary, leaving only the
parameters $q\ge 0$ and $\E\ge 0$. 
From \eqref{Lax} and the
assumption that $(\bar u,\bar z)$ converges as $x\to \pm \infty$, we
find \cite{JY,Z1} that
\begin{equation}\label{const}
(\bar u, \bar z)(x)\equiv (u_+,z_+)=(0,1) \; \hbox{\rm for } \; x> 0
\end{equation}
and also 
(see \eqref{profsoln} below)
$
 u_-= 1+\sqrt{1-2q} ,
$ so that
\begin{equation}\label{range}
1\le u_-\le 2, \quad 0\le q\le \frac{1}{2}.
\end{equation}
with activation energy $\mathcal{E}$ 
varying in the infinite range $0\le \E <+\infty$.

\subsection{Profiles}\label{S.3}

The ZND profile equations may be written, adding $q$ times
the second equation to the first, as
 $-s(\overline{u}+q\overline{z})^{\prime}
+ \left(\frac {\overline{u}^2}{2}\right)^{\prime} = 0 $
and 
  $-s\overline{z}'+k\phi(\overline{u})\overline{z}=0$.
Integrating the first equation from -$\infty$ to $x$ and solving
the resulting quadratic, we obtain (see \cite{JY,Z1} for details)
\be\label{profsoln}
\overline{u}(x)=1+ \sqrt{1-2q(1-\overline{z}(x))}.
\ee
We then obtain $\bar z$ by solving the ODE
\be\label{zODE}
z'=k\phi(u(z))z, \quad z(0)=1.
\ee
In the simplest case $\E=0$, we have the explicit solution $\bar z=e^{kx}$.

\subsection{Square-wave structure and the high-activation energy limit}
\begin{figure}[htbp]
\begin{center}
$
\begin{array}{lcr}
 (a) \includegraphics[scale=.3]{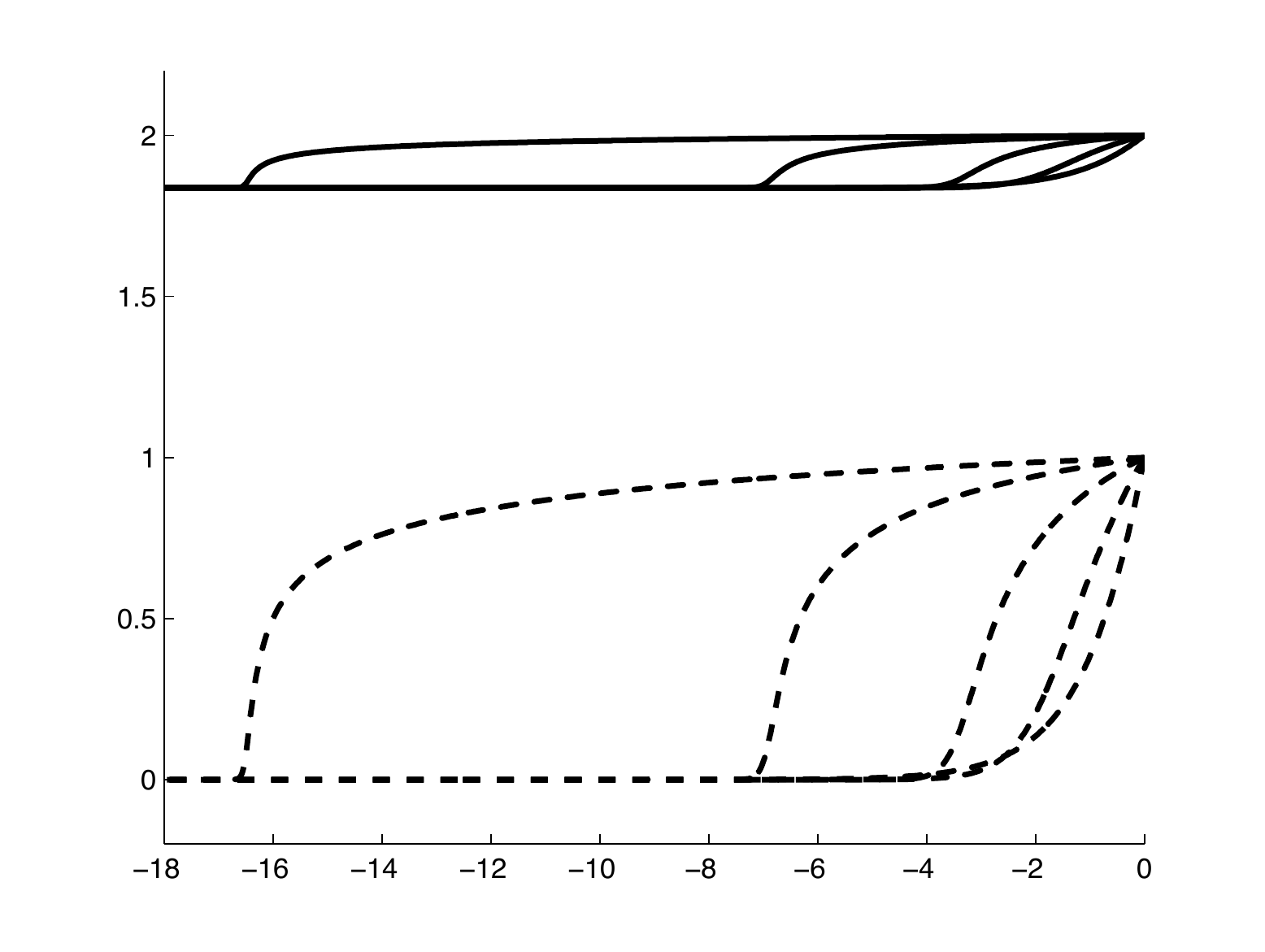}&(b)\includegraphics[scale=0.3]{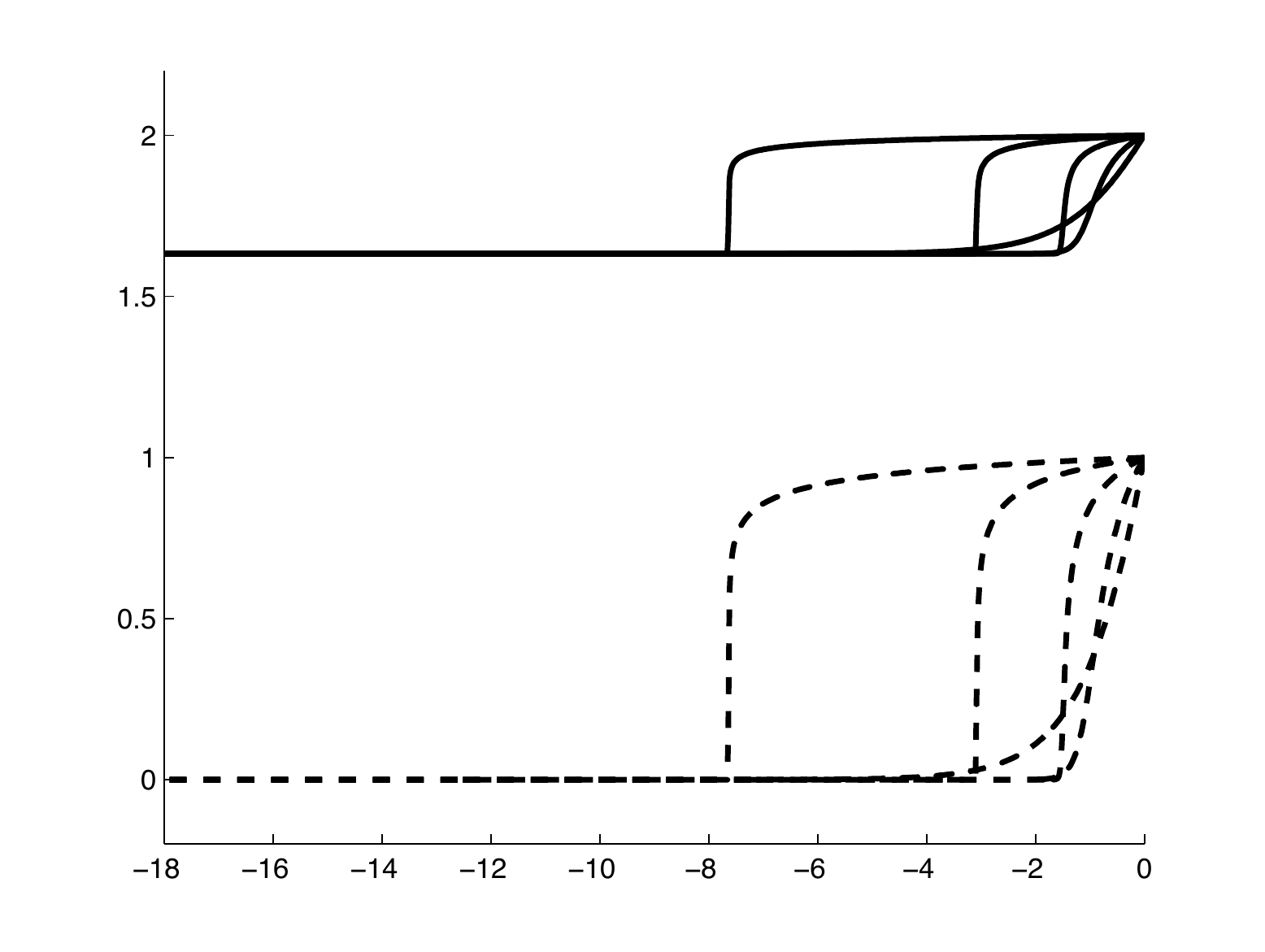} (c) \includegraphics[scale=0.3]{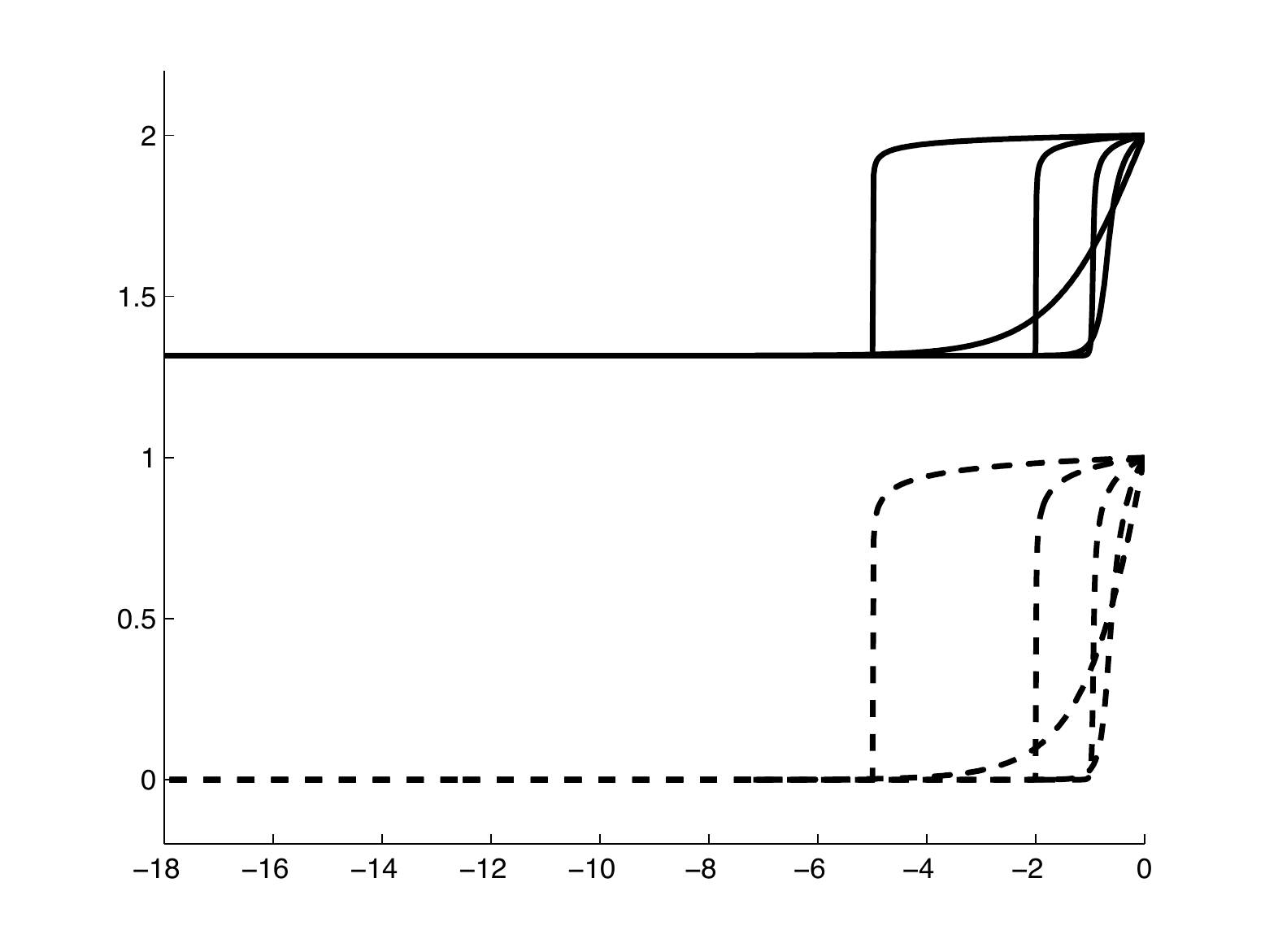}  \end{array}
 $
\end{center}
\caption{
Profile plots for fixed $q$ value and $\E=1,\ 10,\ 20,\ 30,\ 40$. We use a solid line for $\bar u$ and a dashed line for $\bar z$ plotted against $x$. Here $\phi(u)=Ce^{-\mathcal{E}/T(u)}$, $T(u)=1-(1.5-u)^2$, and $C=10^{21\mathcal{E}/40}$. In Figure (a), $q=0.15$; (b), $q=0.3$; (c), $q=0.45$. 
}
\label{profplots}
\end{figure}

\begin{figure}[htbp]
\begin{center}
$
\begin{array}{lcr}
 (a) \includegraphics[scale=.3]{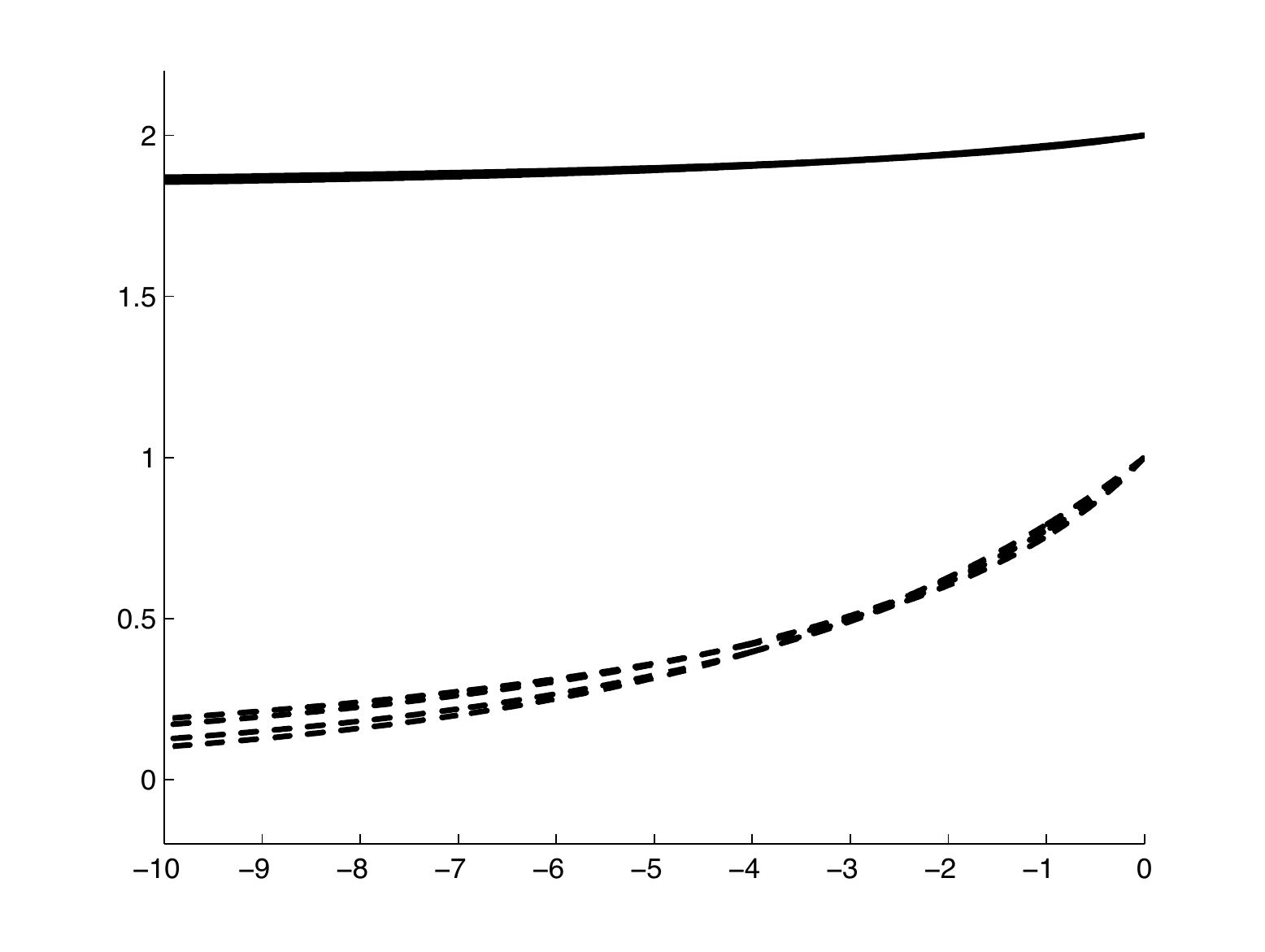}&(b)\includegraphics[scale=0.3]{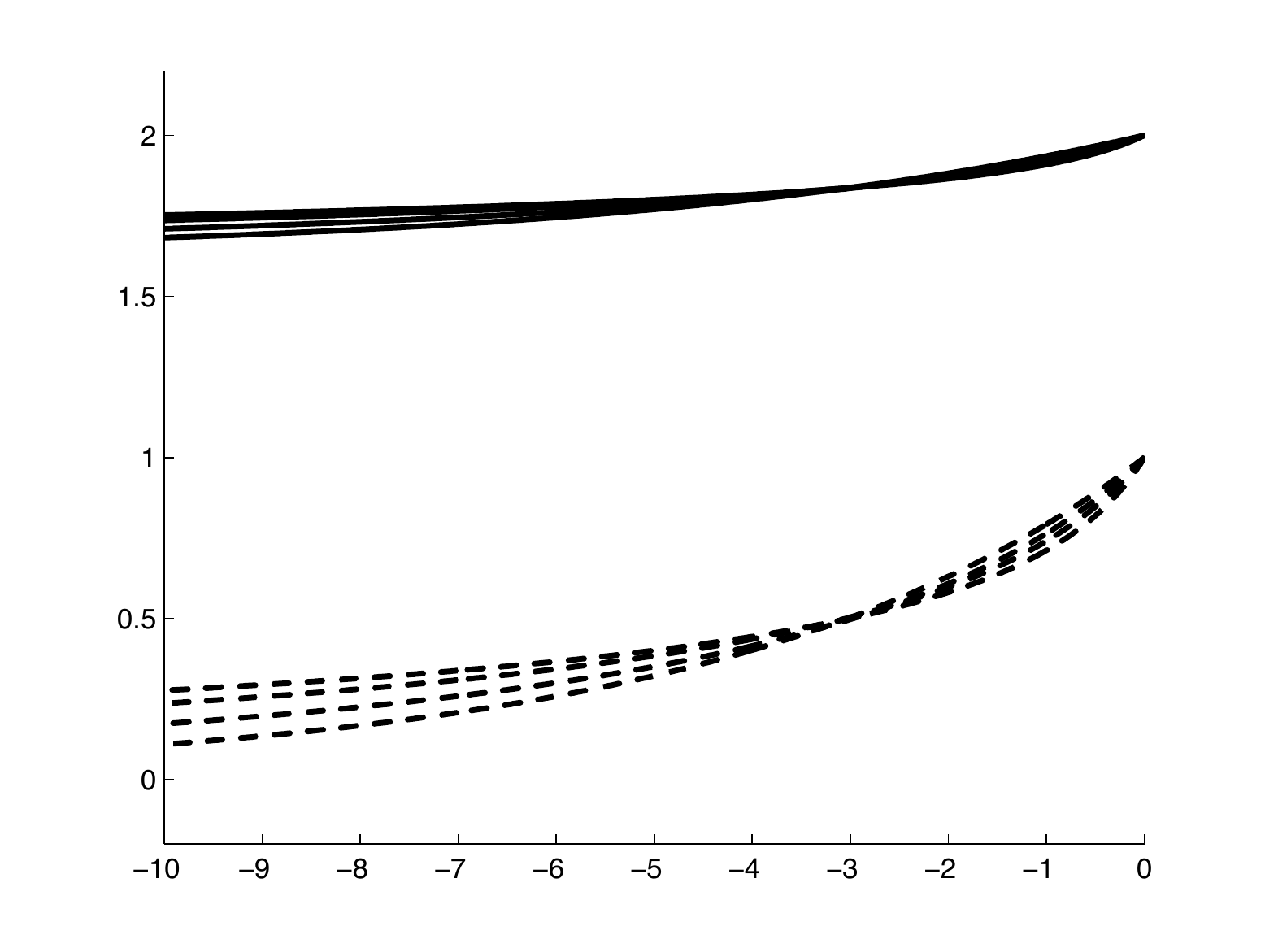} (c) \includegraphics[scale=0.3]{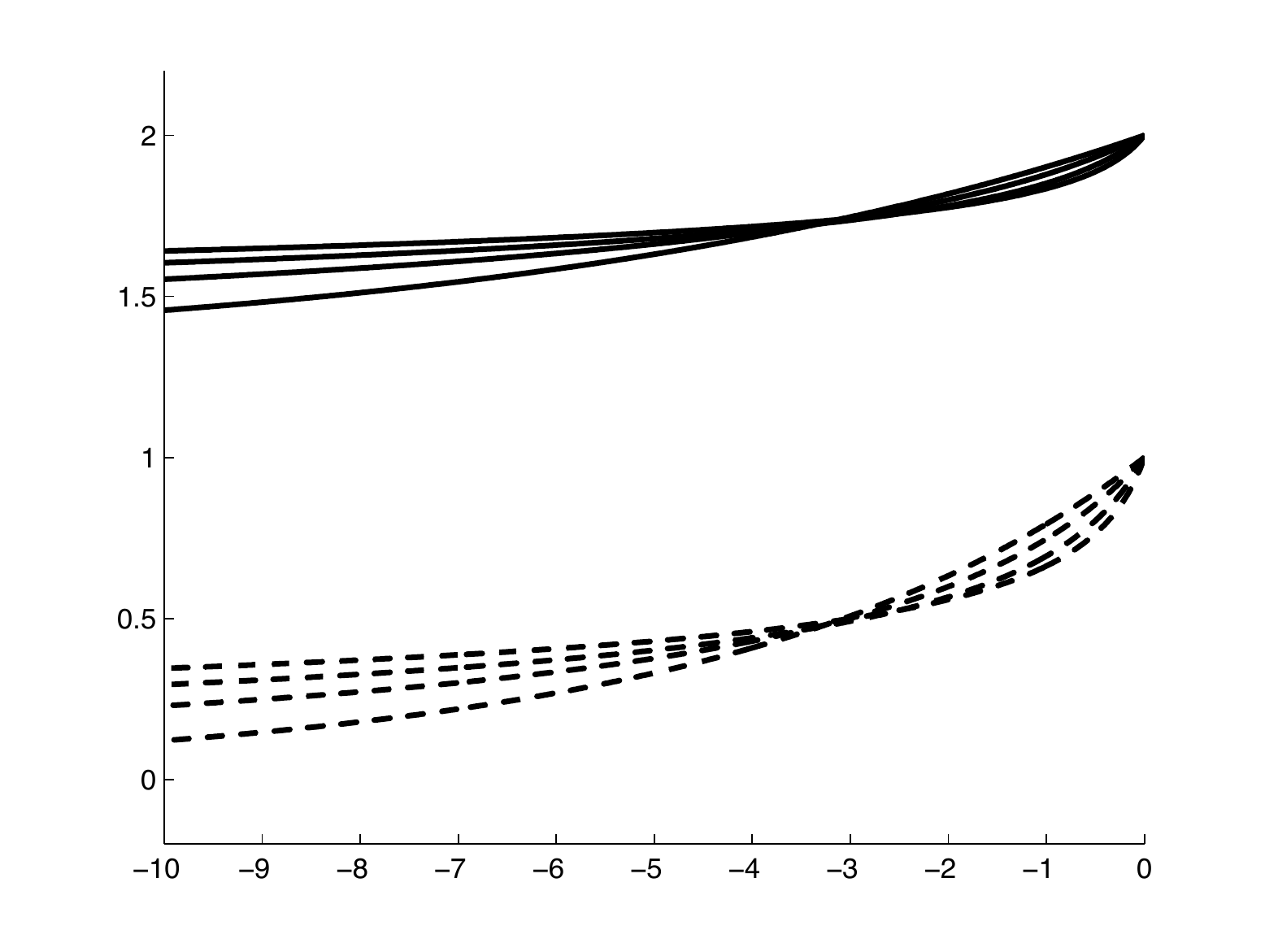}  \end{array}
 $
\end{center}
\caption{
Profile plots for fixed $q$ value and $\E=1,\ 10,\ 20,\ 30$. We use a solid line for $\bar u$ and a dashed line for $\bar z$ plotted against $x$. Here $\phi(u)=Ce^{-\mathcal{E}/u}$, where $C$ is chosen so that $\bar z(-2)=1/2$. In Figure (a), $q=0.15$; (b), $q=0.3$; (c), $q=0.45$. 
}
\label{fig3}
\end{figure}

In the high-activation energy limit $\E\to 0$, the profile ODE becomes
singular, with variation in $\bar z$ concentrated near the value
$u_{max}$ for which the ``temperature function''
$T(u)$ is maximized.
For the Arrhenius case $T(u)=u$, this means concentration near the
value $u_{max}=\bar u(0^-)=2$ at the right endpoint $x=0$ of the
profile, and results in a sharpened reaction spike near that point;
see Figure \ref{fig3}.
Here, following the standard normalization of \cite{Er2,LS}, we have
chosen $C$ in \eqref{Ar} so that $\bar z(-2)=1/2$, that is, the
half-reaction occurs at a specified spatial point $x=-2$.
For the modified Arrhenius case $T(u)=1-(u-1.5)^2$, rather,
$\phi(2)\to 0$ exponentially in $\E$, and so the profile takes on
a characteristic ``square-wave'' shape similarly as for the full
reacting Euler equations, consisting of a long flat tail from
$x=-\infty$ to an intermediate point $x_0$ for which $\bar u(x_0)=u_{max}=1.5$, 
a rapid change in
the vicinity of $x_0$, and a second long flat region from $x=x_0$
to  $x=0$;
see Figure \ref{profplots}.
Here, we have have chosen a bit less carefully the normalization
$C=e^{21\E/40}$ in order to keep in frame the value of $x$ at
which $z=1/2$. 
(Recall that change in $C$ amounts to rescaling in $x$, so is
not essential to our analysis.)
For the full equations, square-wave structure and the high-activation energy
limit are associated with Hopf bifurcation and
transition to instability \cite{Er1,Er2,LS}.

\section{Linearized stability and the Evans--Lopatinski determinant}\label{evansZND}
We now briefly review the linearized stability
theory of \cite{Er1,JLW,Z1,HuZ1}. 
Shifting to coordinates $\tilde x=x-st$ moving with the background
Neumann shock, write \eqref{1.1} as
$W_t + F(W)_x=R(W)$, where
\begin{equation}\label{abcoefs}
\begin{aligned}
W:=\begin{pmatrix} u\\z \end{pmatrix}, \quad
F:=\begin{pmatrix} 
u^2/2-su\\
-sz \end{pmatrix},\quad
R:=\begin{pmatrix} qkz\phi(u)\\-kz\phi(u) \end{pmatrix}.\\
\end{aligned}
\end{equation}
To investigate solutions in the vicinity of a discontinuous 
detonation profile, we postulate existence of a single shock
discontinuity at location $X(t)$, and reduce to a fixed-boundary
problem by the change of variables $x\to x-X(t)$.
In these coordinates, the problem becomes
$W_t + (F(W) - X'(t) W)_x=R(W)$, $x\ne 0$,
with jump condition
$X'(t)[W] - [F(W)]=0$, 
$[h(x,t)]:=h(0^+,t)-h(0^-,t)$ as usual denoting jump 
across the discontinuity at $x=0$.

\subsection{Linearization/reduction to homogeneous form}\label{linearization}
In moving coordinates, $\bar W^0$ 
is a standing detonation, hence 
$(\bar W^0, \bar X)=(\bar W^0,0)$ is a steady solution of 
the nonlinear equations.
Linearizing about $(\bar W^0,0)$, we obtain the {\it linearized equations} 
$(W_t - X'(t)(\bar W^0)'(x)) + (AW)_x =EW,$
with jump condition
$X'(t)[\bar W^0] - [A W]=0$ at $x=0$,
where
$A:= (\partial/\partial W)F$, $E:= (\partial/\partial W)R$.

Reversing the original transformation to linear order, following \cite{JLW},
by the change of variables $W\to W- X(t)(\bar W^0)'(x)$,
and noting that $x$-differentiation of 
the steady profile equation $F(\bar W^0)_x=R(\bar W^0)$ gives
$(A(\bar W^0)(\bar W^0)'(x))_x=E(\bar W^0)(\bar W^0)'(x)$,
we obtain modified, {\it homogeneous} interior equations
$W_t + (AW)_x =EW$
together with a modified jump condition
accounting for front dynamics of
$X'(t)[\bar W^0]-  [A \big(W+ X(t) (\bar W^0)'\big) ]=0$.

\subsection{Evans--Lopatinski determinant}\label{S.4}
Seeking normal mode solutions $W(x,t)=e^{\lambda t}W(x)$,
$X(t)=e^{\lambda t}X$,
$W$ bounded, of the linearized homogeneous equations,
we are led to the generalized eigenvalue equations
$(AW)' = (-\lambda I   + E)W$ for $x\ne 0$, and
 $X(\lambda[\bar W^0]-[A (\bar W^0)']) - [A W]=0$, 
where ``$\prime$'' denotes $d/dx$, or, setting $Z:=AW$, to
\be\label{eig}
Z' = GZ, \quad x\ne 0,
\ee
\begin{equation}\label{eigRH}
\begin{aligned}
X(\lambda[\bar W^0]-[A (\bar W^0)']) - [Z]=0, 
\end{aligned}
\end{equation}
with
\begin{equation}\label{4.3}
A=\left(
    \begin{array}{cc}
      \overline{u}-1 & 0 \\
      0 & -1 \\
    \end{array}
  \right),
\;
  E=\left(
      \begin{array}{cc}
        qkd\varphi(\overline{u})\overline{z} & qk\varphi(\overline{u}) \\
        -kd\varphi(\overline{u})\overline{z} & -k\varphi(\overline{u}) \\
      \end{array}
    \right),
\end{equation}
and 
\be \label{G}
 G=(E-\lambda I)A^{-1}=
\bp
\frac{qk\, d\phi(\bar u) \bar z -\lambda}{\bar u-1}& -qk\phi(\bar u)\\
\frac{-k\, d\phi(\bar u)\bar z}{\bar u-1}& k\phi(\bar u) + \lambda
\ep. 
\ee

The Lopatinski determinant is then defined as
\begin{equation}\label{5.1}
    \displaystyle D_{ZND}(\lambda):=\det (Z^-(\lambda,0),
\; \lambda[\overline{W}]+R(0^-))|_{x=0},
\end{equation}
where $Z^-(\lambda, \xi)$ is a bounded
exponentially decaying solution of \eqref{eig}, analytic in
$\lambda$ and tangent as $\xi \rightarrow -\infty$ to the subspace
of exponentially decaying solutions of the limiting,
constant-coefficient equations $Z^{\prime}=G_-Z$, 
and 
$$
  R(0^-)=\left(
      \begin{array}{cc}
qk\varphi(\overline{u}) \\
-k\varphi(\overline{u}) \\
      \end{array}
    \right),
   \quad
[\bar W]=
\left(
                 \begin{array}{c}
                   -2 \\
                   0 \\
                 \end{array}
               \right).
$$
More precisely, $Z^-\sim e^{\mu_- x}V_-(\lambda)$, 
where $\mu_-$ is the positive real part eigenvalue of $G_-$ and
$V_-$ is an associated eigenvector chosen analytically in $\lambda$
 \cite{Br,HuZ2,Z3}.
Evidently, {\it there exists a normal mode solution with frequency $\lambda$,
$\Re \lambda \ge 0$, if and only if $D_{ZND}(\lambda)=0$}.

\section{Numerical approximation}
To estimate $D_{ZND}$ numerically, we approximate $Z^-$
at a large but finite negative value $x=-M$ with the value
$Z^-(-M)=e^{-\mu_- M}V_-$, and solve from $-M$ to $0$ using
a standard adaptive-step Runga--Kutta scheme.
The vector $V_-(\lambda)$ is evolved analytically by solution
of Kato's ODE $\partial_\lambda V=P'(\lambda)V$, where
$P(\lambda)$ is the eigenprojection associated with $\mu_-$;
see \cite{BrZ,Z2,HuZ2,Z3}.
Numerical convergence and efficiency of general schemes of this
type are discussed in \cite{Br,HuZ2,Z3}.

\br
As discussed in \cite{HuZ1,Z3}, the single most important
factor for efficiency of such computations is to use an adaptive-step
rather than fixed-step ODE solver since traveling front and boundary-layer
solutions inherently involve multiple scales.
As seen here, failure to use adaptive steps can reduce 
efficiency by two or more orders of magnitude.
\er

\medskip
{\bf Detection of roots.}
Zeros of the Evans--Lopatinski determinant can be found through individual
$\lambda$-evaluations by Newton's or other root-finding/following methods
as in \cite{LS}.  Alternatively, as in \cite{Er2,BrZ,HuZ2,Z4}, they
can be detected by a Nyquist diagram, or winding-number, computation,
mapping a large semicircle $S$ contained
in the positive real part half-plane $\Re \lambda\ge 0$ via $D_{ZND}$
and taking the winding number to determine the number of zeros lying
within $S$.
We follow the latter approach here.
We discretize $S$ by adaptive $\lambda$-steps, taking
care that the relative change in $D_{ZND}(\lambda_j)$ is $\le 0.2$
for each step,
thus ensuring an accurate winding number count \cite{Br,BrZ}.

\medskip

{\bf $z$-coordinatization and renormalization.}
For purposes of numerical approximation, is advantageous to use
the coordinatization of Erpenbeck, Lee--Stewart, and others,
by $z$ instead of $x$.  This means replacing $Z'=GZ$ by
\be\label{zversion}
\dot Z= \check G Z,
\ee
where $\dot { }$ denotes $d/dz$ and
\be\label{H}
\check G(z):= \frac{G}{\bar z'} =
\frac{G}{k\phi(u(z))} ,
\ee
then integrating from $z=0$ to $z=1$ instead
of from $x=-\infty$ to $x=0$.
In practice, we integrate from $z=e^{-k\phi(-\infty)M}$ to $z=1$,
where $x=-M$ is our usual starting point in $x$-coordinatized version, 
and initialize $Z^-$ as usual as $e^{-\mu M}V_-$.
This avoids the need to solve the $\bar z$ equation 
$\bar z'=k\phi(\bar u)\bar z$ numerically.

A further (standard; see \cite{Br,HuZ2}) improvement is to renormalize $Z=e^{\int_0^x \mu }Y$, 
dividing out expected growth/decay,
converting $Y'=(G-\mu I)Y$ to
$\dot Y=HY$, where
\be\label{H2}
H(z):= \frac{G-\mu I}{\bar z'} =
\frac{G-\mu I}{k\phi(u(z))} .
\ee
Here, $\mu$ can be either $\mu_-$ 
or (better if computable analytically) the eigenvalue $\mu(x)$
of $G(x)$.

\br\label{effic}
For an adaptive-step ODE solver, we found no mathematical difference
between $x$- and $z$- coordinates, as both required the same number
of mesh points/functional evaluations for a given $x$ (resp $z$) 
integration.\footnote{Additional experiments, not recorded here.}
However, in practice, the change to $z$-coordinates gave an improvement
of $6$ times or more in speed
due to the cost of the interpolation step used to evaluate the numerically
pre-computed profile at the variable points needed for an adaptive step ODE.
This could be improved somewhat by the use of a more efficient interpolation
scheme; however, $z$-coordinates are always preferable for an adaptive-step
ODE solver.
As discussed in \cite{Z3,HuZ1}, renormalization typically improves speed
for a single $\lambda$-evaluation by a factor of $2$ or more.
This is magnified in the high-activation energy limit, for which growth/decay
rates become extreme, and the unrenormalized computation often cannot even be carried out, leaving machine scale and returning NaN errors.
For winding number computations, it is still more important to divide out
expected growth/decay, which otherwise introduces additional winding, decreasing the $\lambda$ stepsize and greatly increasing computational cost.
\er

\begin{figure}[htbp]
\begin{center}
 \includegraphics[scale=.3]{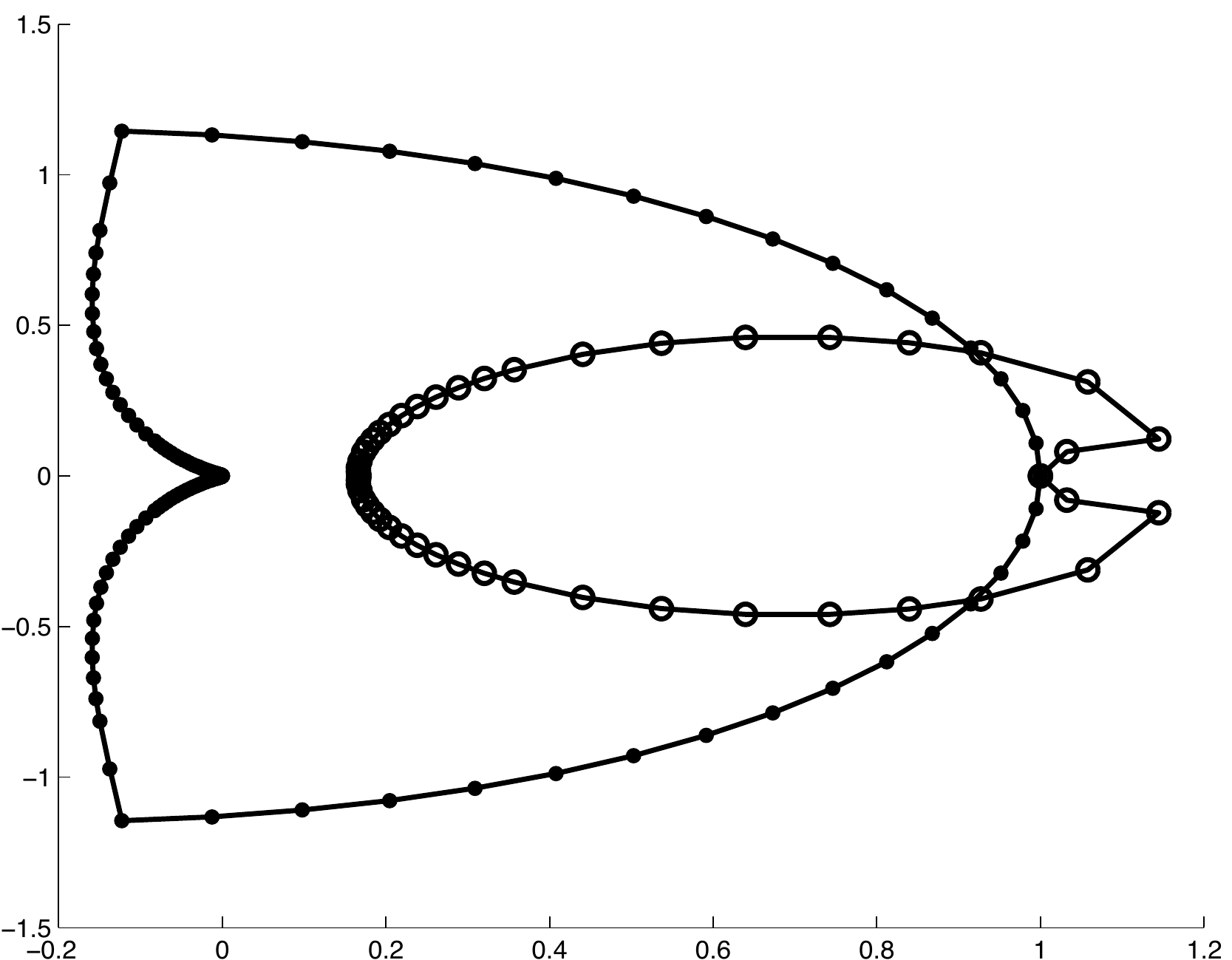}
\end{center}
\caption{
For $\E=10$, $q=0.3$, $C=10^{E/2+E/40}$, and $R=4$, 
which we get from the high-frequency convergence study, it takes 
$2.7524$ seconds and $127$ mesh points to obtain $0.2$ 
tolerance in the $\lambda$ contour for the unreduced Evans--Lopatinski
determinant $D_{ZND}$. 
If we divide by $\lambda$, it takes $0.6858$ seconds and $28$ 
mesh points.  Open circles correspond to the reduced
Evans function for which we divide by $\lambda$.
}
\label{redfig}
\end{figure}

\medskip
{\bf Reduction by $\lambda$.}
By translation-invariance, $D_{ZND}$ has always a root at $\lambda=0$.
Likewise, there is an extra factor $\lambda$ as $\lambda \to \infty$ induced
by the form of factor $(\lambda[\bar W^0]+R)$ in the defining determinant
beyond what is induced by growth/decay of the ODE solution $Z^-$;
see also the more detailed discussion of Section \ref{s:hf}.
For both these reasons, it is advantageous in performing winding-number
computations to work with a {\it reduced Evans--Lopatinski determinant}
$D_{ZND}(\lambda)/\lambda$, effectively removing a single zero, hence
one circuit about the origin, and thereby greatly reducing the number
of $\lambda$-points required for the computation (typically by factor $4$
or so); see Figure \ref{redfig}.
All winding number computations are done, therefore, with respect to
the reduced determinant, throughout the paper.

\section{Alternative formulations}\label{s:alt}
\subsection{The adjoint Evans--Lopatinski determinant}
We may define an alternative stability function (see \cite{HuZ2,CJLW,Z4})
as 
\begin{equation}\label{adevanseqmajda}
\tilde D_{ZND}(\lambda)=
\langle \tilde Z,  (\lambda[\bar W^0]+R) \rangle|_{x=0},
\end{equation}
where $\tilde Z$ denotes the unique (up to scalar multiplier) 
decaying solution of
\be\label{eigmajda}
\tilde Z' = -G^* \tilde Z, \quad x\le 0,
\qquad
 -G^*=
\bp

\frac{-qk\, d\phi(\bar u) \bar z +\lambda^*}{\bar u-1} &
\frac{k\, d\phi(\bar u)\bar z}{\bar u-1}\\
 qk\phi(\bar u)&
-k\phi(\bar u) - \lambda^*\\
\ep. 
\ee
$D_{ZND}$ and $\tilde D_{ZND}$ differ by a nonvanishing
analytic factor $\beta$, by 
duality relations $\tilde Z= \frac{Z^\perp}{\mathcal{W}}$ 
and $\det(v,w)=v^\perp \cdot w$,
where $\mathcal{W}=ce^{\int_0^x \Trace (G)}$ 
is an appropriately normalized Wronskian of \eqref{eig}
and $D_{ZND}=\beta \tilde D_{ZND}$ with $\beta= \mathcal{W}$.

\br\label{erprmk}
The method of reduction to homogeneous form \cite{JLW} is equivalent
to solving the inhomogeneous equation 
\be\label{inhom}
\hat Z'=G\hat Z+ \lambda X \bar W'(x)
\ee
by linear superposition, using the
fact that a particular solution is given by $W_p= X\bar W'(x)$,
hence $\hat W=W^- + X\bar W'(x)$, or $\hat Z=Z^-+XA\bar W'(x)=Z^--XR(x)$,
and substituting into the jump relation
$\lambda X [\bar W^0]+ [\hat Z]=0$.
\er

As discussed further in \cite{HuZ1}, \eqref{eigmajda} may be viewed as
a streamlined version of the original method of Erpenbeck \cite{Er2};
we shall refer to this scheme as the {\it homogeneous Erpenbeck} method.
When renormalized by $\tilde \mu_-$ (resp. $\tilde \mu(x)$), we will
refer to it as the {\it adjoint $\mu$} (resp. {\it adjoint $\mu(x)$})
method.

\subsection{The method of Lee and Stewart}\label{s:LS}

The method of Lee and Stewart consists, rather, of solving
the inhomogeneous equation \eqref{inhom} with $X=1$ from $x=0$
initialized with $\hat Z(0):= \lambda [\bar W]$.
They then take the inner product of $\tilde V^-$ against $Z(-M)$,
with $M$ chosen sufficiently large.\footnote{
More precisely, since
they work like Erpenbeck in $z$ and not $x$ coordinates,
against $Z|_{z=\eps}$, $\eps>0$ sufficiently small.
This is phrased 
in terms of a progress variable $\lambda:=1-z$;
what they call $\alpha$ is our $\lambda$.} 
By duality, this can be seen to be exactly
\be\label{LSevans}
D_{LS}(\lambda)=
\frac{\tilde D_{ZND}(\lambda)}{ e^{-\bar \mu_2(-\infty)(\lambda) M}},
\ee
where $\mu_2(-\infty)$ is the negative eigenvalue of $G_-$.


Expanding $\bar W'=(d\bar W/dz)\bar z'=
(d\bar W/dz)k\phi \bar z$, we may rewrite \eqref{inhom} 
in $z$-coordinates as
\be\label{Zz}
    \dot { \hat Z}=
\frac{G}{k\phi \bar z} \hat Z+  \lambda \bp d\bar u/d\bar z\\1\ep ,
\ee
where $d\bar u/d\bar z$ is determined from the 
profile solution \eqref{profsoln}.
Equivalently
\be\label{simpleLS}
\tilde D_{ZND}(\lambda)= e^{-\bar \mu_2(-\infty)(\lambda) M} D_{LS}(\lambda)=
\langle \tilde V^-,\hat Z_c(-\infty)\rangle,
\ee
where
\be\label{Zzcentered}
    \dot { \hat Z}_c=
\frac{G-\mu_2(-\infty)}{k\phi \bar z} \hat Z_c
+  \lambda \bp d\bar u/d\bar z\\1\ep .
\ee

\medskip
{\bf Homogeneous version.}
The Lee--Stewart method may equally well be
implemented via the homogeneous equations 
$Z'=GZ$, $Z(0)=\lambda[\bar W]+R$.
(Here, we are using $Z-\hat Z\to 0$ as $x\to -\infty$.)
This appears to be slightly faster but essentially equivalent
in practice.

\subsection{The Polar/Drury method}

An alternative to direct renormalization is the polar method 
(or continuous orthogonalization method of Drury) \cite{HuZ2},
which, in $z$ coordinates, appears as
\be\label{checkH}
\dot Y = \frac{(G- YY^* G)Y}{\bar z'} =
\frac{(G- YY^* G)Y}{k\phi( u(\bar z))\bar z} ,
\ee
where $D_{polar}(\lambda):=\det (Y,\lambda [\bar W^0]+R)|_{z=1}$.
This differs from $D_{ZND}$ by a nonvanishing continuous factor $r(\lambda)$,
so has the same winding number and roots.
It has the advantage that $|Y|\equiv 1$, so growth/decay is completely
factored out, while good numerical conditioning is maintained
so long as there remains at least a neutral spectral gap between 
$\mu_-$ and the other eigenvalue of $G_-$ \cite{Z3}, as holds
in this case on a neighborhood of $\Re \lambda\ge 0$ (see Section \ref{s:hf}).
An adjoint version can be computed similarly; we refer to this
scheme as the {\it polar adjoint} method.
The factor $r(\lambda)$ may be computed by 
integrating $\dot{\log r}= \frac{ Y^* G Y-\mu_-}{k\phi( u(\bar z))\bar z} $
from $z=0$ to $z=1$, with $r(0)=1$, and used to recover $D_{ZND}$ \cite{HuZ2}.
We refer to the latter scheme as the {\it polar radial}, and the
analogous dual as the {\it polar adjoint radial} method.

\subsection{Hybrid and limiting methods}

In the high-activation energy limit $\mathcal{E}\to 0$ with modified 
Arrhenius-type ignition function, $\bar u$ approaches the
well-known ``square-wave'' profile familiar from the full ZND equations
\cite{FD,Er1}, which consists approximately of three constant zones
separated by a reaction jump following at some distance the
Neumann shock.
This singular configuration is difficult for the standard numerical
Evans--Lopatinski function to resolve, since the integration passes through an
unexpected second, different, constant region for which it is not
tuned.
Indeed, following the numerical prescription of \cite{Z3}, we should
for best numerical conditioning
rather perform an Evans function type computation integrating from
$-\infty$ forward toward the reaction jump and from $0$ backwards
toward the reaction jump, taking a determinant in between.

This can conveniently be accomplished by a {\it hybrid method}
intermediate to the adjoint and Lee--Stewart methods, in which
we compute from $z=0$ to $z=1/2$ by the adjoint method, in some
optimized version, and from $z=1$ to $z=1/2$ via the Lee--Stewart
method, again in optimized version, taking the inner product at
$z=1/2$.  By duality, this gives the same result as the 
adjoint method with same optimized form.\footnote{
In practice, we solve $Z'=(G-\mu_2(x))Z$, where $\mu_2$ is the negative
eigenvallue of $G$, and the adjoint version $\tilde Z'=-(G-\mu_2(x))^*\tilde Z$.}
{\it This should be done only when $\phi(2)<<\max \phi$, so that
such singular square-wave structure occurs}.  When this is not the
case, the adjoint method is expected to be preferred.  However, the
hybrid method is safe, in that its numerical conditioning should always
lie somewhere between that of the adjoint and Lee--Stewart methods, and
captures some of the gain of the adjoint method even when square-wave
structure does not occur.  


\medskip
{\bf Limiting method.}  
When $\phi(2)<<\max \phi$, and square-wave structure
is truly in effect, we may take a limit {\it before} computing the hybrid
determinant and simply project out decaying modes.  That is, we may
define a projective, limiting, method, in which we substitute for 
initial data $\lambda[\bar W]+R$ the data
$P(\lambda)(\lambda[\bar W]+R)$, where $P(\lambda)$ is the eigenprojection
of $G|_{z=1}$ associated with the negative real part eigenvalue $\mu_2$.
For related analytical results in a similar situation, see \cite{Z5}.
We will refer to this version as the {\it Evans function} method.

\section{High frequency asymptotics}\label{s:hf}
Making the change of coordinates $x\to \hat x:= |\lambda|x$,
$\lambda \to \hat \lambda:=\lambda/|\lambda|$, we convert
\eqref{eigmajda} to the approximately diagonal system
\be \label{hfeig}
\dot { \tilde Z}=A(\eps \hat x)\tilde Z+ \eps B(\eps \hat x)\tilde Z,
\ee
where $|\hat \lambda|=1$, $\Re \hat \lambda \ge 0$,
 $\eps:= |\lambda|^{-1}$, $\dot{ }$ denotes $d/d\hat x$, and
\be\label{AB}
 A=
\bp
\frac{ \hat \lambda^* -\eps qk\, d\phi(\bar u) \bar z }{\bar u-1} & 0\\
0& - \lambda^*  -\eps k\phi(\bar u) \\
\ep,
\qquad 
 B=
\bp
0& \frac{k\, d\phi(\bar u)\bar z}{\bar u-1}\\
 qk\phi(\bar u)& 0\\
\ep. 
\ee

Making the further change of coordinates $\tilde Z=TY$,
$$
T=\bp 1 & \eps a\\
\eps b
& 1
\ep,
\qquad
T^{-1}=
(1-\eps^2ab)^{-1}
\bp 1 & -\eps a \\
-\eps b& 1
\ep,
$$
\be\label{abc}
a=\frac{\ k\, d\phi(\bar u)\bar z}{
-\bar u \hat \lambda^* -\eps(k\phi(\bar u)+ qk\, d\phi(\bar u)\bar z)},
\qquad
b=\frac{qk\phi(\bar u)(\bar u-1)}
{ \bar u \hat \lambda^* +\eps(k\phi(\bar u)+ qk\, d\phi(\bar u)\bar z)},
\ee
we obtain 
\be\label{Yeq}
\dot Y= T^{-1}(A+\eps B)TY -T^{-1}\dot T  Y
= A_1Y + \eps^2 B_1 Y + \eps^3 B_2Y,
\ee
where
\be\label{A1}
 A_1=
\bp \alpha_1&0\\0& \alpha_2\ep
:=
\bp
\frac{ \hat \lambda^* -\eps qk\, d\phi(\bar u) \bar z }{\bar u-1} & 0\\
0& - \lambda^*  -\eps k\phi(\bar u) \\
\ep,
\ee
\ba\label{B1}
 B_1&= 
\bp -a qk\phi(\bar u)
+\frac{bk\, d\phi(\bar u)\bar z}{\bar u-1} -ab(\lambda^*  +k\phi(\bar u)) 
& -a^2qk\phi(\bar u)+

\frac{ a_x}{ 1-\eps^2ab}\\
\frac{ b_x}{ 1-\eps^2ab}+
 \frac{b^2k\, d\phi(\bar u)\bar z}{\bar u-1}&
a qk\phi(\bar u) -\frac{bk\, d\phi(\bar u)\bar z}{\bar u-1} 
-ab \frac{ \hat \lambda^* -qk\, d\phi(\bar u) \bar z }{\bar u-1}
\ep,\\
B_2&=
\bp
\frac{ -ab_x}{ 1-\eps^2ab}&0\\
0
& \frac{ ba_x}{ 1-\eps^2ab}\\
\ep.
\ea


\subsection{Estimation of $Z^-$}\label{s:I}
Change now to $\hat y$-coordinates, where 
$d\hat x/d\hat y= \phi (\bar u(\hat x)$, 
i.e., $d\bar z/d\hat y= \eps k\bar z$, or $\bar z(y)=e^{k\eps \hat y}$.
Defining $Y(\hat y)=e^{\int_0^{\hat y}\alpha_1(\hat z)d \hat z} V(\hat y)$, 
and simplifying, we obtain as usual the Duhamel representation:
\be\label{Duhamel}
V(\hat y)= \mathcal{T}V(y):= V_-+ \int_{-\infty}^{\hat y}
e^{\int_{\hat z}^{\hat y} (\hat A_1-\hat \alpha_1 I)(\hat s)d\hat s}
(\eps^2 \hat B_1+\eps^3\hat B_2)V(\hat z)  d\hat z,
\ee
where $\hat A=A/ \phi$, $\hat \alpha_j=\alpha_j/\phi$,
$\hat B_j=B_j/ \phi$, and 
$|\hat B_1(\hat y)|\le C_1e^{-k \eps |\hat y|}$,
$|\hat B_2 (\hat y)|\le C_2e^{-k \eps |\hat y|}$
for some $C_j>0$,
and $V_-=V_-(\hat \lambda)=c(\hat \lambda)(1,0)^T$ is
the asymptotic limit at $\hat x\to -\infty$, $c\ne 0$.
Here, by the Mean Value Theorem,
$C_j\le \sup |d\hat B_j/dz|$.
From
$
\alpha_2-\alpha_1=  \frac
{- \bar u \hat \lambda^* -\eps k\phi(\bar u)+ \eps qk\, d\phi(\bar u)\bar z)}
{ \bar u-1}, 
$
we have by positivity of $\eps k\phi(u_-)$, $\Re \hat \lambda \bar u$,
 and $\bar u-1$, the bound
\ba\label{expA}
|e^{\int_{\hat z}^{\hat z} (\hat A_1-\hat \alpha_1)(\hat z)d\hat z}|
\le 
e^{\int_{\hat z}^{\hat y} \Re (\hat \alpha_2-\hat \alpha_1)(\hat z)d\hat z}
\le 
e^{ \gamma \int_{\hat z}^{\hat y} \frac{d \bar z}{d\hat z}
d\hat z}
\le e^{\gamma \bar z(\hat y)}
\le e^{\gamma } ,
\ea
\be\label{gamma}
\gamma:=
q \sup_x  \Big( \frac{ (d\phi /\phi)^{\mathcal{P}}} {\bar u-1}\Big),
\ee
where $f^{\mathcal{P}}$ denotes the positive part of $f$.
For $\Re \lambda \ge \max \frac{ qkd\phi\bar z-k\phi}{\bar u}$, we have simply
\ba\label{expA2}
|e^{\int_{\hat z}^{\hat z} (\hat A_1-\hat \alpha_1)(\hat z)d\hat z}|
&\le 
e^{\int_{\hat z}^{\hat y} \Re (\hat \alpha_2-\hat \alpha_1)(\hat z)d\hat z}
\le  1.
\ea

\br\label{gammarmk}
For the Arrhenius-type ignition function $\phi(u)=Ce^{\frac{-\E}{u}}$,
we have 
\be\label{Arr_ration}
d\phi/\phi=\E/u^2>0, 
\quad
\hbox{\rm hence $\gamma= q\E/u_-^2\le q\E$.}
\ee
For the modified Arrhenius ignition function $\phi(u)=Ce^{-\frac{\E}{T(u)}}$,
where $T(u)= 1-(u-1.5)^2$, we have
\be\label{mod_ration1}
d\phi/\phi=-2(\E/T^2)(u-1.5),  
\quad
\hbox{\rm hence $\gamma \le (4/3)q\E$ for $q\ge .375$}
\ee
and
\be\label{mod_ration}
\hbox{\rm $ \gamma\equiv 0$ for $q\le .375$, 
in which case $\bar u\ge u_-\ge 1.5$ and $d\phi/\phi \le 0$.}
\ee
A sufficient condition for 
$\Re \lambda \ge \max \frac{ qkd\phi\bar z-k\phi}{\bar u}$,
is
\be\label{suff}
\Re \lambda \ge k\gamma \max \phi \sim k\E\max \phi.
\ee
\er

\bl\label{cont}
For 
$\eps\le  \frac{k}{2e^{\gamma} (C_1+\eps C_2)}$ or
$\eps\le  \frac{k}{2 (C_1+\eps C_2)}$ and 
$\Re \lambda \ge \max \frac{ qkd\phi\bar z-k\phi}{\bar u}$,
$\mathcal{T}$ is a contraction on $L^\infty(-\infty,0]$,
and $\tilde Z^-|_{\hat x=0}$ 
is proportional to $(1,\omega_2)^T$,
where 
$ \omega_2:=\frac{\eps b(0) + \omega}{1+\eps a(0)\omega} $
and $|\omega|\le 
 \frac{\eps e^{\gamma}}{k} (C_1+\eps C_2)$.
\el

\begin{proof}
From \eqref{Duhamel} and the stated bounds, the Lipshitz norm
of $\mathcal{T}$ is bounded by
\ba\label{Duhamelbd}
\int_{-\infty}^{\hat y}
 e^{\gamma}
\frac{\eps^2(C_1+\eps C_2)}{\phi(\hat z} e^{-k \eps |\hat z|} d\hat z
&\le
 \eps e^{\gamma} (C_1+\eps C_2)
\int_{-\infty}^{\hat y}
e^{-k\eps |\hat z|} d\hat z\\
&\le 
 \eps e^{\gamma} (C_1+\eps C_2) 
\int_{-\infty}^{\hat y}
e^{-k \eps |\hat z|} d\hat z\\
&\le
 \frac{\eps e^{\gamma}}{k} (C_1+\eps C_2) e^{-k |\hat x|},
\ea
hence it is contractive for
 $\frac{\eps e^{\gamma}}{k } (C_1+\eps C_2)\le 1/2$, or
$\eps\le  \frac{k}{2 e^{\gamma} (C_1+\eps C_2)}$,
with relative error at $\hat x=0$ from the diagonal solution 
(i.e., with $\eps$ set to zero) given by
 $\frac{\eps e^{\gamma}}{k} (C_1+\eps C_2).$

That is, $Y(0)$ is proportional to $(1,\omega)^T$,
with
$|\omega|\le 
 \frac{\eps e^{\gamma}}{k} (C_1+\eps C_2)$.
Transforming back by $\tilde Z^-=TY$, with
$ T=\bp 1 & \eps a\\ \eps b & 1 \ep $, we thus find that
$$
\tilde Z^-|_{\hat x=0}=
 T=\bp 1 & \eps a\\ \eps b & 1 \ep \bp 1\\ \omega\ep=
\bp 1+\eps a \omega \\ \eps b+  \omega \ep
$$
is proportional to $(1,\omega_2)^T$,
where 
$ \omega_2:=\frac{\eps b(0) + \omega}{1+\eps a(0)\omega} $.
\end{proof}

\br\label{vargaprmk}
Lemma \ref{cont} may be recognized as a variable-coefficient analog
of the gap lemma of \cite{GZ}, and a quantitative version of the abstract
Lemma A.1 established for the full ZND system in \cite{Z4}.
Though we did not state it, 
the argument implies also that
$\tilde Z^-$ converges exponentially in relative error as $x\to -\infty$
to the solution of the diagonal system $\tilde Z'=A_1\tilde Z$, $A_1$ 
as in \eqref{A1}.
A similar argument applies for the related polar method \cite{HuZ2}.
\er

\bc\label{corI}
$D_{ZND}(\lambda) \ne 0$ for $0\le \Re \lambda$
and 
$$
 |\lambda|\ge \max \{
 \frac{2e^{\gamma}}{k} (C_1+\eps C_2),
\,
\frac{k\phi(2)}{2}(q  +|\omega_2|)\},
\qquad
 \omega_2:=\frac{\eps b(0) + \omega}{1+\eps a(0)\omega} \, .
$$
or
$$
\Re \lambda \ge \max \frac{ qkd\phi\bar z-k\phi}{\bar u}
\;
\hbox{ and } \;
 |\lambda|\ge 
\frac{k\phi(2)}{2}(q  +|\omega_2|),
\qquad
 \omega_2:=\frac{\eps b(0) + \omega}{1+\eps a(0)\omega} \, .
$$
\ec

\begin{proof}
Computing
$
  \lambda[\bar {W}]+R(\bar W(0^-))= (-2\lambda+qk\phi(2) , -k\phi(2) )^T, $
we have that $D_{ZND}(\lambda)$ is proportional by a
nonvanishing analytic factor to
$$
\begin{aligned}
\langle \tilde Z^-(\lambda,0), \; 
(-2\lambda+qk\phi(2) , -k\phi(2) )^T\rangle|_{ x=0}
&=
(1,\omega_2)^T \cdot (-2\lambda+qk\phi(2) , -k\phi(2) )^T,
   \end{aligned}
$$
and thus to
$
(1,\omega_2)^T \cdot (-2\hat \lambda+\eps qk\phi(2) , -\eps k\phi(2) )^T=
-2\hat \lambda+\eps k\phi(2)(q  -\omega_2),
$
which is nonvanishing provided
$ \eps k\phi(2)(q  +|\omega_2|)\le 2,$
or
$ \eps\le \frac{2}{ k\phi(2)(q  +|\omega_2|)}$.
\end{proof}

\br
When $|\lambda|\ge 
 \frac{2e^{\gamma}}{k} (C_1+\eps C_2)$,
$|\omega|\le 1$,
so that $|\omega_2|$ is roughly $1$ as well, and
$\frac{k\phi(2)}{2}(q  +|\omega_2|)\lesssim \frac{3k\phi(2)}{4}$.
\er

\br\label{turnrmk}
Corollary \ref{corI} is efficient for $\lambda$ near the real
axis.  However, for the Arrhenius-type ignition function
$T(u)=u$ it can be quite inefficient away from the real axis,
where it does not make sufficient use of spectral separation
of modes as opposed to spectral gap, hence requires the unusable
bound $|\lambda| \ge e^{\gamma}\sim e^{\E}$ as $\mathcal{E}\to \infty$.  
To obtain practical (guaranteed) bounds by use of 
more modern turning-point theory
is an important direction for future investigation.\footnote{
For the Arrhenius-type ignition function, actual instabilities
for $\mathcal{E}\to \infty$ are expected to appear (if indeed they do)
only in the region $|\lambda| \le \max qkd \phi \bar z \sim E\max \phi$
where spectral separation $|\alpha_1-\alpha_2|>>1$ fails.
}
\er

\br\label{squareok}
For the modified Arrhenius ignition function (quadratic version),
our estimates are efficient for $q\le .375$, where $\gamma=0$: 
in particular for
the ``difficult'' case $q=0.3$ we have much investigated; see Remark
\ref{gammarmk}.
\er

\br
For $q$ bounded from $1/2$, $k=1$ and $\E \equiv 0$, $C=1$, we have $\gamma=0$, 
$a\equiv 0$, and $b$, $b_x$ roughly $1$,
giving $C_2\sim C_1\sim 1$.
Thus, $\eps\lesssim 1$ is sufficient in this case for contraction,
or $|\lambda|\gtrsim 1$.
Our bounds blow  up in the CJ limit $q\to 1/2$, 
for which $(\bar u-1)\to 0$ as $x\to -\infty$.
(However, note that $b\to 0$ as $q\to q_{cj}$, partly compensating 
for badness of this limit.)
This case would be interesting for further investigation.
\er


\subsection{Limiting behavior}\label{s:limit}
By arguments like those of Section \ref{s:I} and
especially Remark \ref{vargaprmk}, we find that 
$$
Z^-(0)\sim 
\bp e^{\int_{-\infty}^0 (\beta_2(\hat y)-\beta_2(-\infty) d\hat y} \\ 0\ep
=
\bp e^{k\int_{-\infty}^0 (\phi(\bar u(y))-\phi(\bar u(-\infty))) dy} \\ 0\ep
$$
as $\eps \to 0$, where
$\beta_2:=-\alpha_2^*= \lambda  + \eps k\phi(\bar u) $,
$\alpha_2$ as in \eqref{A1}, yielding the following asymptotic behavior.

\begin{proposition}\label{hfbehavior}
For $\Re \lambda \ge 0$,
\be\label{hfconv}
D_{ZND}(\lambda) = C\lambda(1+O(\lambda^{-1})) ,
\qquad 
\tilde D_{ZND}(\lambda)= C_1 \lambda e^{C_2 \lambda}(1+O(\lambda^{-1})).
\ee
as $|\lambda|\to \infty$ ,
where $C:= 2 e^{k\int_{-\infty}^0 (\phi(\bar u(y))-\phi(\bar u(-\infty))) dy} $,
$C_1$, and $C_2$ are real constants.
\end{proposition}

\br\label{hfrmk}
\textup{
Higher order approximants
$
D_{ZND}(\lambda)=  e^{C_1 \lambda +C_0 + C_{-1}\lambda^{-1})}\lambda 
(1 + O(\lambda^{-2})),
$
etc., may be obtained by further 
diagonalizations as discussed in \cite{Z4,MaZ3}.}
\er

In particular, $D\sim 2\lambda$ for $\E=0$.
These relations can be used to determine a maximum radius through a convergence
study, which, as pointed out in \cite{HLZ,HLyZ},
in practice typically gives better bounds than those obtained
by rigorous tracking/conjugation bounds.

\section{Numerical experiments}


\subsection{Initialization at $-\infty$}
In our numerical experiments we compute the Lopatinski determinant given in 
\eqref{5.1}.
Following \cite{Z2,Z3}, we use the ODE 
\be\label{Kode}
S'=P'S
\ee
 of Kato to explicitly compute an analytically varying initializing eigenvector
$\tilde V^-$ for the Evans function, where $'$ denotes $d/d\lambda$,
$S(\lambda)$ is the desired eigenvector of $G_{-\infty}$, and $P(\lambda)$ is the
associated eigenprojection.
Preserving analyticity in this manner, we are able to employ the argument principle thus determining the number of zeros of the Evans function inside a contour by computing the winding number. At $x=-\infty$,
\be
G_{-\infty}=\begin{pmatrix}
a \lambda &b\\ 0&c+\lambda
\end{pmatrix},
\ee
where $a=\frac{-1}{\bar u-1}=\frac{-1}{\sqrt{1-2q}}$, $b=-qk\phi_-$, and $c=k\phi_-$. Then right and left eigenvectors of $G_{-\infty}$ corresponding to the eigenvalue $c+\lambda$, which satisfies $\Re(c+\lambda)>0$ for $\Re \lambda>0$, are respectively
$r_+= \begin{pmatrix} \frac{-b}{(a-1)\lambda -c} \\1 \end{pmatrix}$, 
$l_+=\begin{pmatrix}  0&1\end{pmatrix}$.
We then have the projection $P=\frac{r_+l_+}{l_+ r_+}$, and it's derivative (with respect to $\lambda$) $P'$,
$P=\begin{pmatrix} 0&\frac{-b}{(a-1)\lambda-c}\\0&1 \end{pmatrix}$,
$P'=\begin{pmatrix} 0&\frac{b(a-1)}{((a-1)\lambda -c)^2}\\
0&0 \end{pmatrix}$.
Solving the ODE $S'=P'S$, we find that
\be\label{Sform}
S(\lambda)=\alpha_2\begin{pmatrix} \frac{-b}{(a-1)\lambda -c}\\ 1 \end{pmatrix}+\alpha_1\begin{pmatrix}1\\0\end{pmatrix},
\ee
where we take $\alpha_1=0$ and $\alpha_2=1$.
The corresponding initializing vector for the adjoint method is
$\tilde V_-=\begin{pmatrix} 1\\ \frac{b}{(a-1)\lambda^*-c}\end{pmatrix}$.

\subsection{Computation of $Z^-$/evaluation of $D$}



Following the general approach of \cite{HuZ2,HLZ,BHZ},
the ODE calculations for individual $\lambda$ are carried out using \textsc{MATLAB}'s {\tt ode45} routine, which is the adaptive 4th-order Runge-Kutta-Fehlberg method (RKF45).  This method is known to have excellent accuracy with automatic error control.  Our standard
error tolerance setting is {\tt AbsTol = 1e-6} and {\tt RelTol = 1e-8}
unless otherwise mentioned, and 
the value of approximate minus spatial infinity $-M$ 
is determined experimentally by the requirement that the absolute error 
$$
|(\bar u,\bar z)(-M)-(u,z)_{-\infty}|\leq TOL
$$
be within a prescribed tolerance, 
say $TOL=10^{-3}$, where $(\bar u,\bar z)$ and $(u,z)_{-\infty}$ 
are the profile and limiting endstates of Section \ref{S.3}.
In the rescaled $z$ coordinates where we do not even need the profile to compute, we perform a convergence study of the Evans--Lopatinski determinant for $\lambda=2i$ to determine numerical infinity, requiring relative error between output be less than $10^{-3}$ between successive computations for numerical infinity decreasing as negative powers of 2.
For a theoretical
convergence analysis in the simple case considered here, see \cite{Br}; 
for a more general treatment, see \cite{Z3}.
For the results of a numerical convergence study, 
see Fig. \ref{domainconv} below.

\medskip
{\bf Fixed-mesh Lee-Stewart implementation.}
For comparison purpose, we carry out also experiments using
a fixed-mesh grid as prescribed in \cite{LS},
taking $\eps=1/N$,
and using a uniform mesh for $z\in [0,1]$ of mesh points $z_j=j/N$,
then choosing $N$ large enough to get 
a prescribed level of convergence, as determined by numerical convergence
study.

\begin{center}
\begin{figure}[htbp]
\begin{center}
$
	\begin{array}{cc}
(a)\includegraphics[scale=.3]{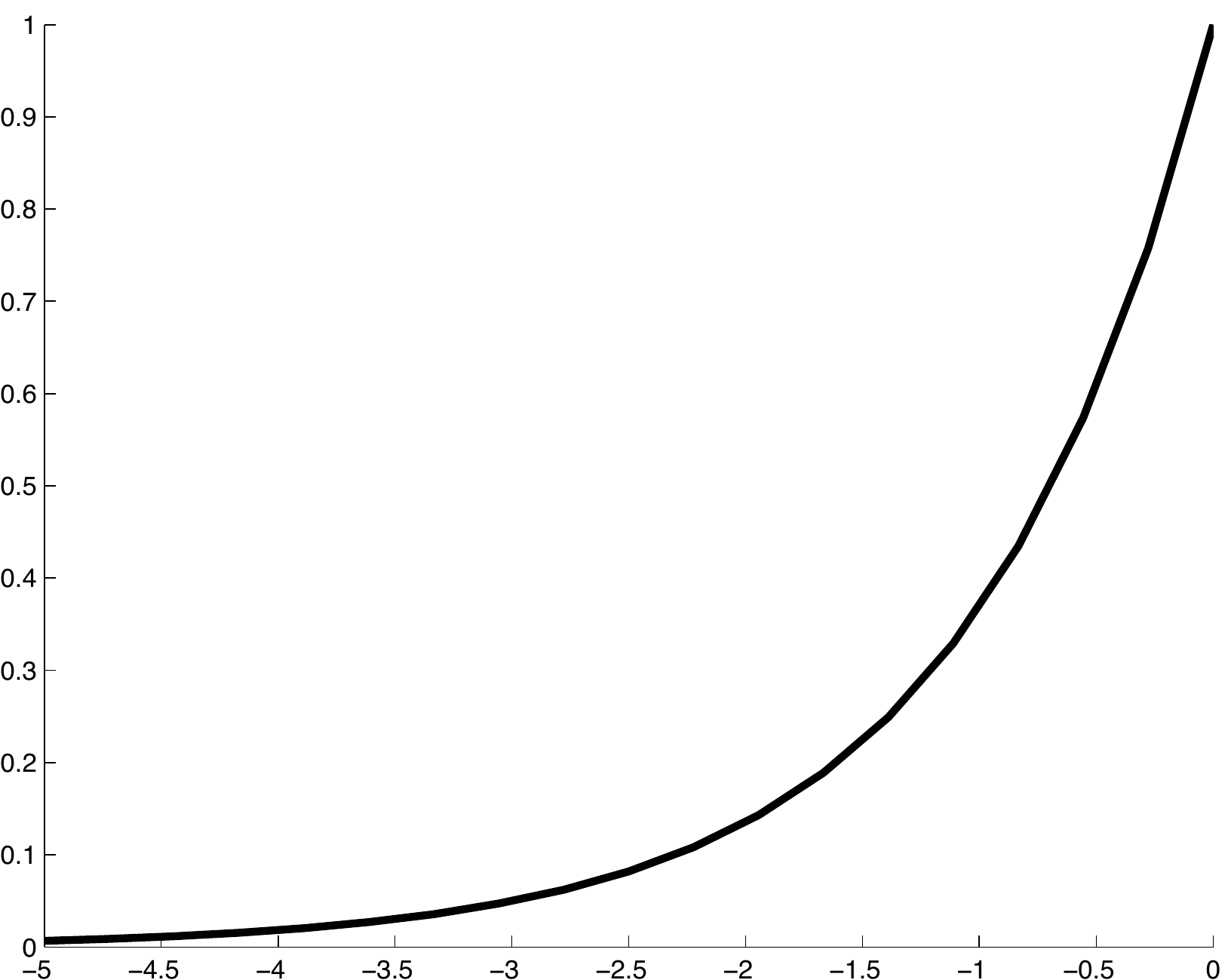}&
   		(b) \includegraphics[scale=.3]{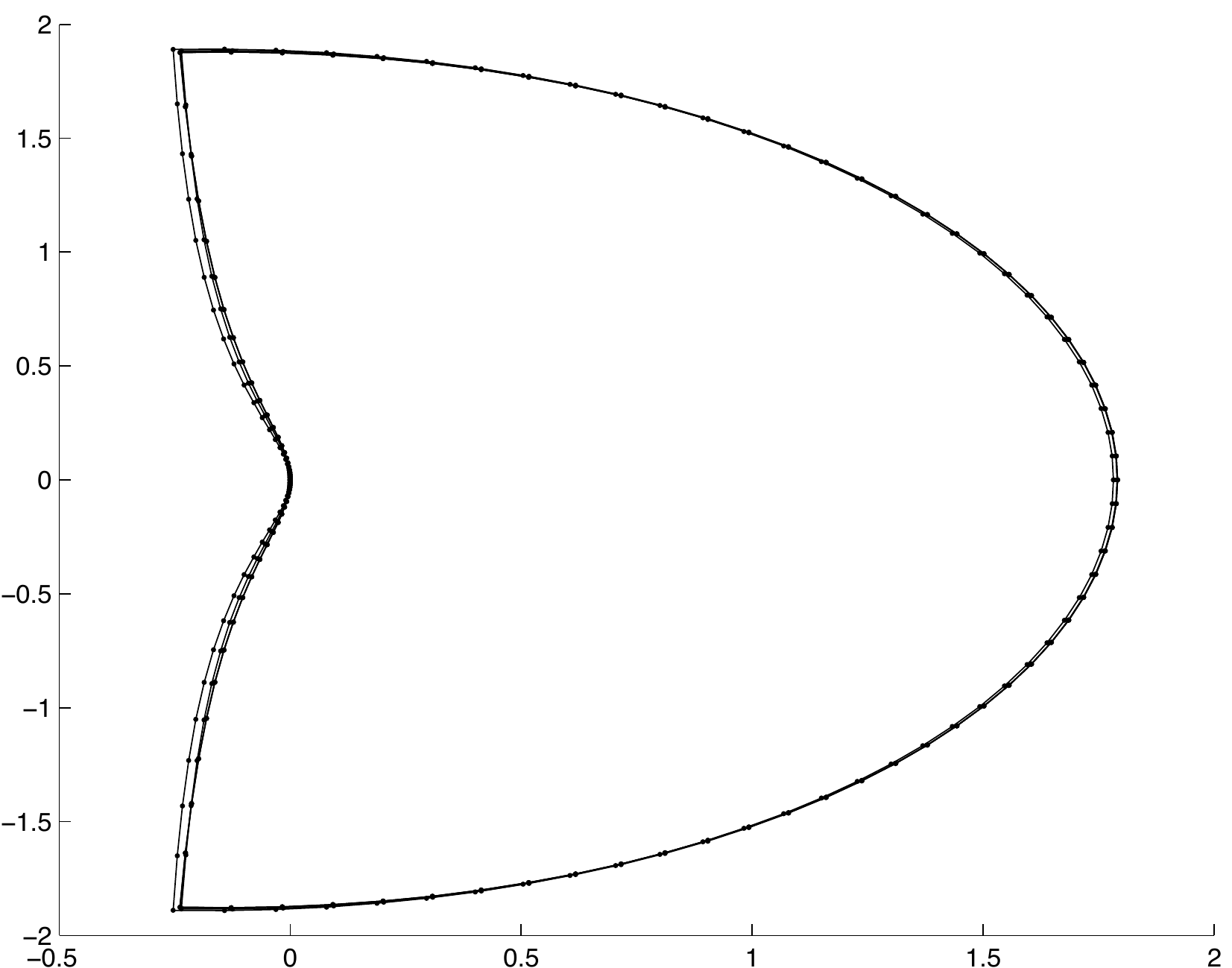}
    	\end{array}
$
\end{center}
 \caption{
a) Convergence of the profile. b) Convergence of Evans approximations (regular method with $\mu(x)$ scaled out)
on a contour with $R=1$, for $\phi\equiv 1$, $q=0.3$, $k=1$ and
approximate spatial infinity $M=0.9,\ 0.5,\ 10^{-1},\ 10^{-2},\ 10^{-3}$. Run times were respectively 1.08, 0.40, 0.76, 1.06, 1.33 seconds.
}
\label{domainconv}
\end{figure}
\end{center}

\subsection{Verification: comparison with exact solution}

\begin{center}
\begin{figure}[htbp]
\begin{center}
$
	\begin{array}{cc}
(a)\includegraphics[scale=.3]{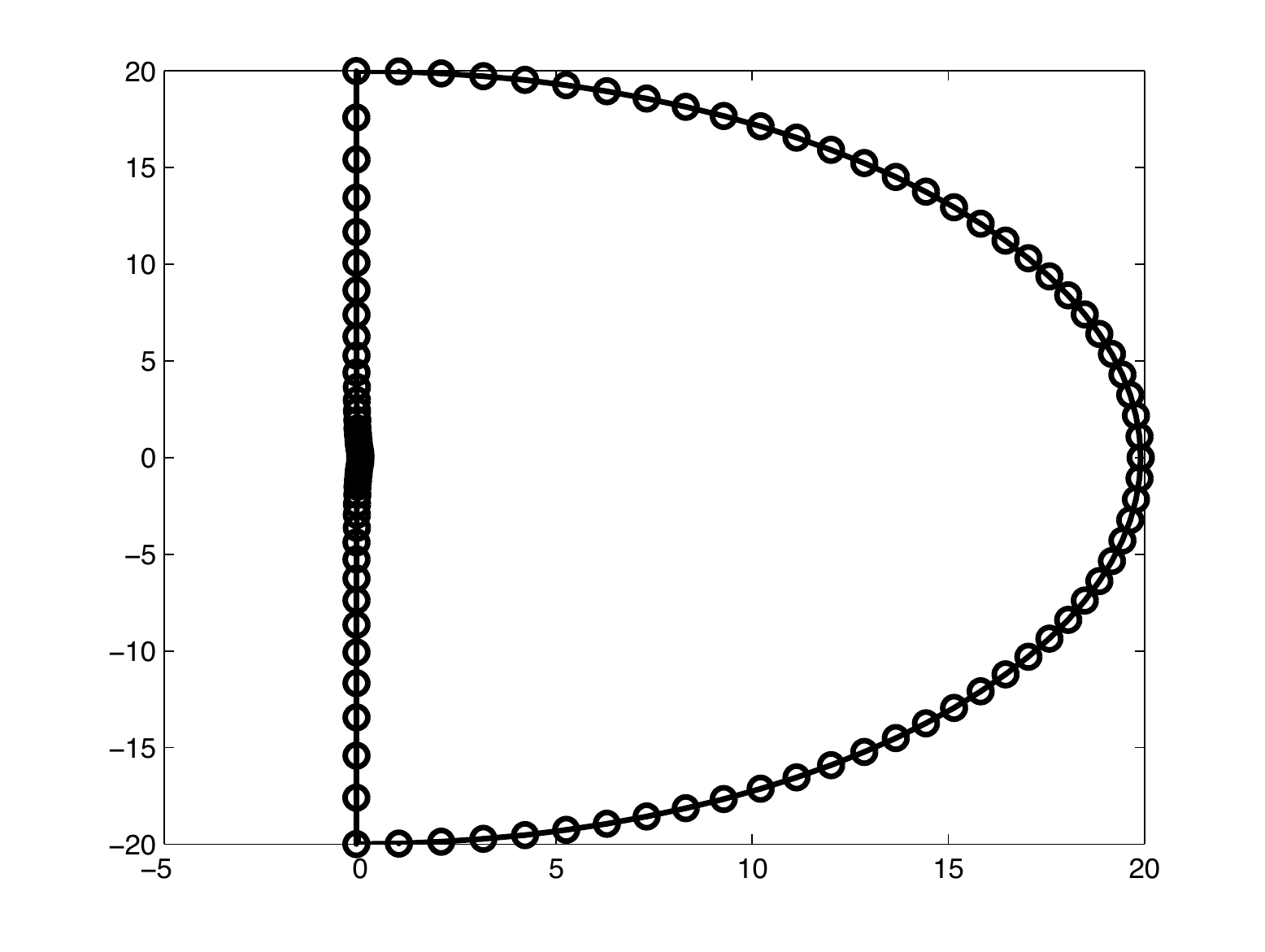}&
   		(b) \includegraphics[scale=.3]{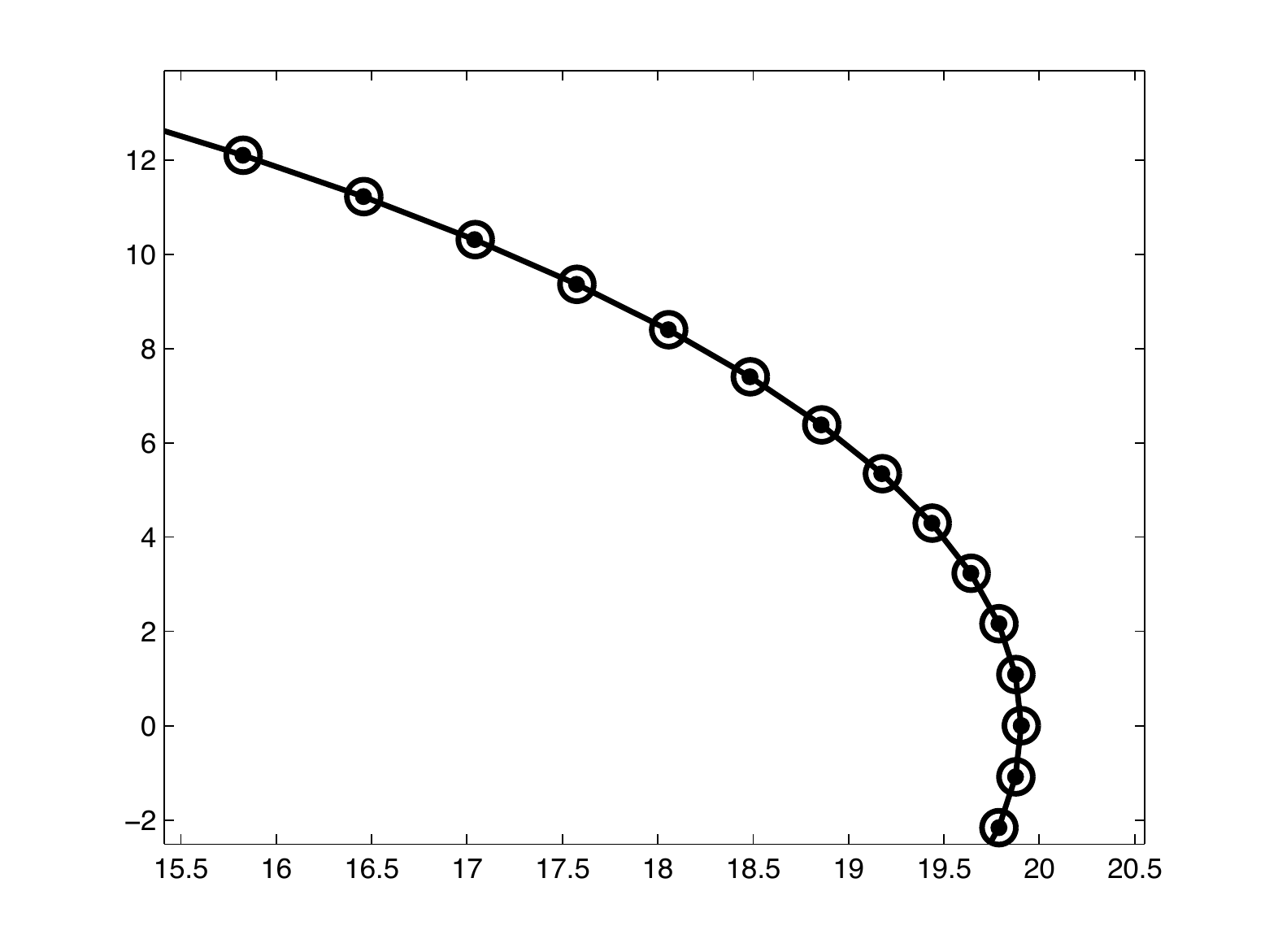}
    	\end{array}
$
\end{center}
\caption{
Comparison in the complex plane to the exact solution of \cite{JY} evaluated on a semicircle of
radius $R=10$. a) full size. b) zoom top right. Open circles correspond to the exact solution and the line with closed points to the $\mu(x)$ method. Here $\E=0$ so that $\phi\equiv 1$, $k=1$, and $q=0.3$.
}
   \label{comparison}
\end{figure}
\end{center}

For comparison, we test our code against an exact solution
found in \cite{JY} for $\phi\equiv 1$ of
\begin{equation}\label{5.5}
\begin{aligned}
    D_{ZND}(\lambda)&=
 \big( 2\lambda+(2-q-qk\Psi) \big) \left(\frac{\lambda}{k+\lambda}\right),
   \end{aligned}
\end{equation}
where 
\begin{equation}\label{5.4}
\Psi
:=\int_{-\infty}^{0}e^{-\int_{y}^{0}
\frac{\lambda}{\sqrt{1-2q(1-e^{ks})}}
ds}\frac{e^{(k+\lambda)y}}{\sqrt{1-2q(1-e^{ky})}}dy,
\qquad
P(\xi):=\frac{\lambda}{\sqrt{1-2q(1-e^{k\xi})}}.
\end{equation}
This formula results from a choice of $Z^-$
asymptotic to $e^{(\lambda + k)x}(*,1)^T$ as $x\to =\infty$
in agreement with \eqref{Sform}, hence should agree with the
results of our code.
As seen in Figure \ref{comparison}, the agreement of code with
exact formulat is excellent.

\subsection{Winding number computations}
To check stability, finally, we determine the number of roots
within a semicircle $S=\partial B(0,R)\cap \{\Re \lambda \ge 0\}$
of large radius $R$ by computing the winding number of the
image curve $D(S)$ as $S$ is traversed counterclockwise, which,
by the Principle of the Argument, is equal to the number of roots within.
We compute this winding number by varying values of $\lambda$ around 
$S$ along 
$20$ points 
of the contour, with mesh size taken quadratic in modulus to concentrate sample points near the origin where angles change more quickly,
   and summing the resulting changes in ${\rm arg}(D(\lambda))$, using $\Im \log D(\lambda) = {\rm arg} D(\lambda) ({\rm mod} 2\pi)$, available in \textsc{MATLAB} by direct function calls.
As a check on winding number accuracy, we test a posteriori that the change in argument of $D$ for each step is less than $0.2$, and add mesh points, as necessary to achieve this.  Recall, by Rouch\'e's Theorem, that accuracy is preserved so long as the argument varies by less than $\pi$ along each mesh interval.

Computations were carried out within the MATLAB-based STABLAB code
developed by J. Humpherys with help of the authors.
Using MATLAB's parallel computing toolbox on
an 8-core Power Macintosh workstation, we were able to 
achieve a speedup of over 600\%, similarly as in the previous numerical
studies \cite{BHZ,BLZ,BLeZ}.

\medskip
{\bf Determination of $R$.}
The radius $R$ is chosen so large that there exist no unstable roots
outside $B(0,R)$, ensuring that the roots found by our winding number
computation are the only possible ones on the nonstable half-plane
$\Re \lambda \ge 0$.
This could be done analytically as described in Section \ref{s:hf}
(and carried out, for example, in \cite{HLyZ}).
Here, following instead \cite{BLeZ,BHZ}, 
we use the asymptotics described in Proposition \ref{hfbehavior}
Remark \ref{hfrmk}, to match the reduced Lopatinski determinant
$D(\lambda)/\lambda$ to a first-order or higher-order approximant 
$e^{C\lambda}$ 
or $e^{C_1 \lambda +C_0 + C_{-1}\lambda^{-1})}$,\footnote{
For winding number computations, we find that it is
important to divide out as much behavior as we can, to avoid
excess winding.  Dividing by $\lambda$ is crucial, $C_1 e{C_2\lambda}$
still better.  Eventually, there is a break-even point in complexity
(and coefficient size of remainder) vs. power of $1/|\lambda|$,
here occurring at first order.
}
carrying out a convergence study to determine when $D(\lambda)/\lambda$
has sufficiently converged.
%

This was accomplished by, first, requiring that the relative error
between $D(\lambda)/\lambda$ and a best-fit value of $Ce^{C_1/\lambda+C_2/\lambda^2}$
be less than or equal to $0.2$, and, second,
that doubling the radius results in reduction of the error 
by a factor of approximately two, in accordance with the linear error
dependence predicted by analytical theory.

\begin{figure}[htbp]
\begin{center}
$
\begin{array}{lcr}
 (a) \includegraphics[scale=.25]{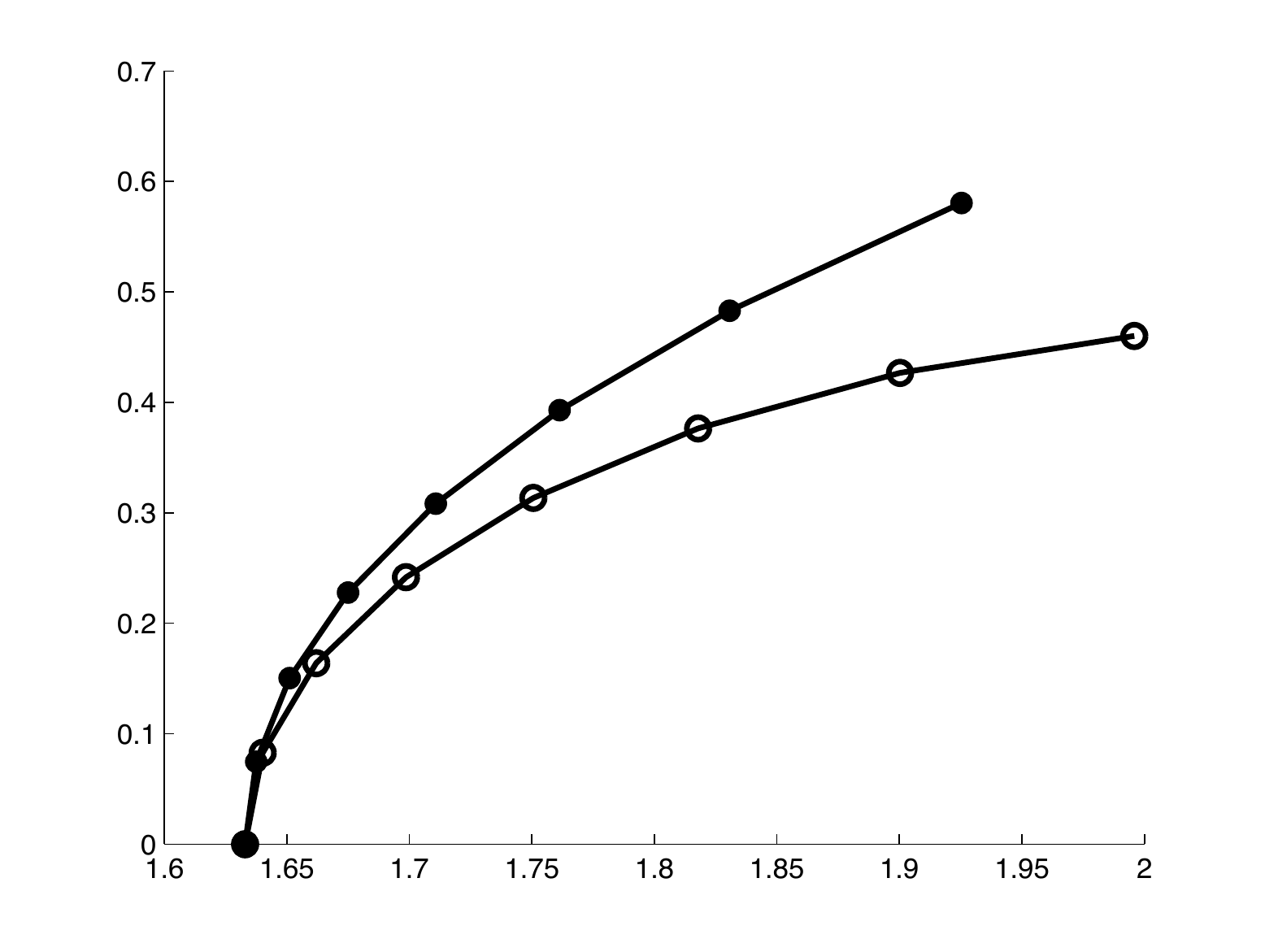}&(b)\includegraphics[scale=0.25]{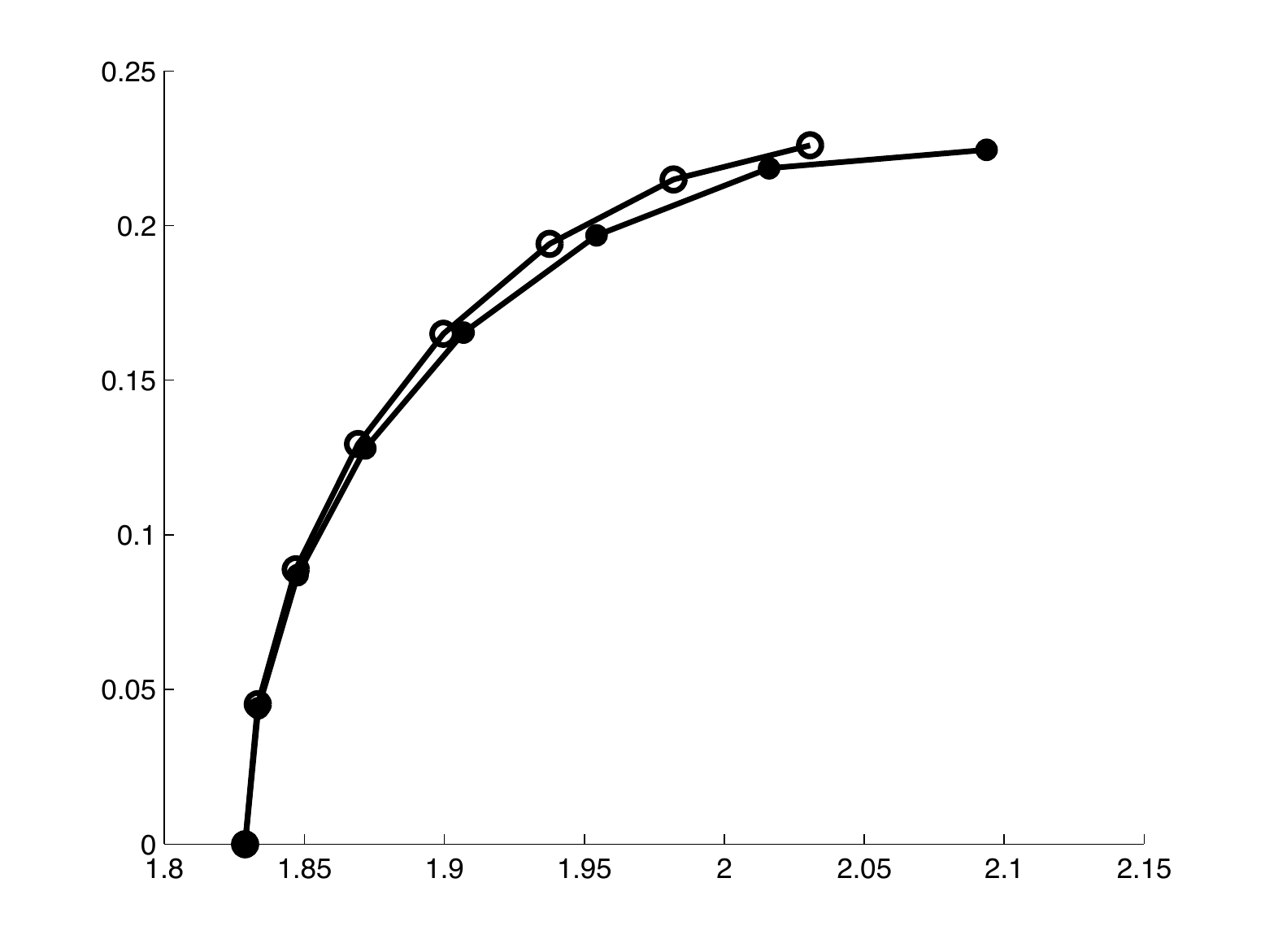}&(c)\includegraphics[scale=0.25]{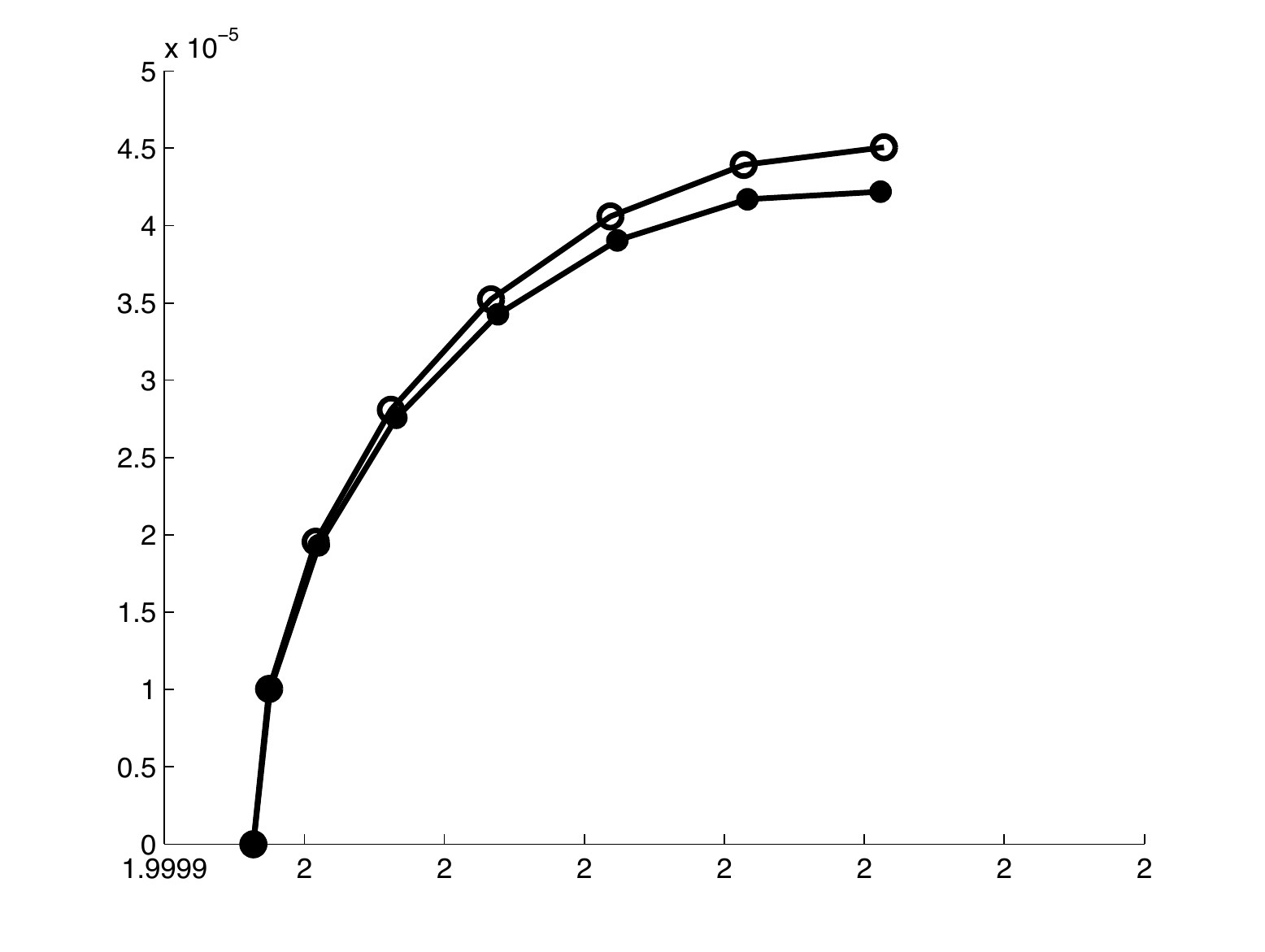}\\
 (d) \includegraphics[scale=.25]{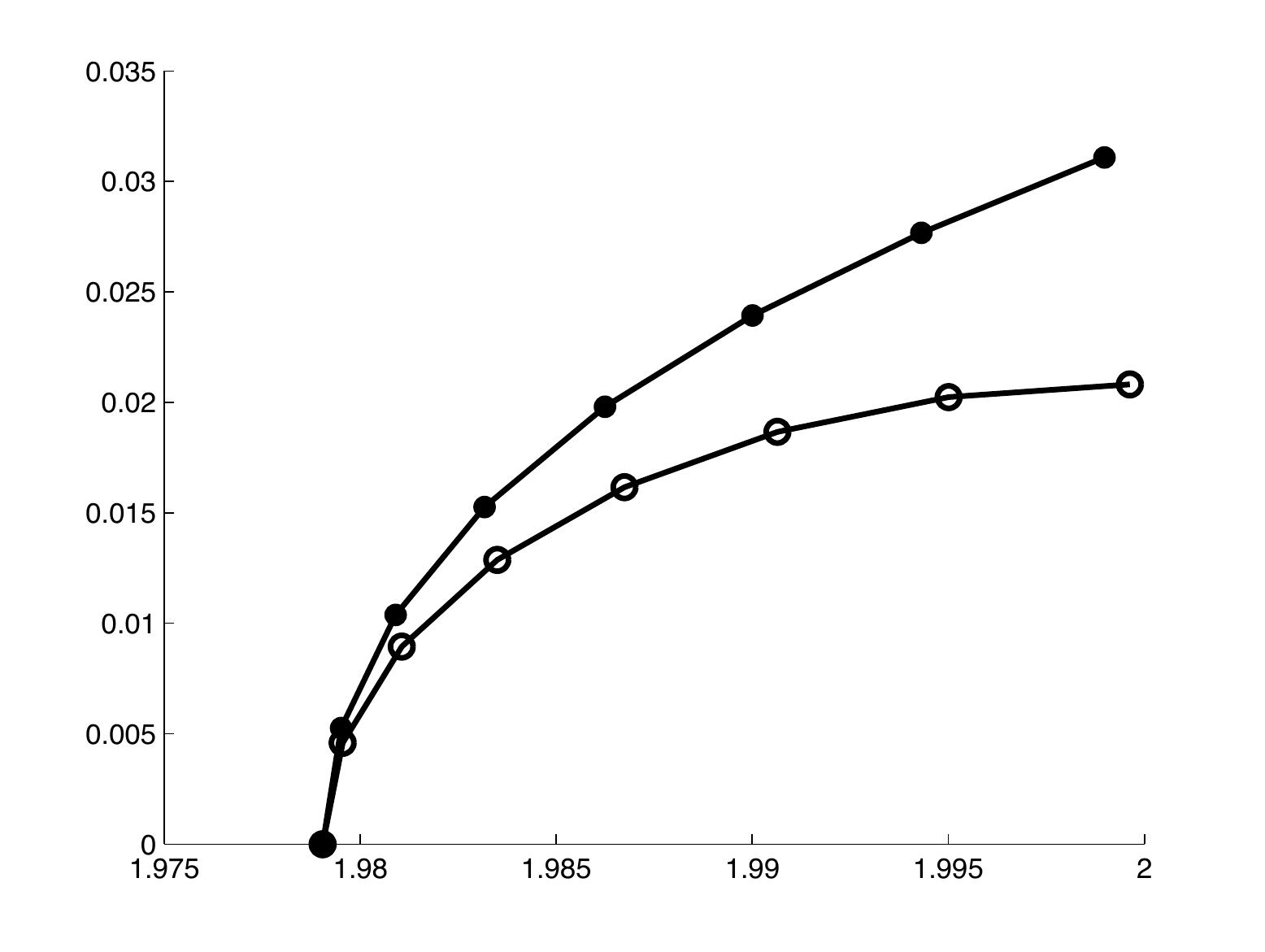}&(e)\includegraphics[scale=0.25]{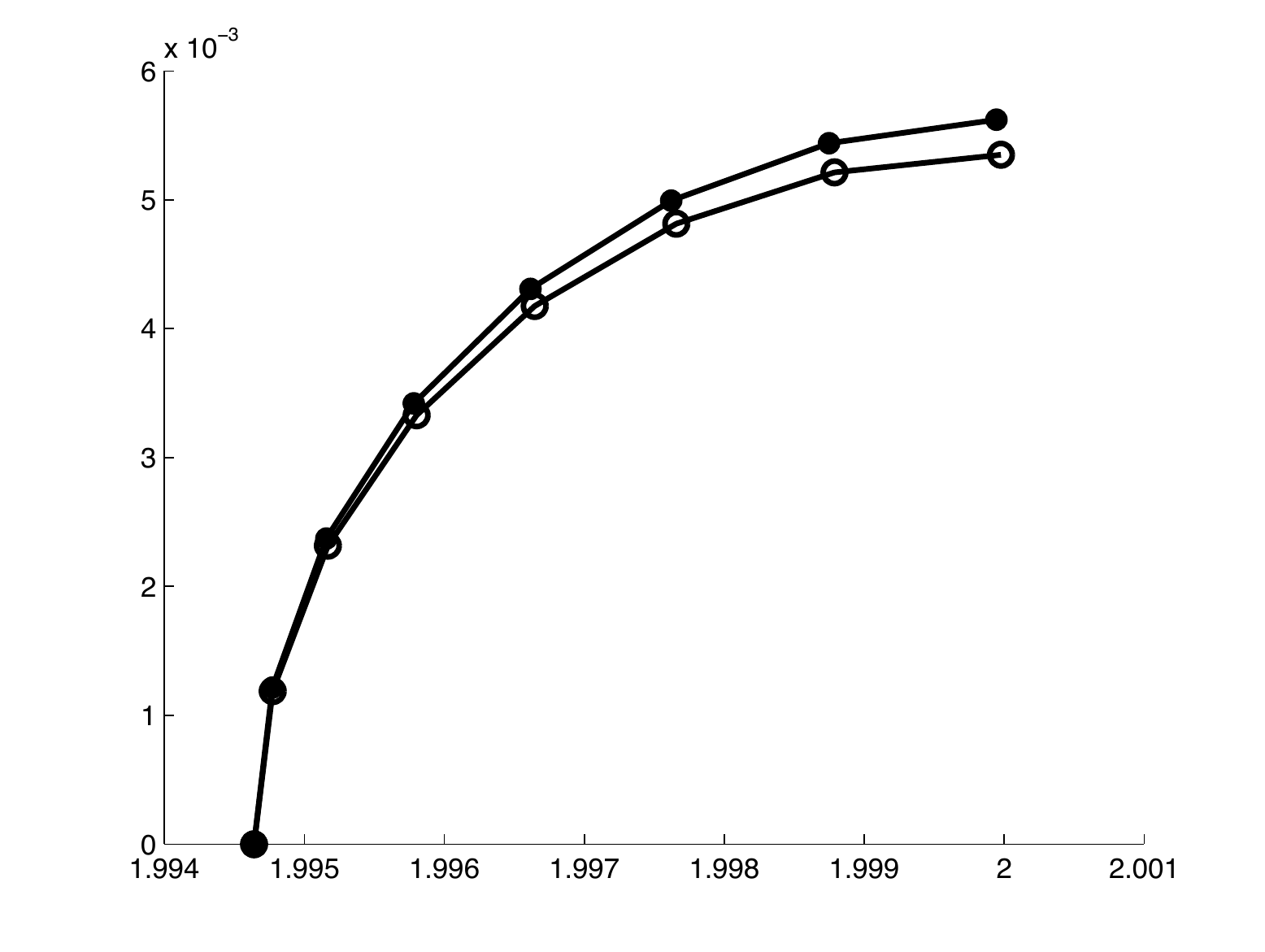}&(f)\includegraphics[scale=0.25]{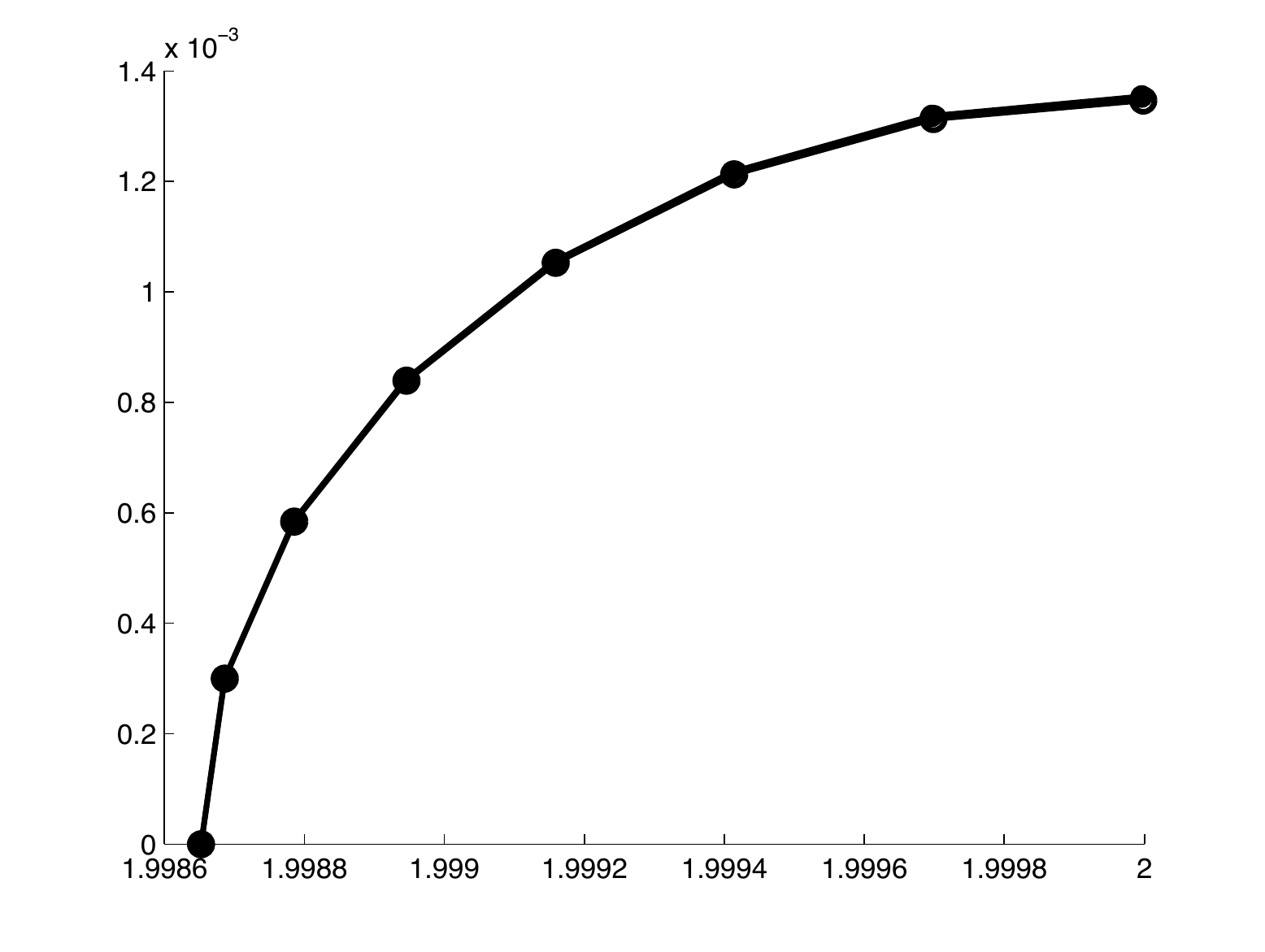}\\
 \end{array}
 $
\end{center}
\caption{
Here we use solid dots to plot the Lopatinski determinant evaluated on a quarter arc, and we use open circles to graph the approximating function, $C_1\exp(C2/\lambda)$. For the $\mu(x)$ method, we have Figures a-c corresponding respectively to $R=2$, $R=4$, and $R=2048$. For the polar method, we have Figures d-f corresponding respectively to $R=4$, $R=16$, $R=64$. Parameters are $\E=10$, $q=0.3$, $\phi(u)=Ce^{-\mathcal{E}/T(u)}$, $C=10^{21\mathcal{E}/40}$. 
}
\label{fig2}
\end{figure}

Typical convergence studies are illustrated in Figure \ref{fig2} and Tables \ref{bestfit} and \ref{bestfit2}.
A comparison of the Evans function vs. reduced Evans function is
given in Figure \ref{redfig}.

\begin{table}[!th]
\begin{center}
\begin{tiny}
\begin{tabular}{|c|c|c|c|c|}
\hline
Radius&Relative Error&K0&K1&K2\\ \hline
2&    0.00387783 & 0.693873 &   -0.157136   & 0.0202546 \\ \hline
4   & 0.000656158   & 0.693278 &   -0.152452 &  0.0110315  \\ \hline
8   & 8.88927e-05  &  0.693167    & -0.150713 &  0.00427183  \\ \hline
16 &   9.83384e-06 &    0.69315  & -0.150192 &   0.000241916  \\ \hline
32 &   6.99413e-07  &  0.693148  &   -0.15005 &   -0.00195688  \\ \hline
64  &  1.26154e-07 &  0.693147   & -0.150013 &   -0.003106  \\ \hline
128 &   5.78215e-08 & 0.693147  &  -0.150003  &  -0.00369357  \\ \hline
\end{tabular}
\end{tiny}
\end{center}
\captionsetup{font=scriptsize}
\caption{
For $\E=10$, $q=0.3$, $C=e^{\mathcal{E}/2}$, and $\phi(u)=e^{-\mathcal{E}/u}$ we examine the high frequency convergence of the adjoint $\mu(x)$ method with it's best fit with $e^{K0+K1/\lambda+K2/\lambda^2}$. 
\label{bestfit}
}
\end{table}

\begin{table}[!th]
\begin{center}
\begin{tiny}
\begin{tabular}{|c|c|c|c|c|}
\hline
Radius&Relative Error&K0&K1&K2\\ \hline
2  &  0.0436173      &   0.72503 &  -0.501318  &  0.0641958\\ \hline 
4   & 0.0294181   &  0.706568    & -0.347815  &  -0.254427 \\ \hline
8   & 0.0100462      &   0.696782   &  -0.190304 &   -0.88822 \\ \hline
16   & 0.00172202      &   0.693867  &   -0.0981737  &  -1.61602 \\ \hline
32    &0.000146306    &    0.693263   &  -0.0604492 &   -2.20532 \\ \hline
64   & 6.21728e-05      &  0.693164   &  -0.0480715 &  -2.58952 \\ \hline
128  &  2.54483e-05      &  0.693149   &  -0.0444984  &  -2.81055 \\ \hline
256   & 7.81771e-06      & 0.693147    & -0.0435363   & -2.92936 \\ \hline
512   & 2.14819e-06     &   0.693147   &  -0.0432867   & -2.99093 \\ \hline
1024   & 5.61212e-07 &      0.693147   &  -0.0432227 &   -3.02249 \\ \hline
\end{tabular}
\end{tiny}
\end{center}
\captionsetup{font=scriptsize}
\caption{
For $\E=10$, $q=0.3$, $\phi(u)=Ce^{-\mathcal{E}/T(u)}$, $T(u)=1-(1.5-u)^2$, $C=10^{21\mathcal{E}/40}$ we examine the high frequency convergence of the adjoint method scaled by $\mu(x)$ to it's best fit with $e^{K0+K1/\lambda+K2/\lambda^2}$. 
\label{bestfit2}
}
\end{table}

\medskip
{\bf Results.}
We computed, using the Evans--Lopatinski determinant (our basic algorithm)
rescaled by $\mu(x)$, a batch job for the parameter values,
$\{\mathcal{E},q\}=\{0.01:0.01:0.49\} \times \{ 0:0.1:5,\ 5.2:0.2:10,\ 12,\ 15,\ 20,\ 30,\ 40\}$, with the exception of a few numerically challenging parameters $\{E,q\}=\{20\}\times\{0.49\},\ \{30\}\times\{0.48,\ 0.49\},\ \{40\}\times\{0.47,\ 0.48,\ 0.49\}$, { with $C=10^{21\mathcal{E}/40}$, for ignition function $\phi(u)=Ce^{-\mathcal{E}/T(u)}$ where $T(u)=1-(u-1.5)^2$. Computational statistics are given in Table \ref{batchjobstat2}. For ignition function $\phi(u)=e^{\mathcal{E}/2}e^{-\mathcal{E}/u}$ we computed the Evans--Lopatinski determinant for 
$\{\mathcal{E},q\}=\{0:0.1:5,\ 5.2:0.2:10,\ 12,\ 15\}\times\{0.01:0.01:0.37,\ 0.375,\ 0.38:0.01:0.49\}\cup\{20\}\times\{0.01:0.01:0.37,\ 0.375,\ 0.38:0.01:0.47\}\cup \{25\}\times \{0.01:0.01:0.37,\ 0.375,\ 0.38:0.01:0.45\}\cup\{30\}\times\{0.01:0.01:0.37,\ 0.375,\ 0.38:0.01:0.40\}$. Computational statistics are given in Table \ref{batchjobstat1}.
All computations yielded winding number zero consistent with stability. 

\begin{table}[!th]
\begin{tabular}{|c|c|c|c|c|}
\hline 
q&E=0&E=20&E=30&E=40\\ 
\hline
0.01&(1e-4,4,10,0.0086,1) & (1e-4,40,10,0.018,1e1) & (1e-4,60,10,0.023,1e1) & (1e-4,80,10,0.026,2e1)\\ 
\hline
0.1&(1e-4,4,10,0.1,1) & (1e-4,40,21,0.2,3) & (1e-4,60,33,0.19,3) & (1e-4,80,40,0.2,4)\\ 
\hline
0.2&(1e-4,4,11,0.15,1) & (1e-4,40,39,0.19,3) & (1e-4,60,54,0.2,4) & (1e-4,80,74,0.2,7)\\ 
\hline
0.3&(1e-4,4,12,0.17,2) & (1e-4,40,49,0.19,3) & (1e-4,60,75,0.19,6) & (1e-4,256,89,0.19,2e1)\\ 
\hline
0.4&(1e-4,4,17,0.16,2) & (1e-4,40,62,0.2,4) & (1e-4,60,84,0.2,8) & (1e-4,512,107,0.2,4e1)\\ 
\hline
\end{tabular}
\caption{For $\phi=Ce^{-\mathcal{E}/T(u)}$ where $T(u)=1-(u-1.5)^2$ and $C=e^{21\mathcal{E}40}$ we record statistics for our Evans--Lopatinski determinant computations. The data represented are ($L$,$R$,$error$,$time$) where $L$ is the numerical value of infinity in the spatial $z$ domain, $R$ is the radius of the domain contour, $error$ is the maximum relative error between $\lambda$ contour output, and $time$ is the time it took to compute.}
\label{batchjobstat2}
\end{table}

\begin{table}[!th]
\begin{tabular}{|c|c|c|c|c|}
\hline 
&E=0&E=10&E=20&E=30\\ 
\hline
q=0.01&(0.1,4,0.004,0.4) & (0.1,20,0.0063,0.8) & (0.1,40,0.0075,1) & (0.1,60,0.0081,2)\\ 
\hline
q=0.1&(0.1,4,0.045,0.7) & (0.1,20,0.073,2) & (0.1,40,0.088,4) & (0.1,60,0.096,6)\\ 
\hline
q=0.2&(0.1,4,0.1,1) & (0.1,20,0.18,3) & (0.06,40,0.11,1e+01) & (0.06,60,0.13,2e+01)\\ 
\hline
q=0.3&(0.1,4,0.18,1) & (0.06,20,0.17,7) & (0.06,40,0.18,2e+01) & (0.03,60,0.18,7e+01)\\ 
\hline
q=0.4&(0.06,4,0.16,2) & (0.06,20,0.15,1e+01) & (0.03,40,0.17,8e+01) & (0.03,60,0.19,4e+02)\\ 
\hline
\end{tabular}
\caption{For $\phi=e^{\mathcal{E}/2}e^{-\mathcal{E}/u}$ we record statistics for our Evans--Lopatinski determinant computations. Here we use the polar adjoint method with ode solver tolerance set at $10^{-12}$. We set the radius of the domain contour to be the max of 2$\mathcal{E}$ and the radius obtained by the best curve fit of the Evans--Lopatinski determinant with  $e^{E_0+E_1/\lambda+E_2/\lambda^2}$. The data represented are ($L$,$R$,$error$,$time$) where $L$ is the numerical value of infinity in the spatial $z$ domain, $R$ is the radius of the domain contour, $error$ is the maximum relative error between $\lambda$ contour output, and $time$ is the time it took to compute.}
\label{batchjobstat1}
\end{table}



\section{Performance comparisons}
In the remainder of the paper, we collect a number of data comparing 
performance of the various methods.





In Figure \ref{fig1} we demonstrate agreement between methods that should yield the same output providing a verification of the correctness of our code. This also provides a visual example of which methods require a tighter $\lambda$ mesh in order to obtain relative error within tolerance.

\begin{figure}[htbp]
\begin{center}
$
\begin{array}{lcr}
 (a) \includegraphics[scale=.3]{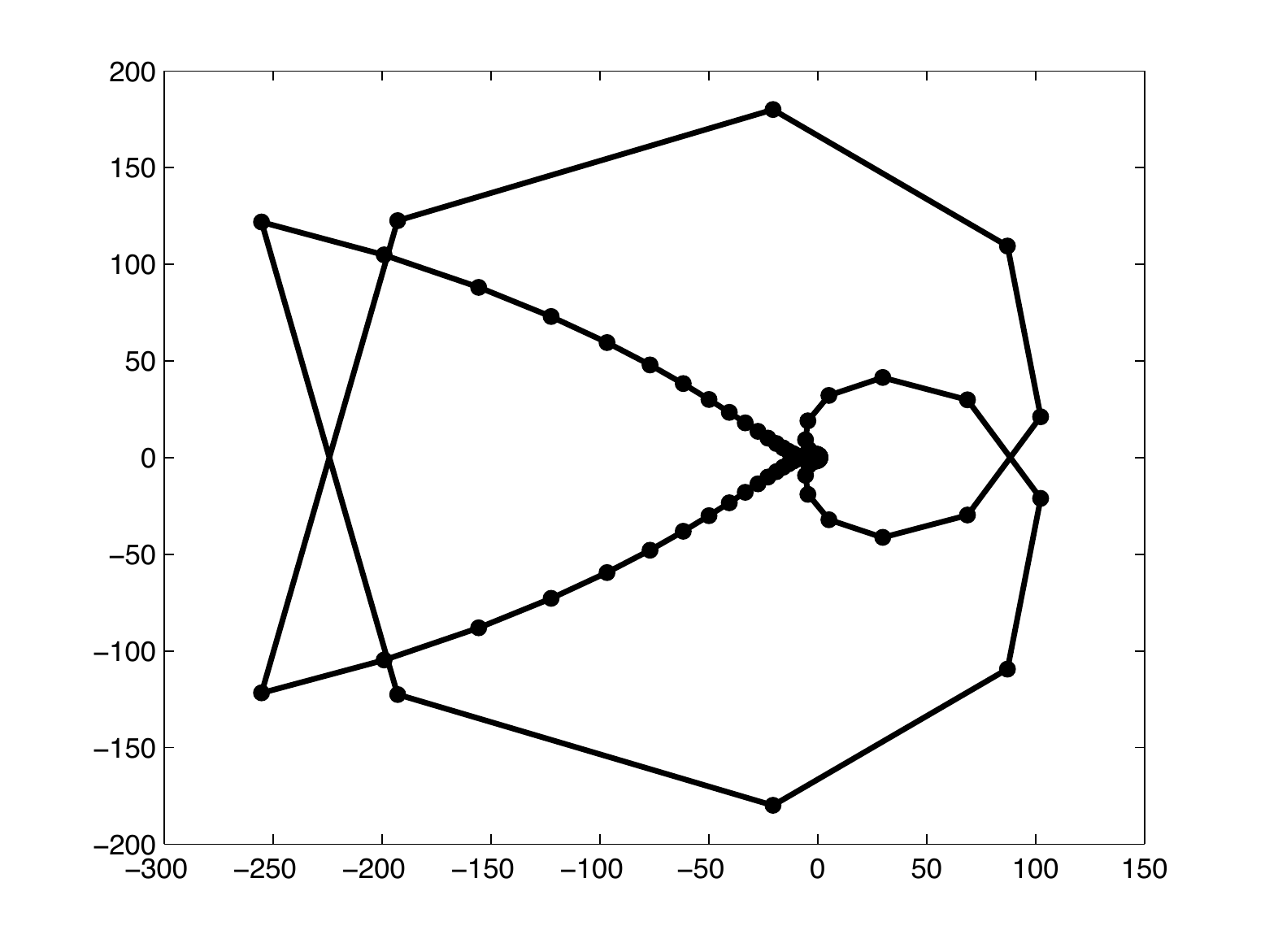}&(b)\includegraphics[scale=0.3]{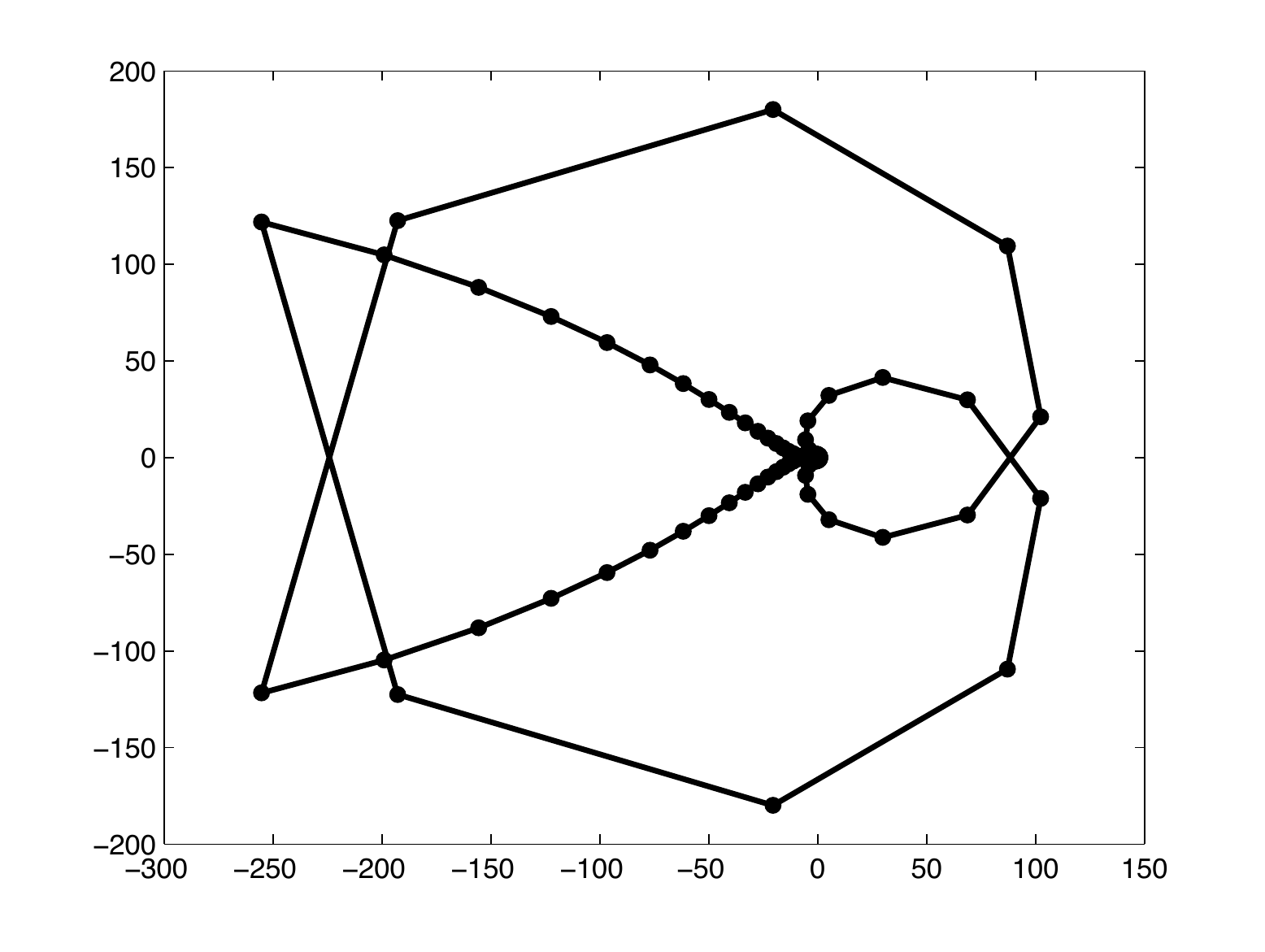}&(c)\includegraphics[scale=0.3]{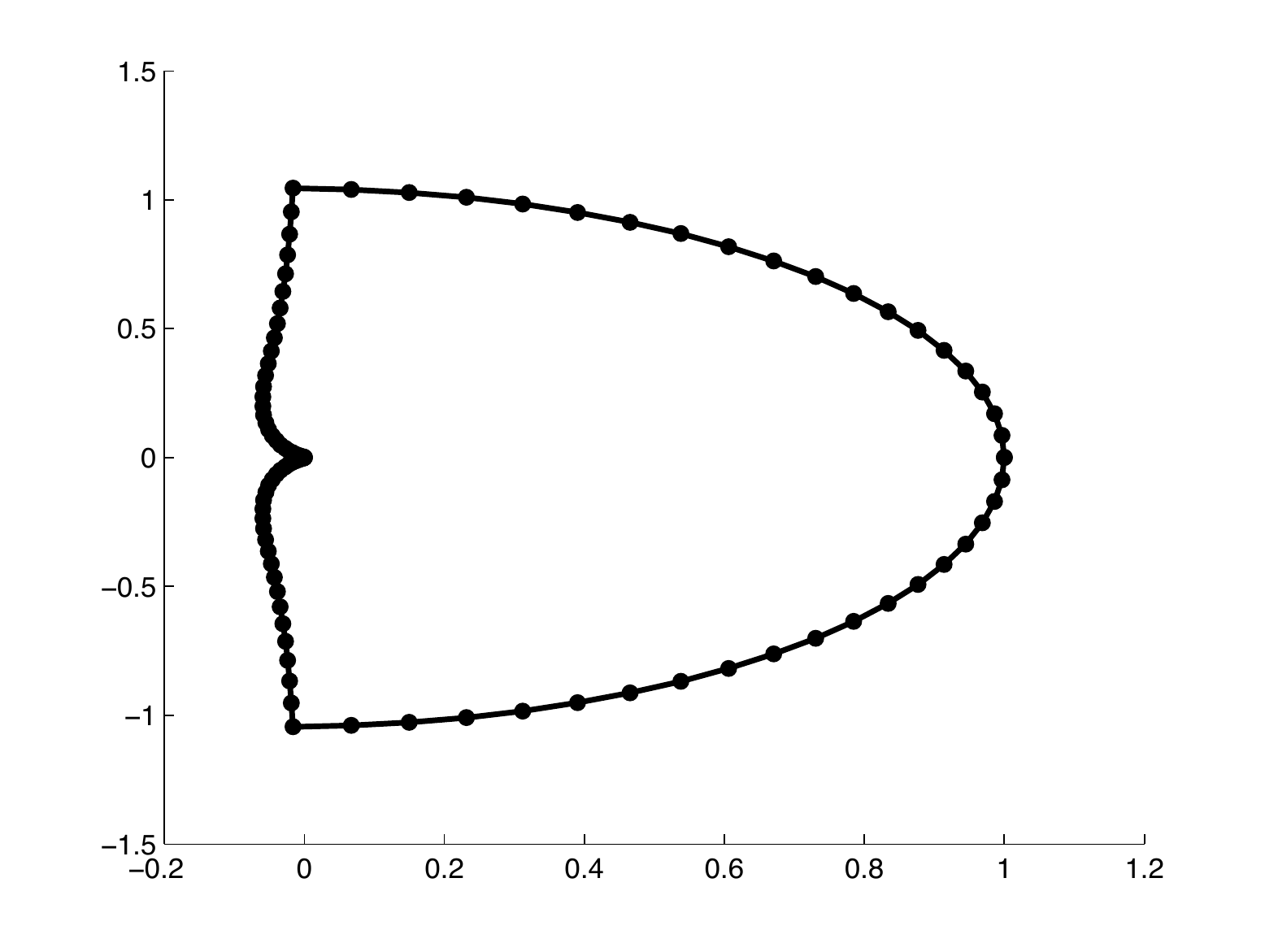}\\
 (d) \includegraphics[scale=.3]{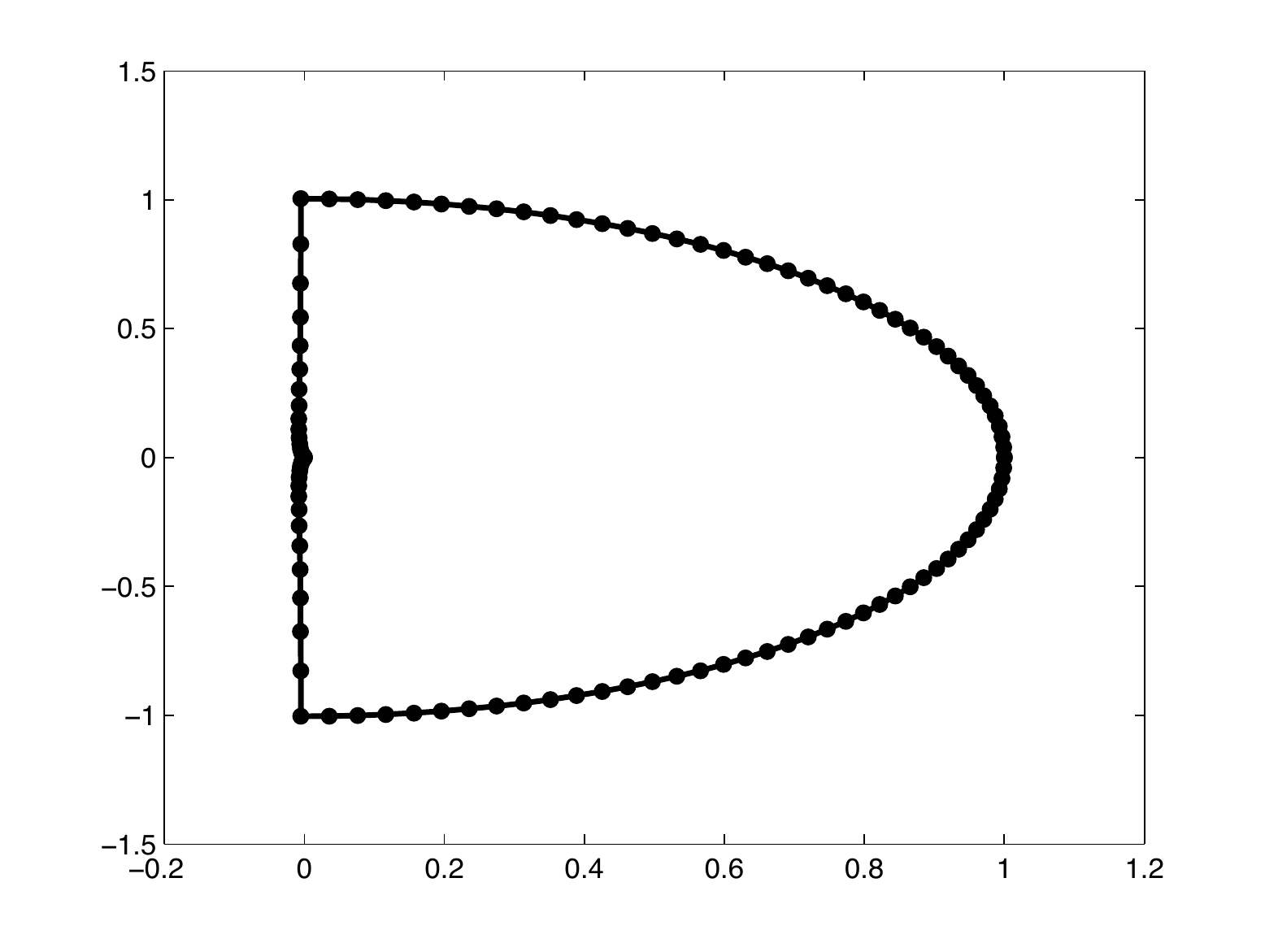}&(e)\includegraphics[scale=0.3]{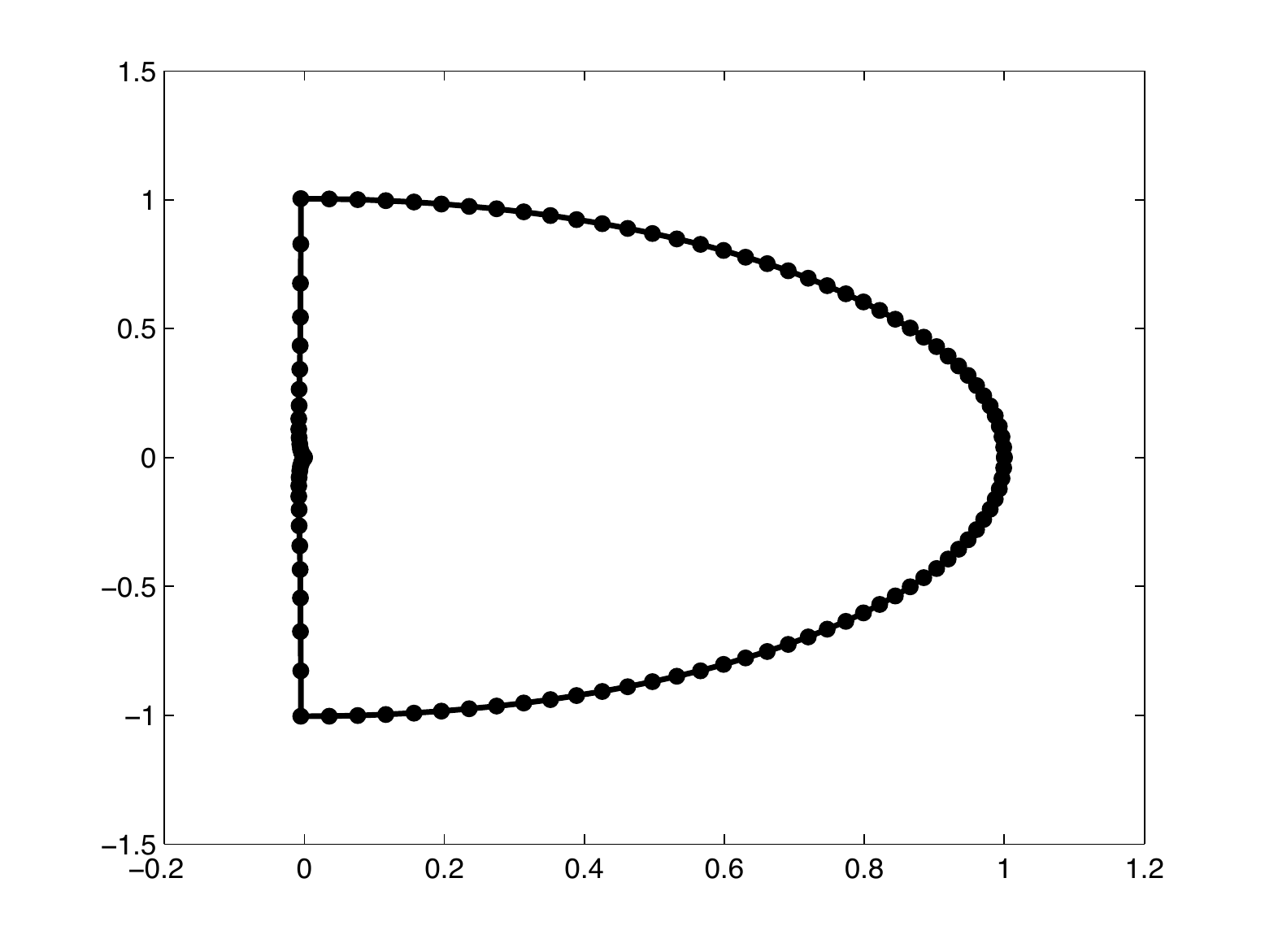}&(f)\includegraphics[scale=0.3]{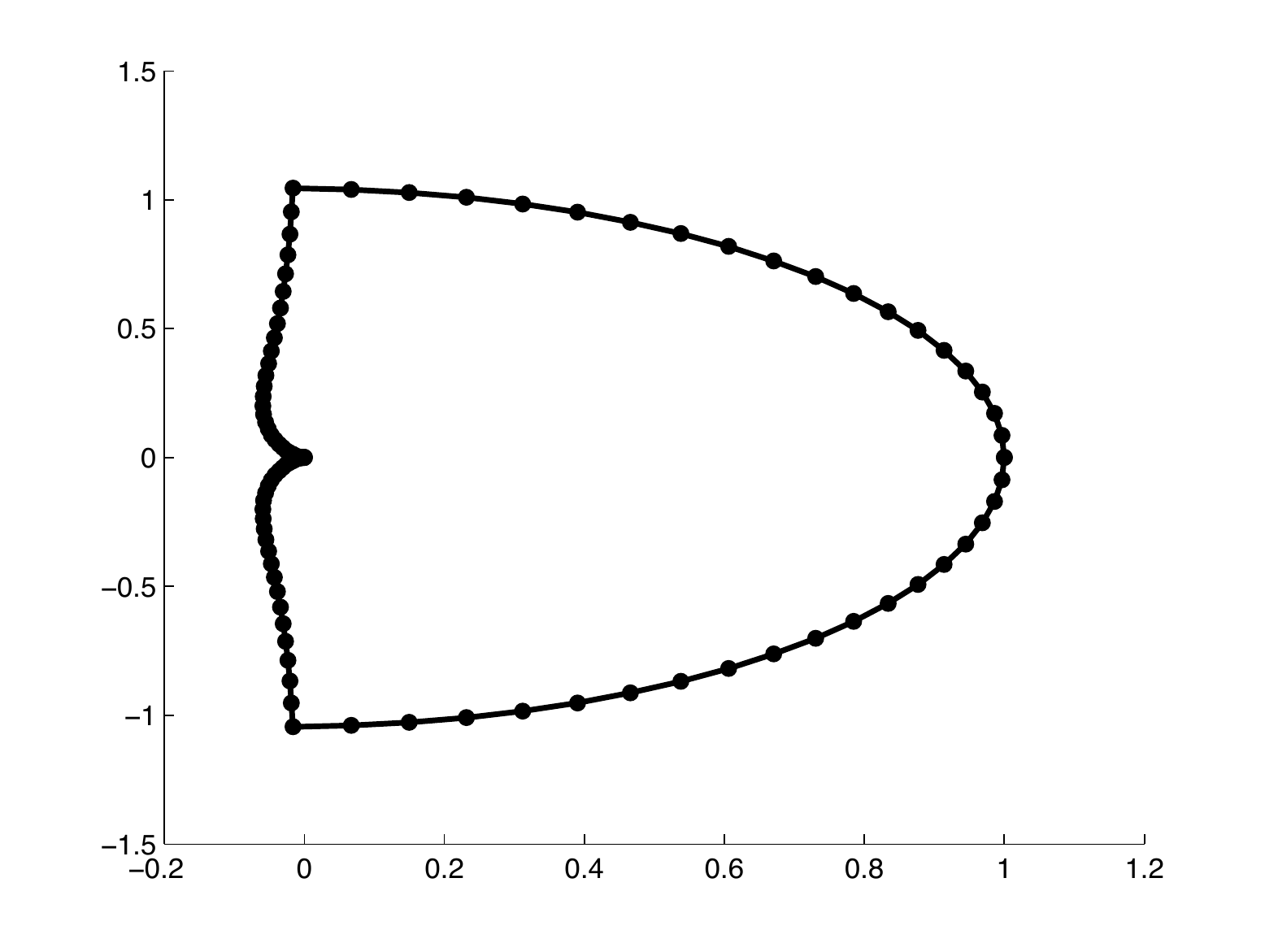}\\
 (g) \includegraphics[scale=.3]{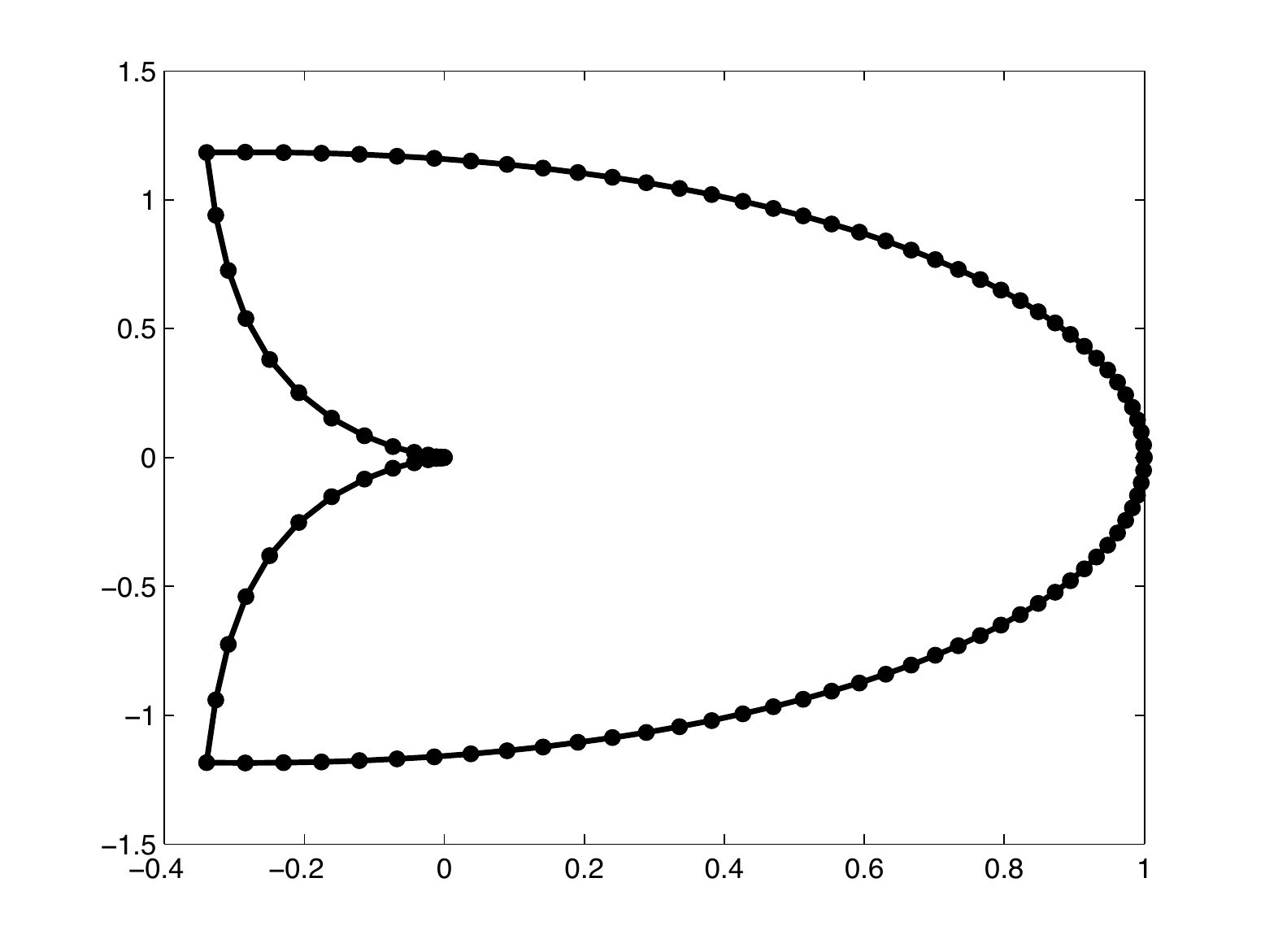}&(h)\includegraphics[scale=0.3]{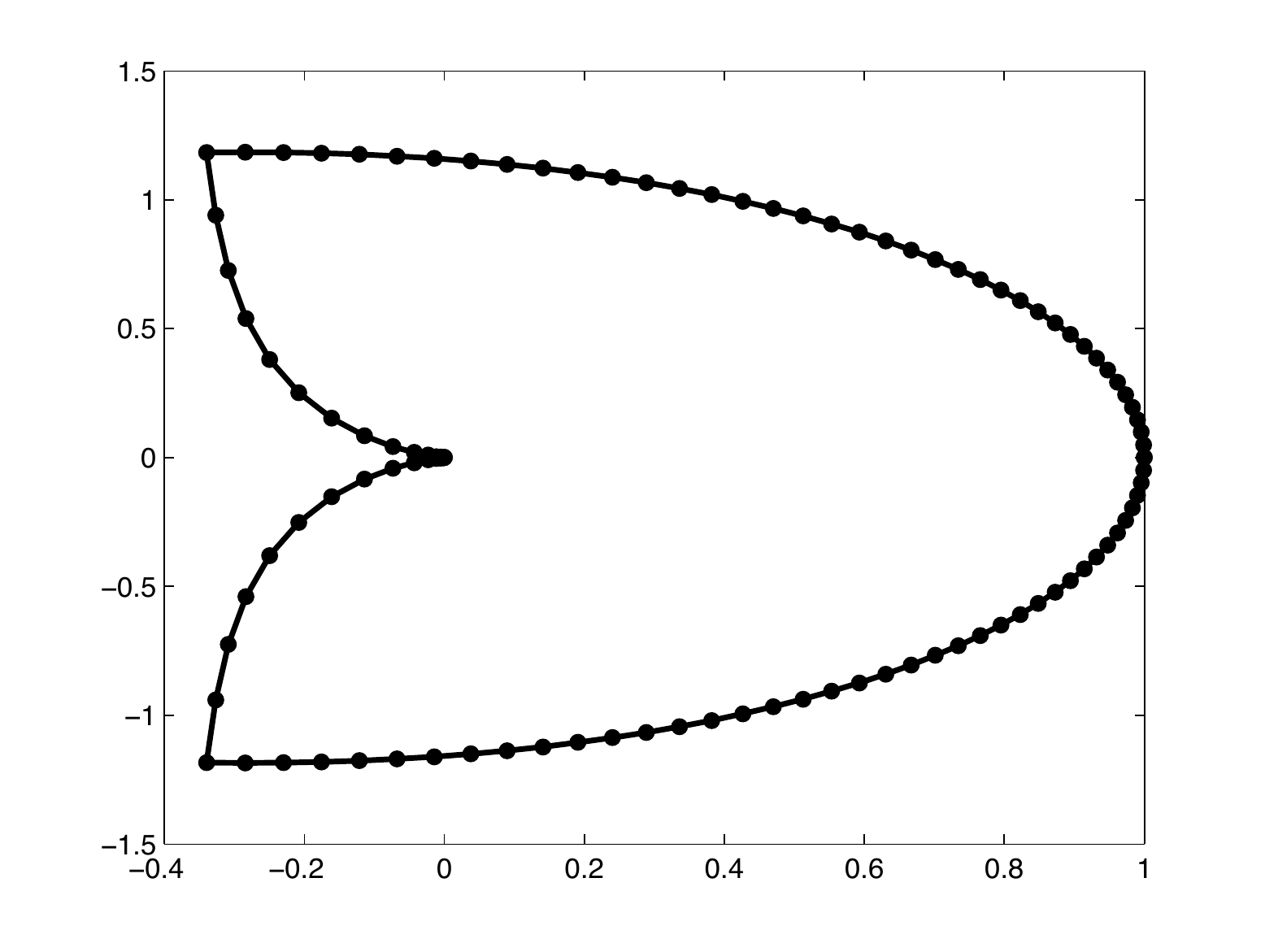}&(i)\includegraphics[scale=0.3]{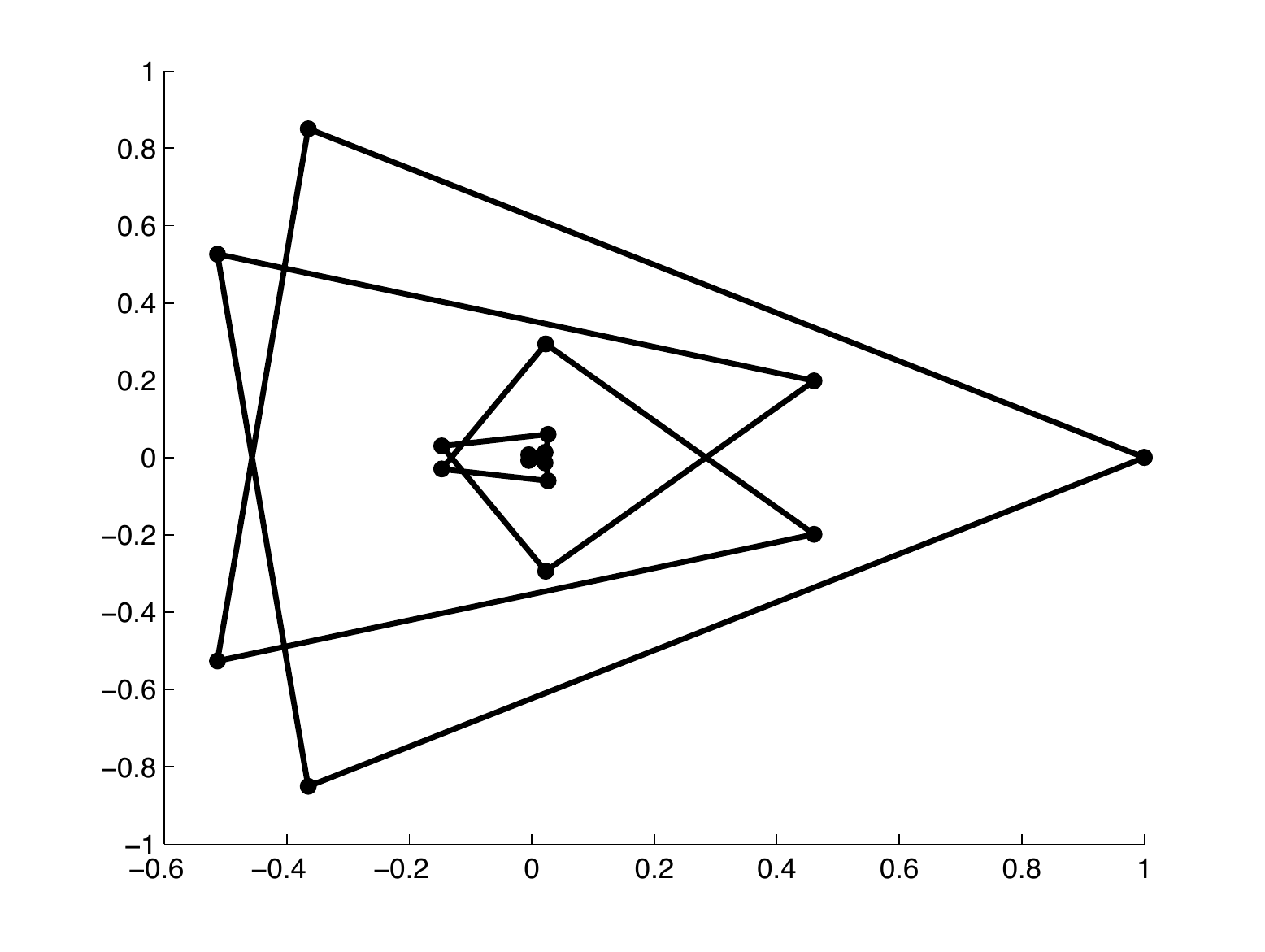}\\
  (j) \includegraphics[scale=0.3]{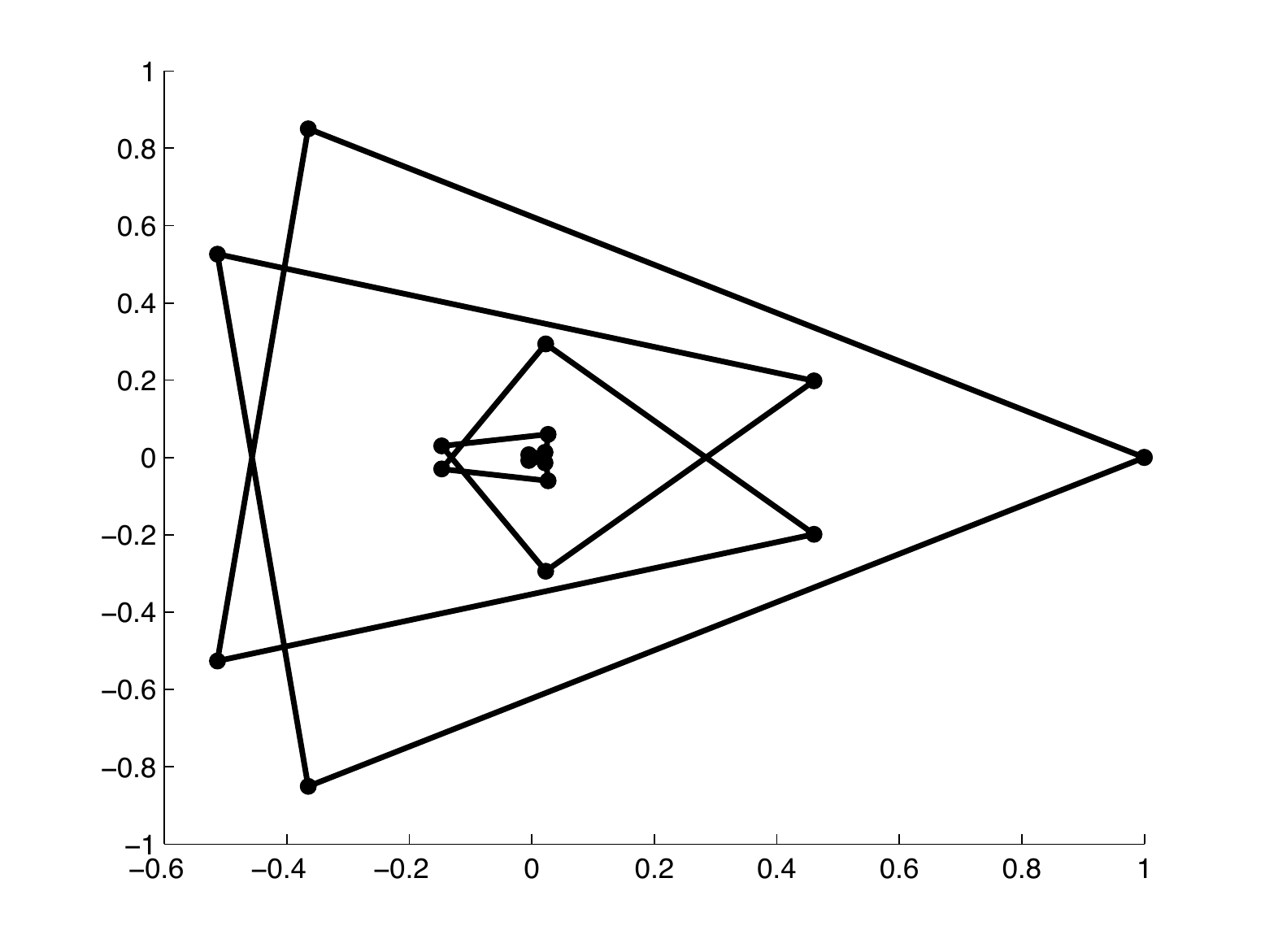}&(k)\includegraphics[scale=0.3]{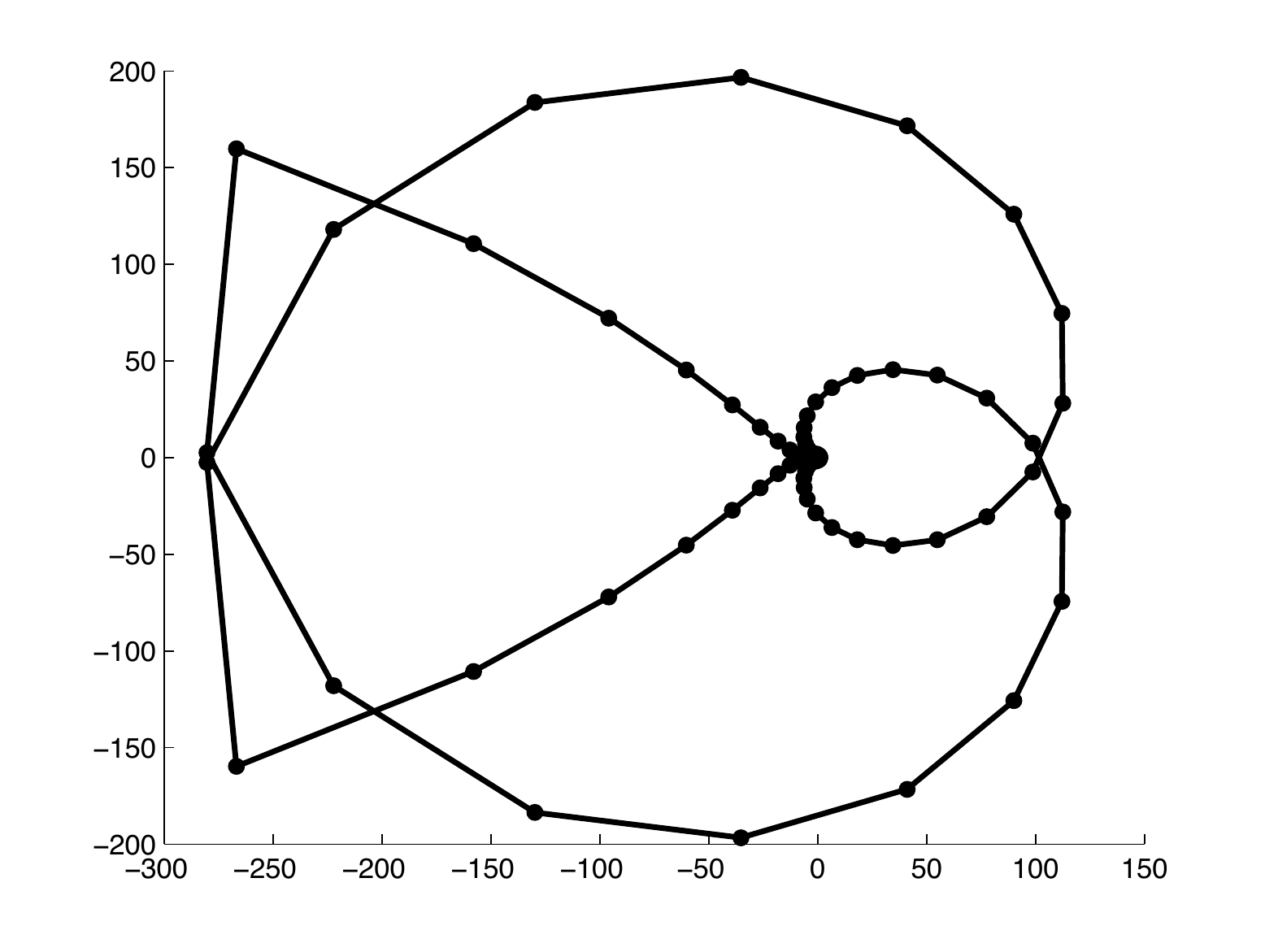}&\\
 \end{array}
 $
\end{center}
\caption{ We compare different methods of integration for $\E=10$, $q=0.3$, $\phi(u)=Ce^{-\mathcal{E}/T(u)}$, $T(u)=1-(1.5-u)^2$, and $C=10^{21\mathcal{E}/40}$, $R=10$, 
with 65 fixed domain mesh points. These comparisons serve as code verification, showing agreement between methods that should give the same output, and contrast the number of points needed by the various methods to obtain output with relative error within tolerance. (a) polar adjoint radial method, (b) adjoint $\mu$ method, (c) 
adjoint $\mu(x)$ method, (d) polar adjoint method, (e) polar method, (f) $\mu(x)$ method, (g) polar radial method, 
(h) $\mu$ method, (i) no $\mu$ method, (j) homogeneous Erpenbeck method (k) the homogeneous method of Lee and Stewart
}
\label{fig1}
\end{figure}

In Table \ref{comp1} we demonstrate computational time for the various methods for a fixed $\lambda$ mesh.
In this experiment, we simply set tolerance at $TOL=10^{-12}$ for
relative and absolute truncation error in the adaptive RK45 ODE solver
and integrated around a fixed radius $10$
semicircular contour in the complex $\lambda$-plane of the type used in our 
winding number computations for stability, with a small semicircle
of radius $10^{-4}$ removed around the origin.
The time taken to complete this computation gives a rough ``average''
measure of performance over different $\lambda$ regimes.
However, it is a bit conservative, as it measures truncation and
not convergence error, so does not reflect the expected better numerical
conditioning of the Humpherys--Zumbrun vs. other schemes.
For this reason, we have performed a series of more careful experiments
at individual $\lambda$-values, specifying convergence rather than
truncation error, described below.
Note finally that the standard fixed-mesh Lee-Stewart method described in
\cite{LS} is not included in the comparisons of Table
\ref{comp1}, since, not being adaptive, it does not allow specification
of truncation error in this simple way.
The important comparison to this scheme is 
carried out in the convergence error-based study below.

For the convergence error-based study,
 we evaluate the various methods at the representative values $\lambda=1,\ i,\ 10,\ 10i$ and determine the relative tolerance with which the methods must be solved in order for the output to converge to a specified tolerance. For the methods using an adaptive mesh in $z$, we specify absolute tolerance to be $10^{-15}$ and increase relative tolerance by powers of 10 starting at $10^{-1}$ until consecutive output is less than $10^{-6}$ in relative error. We use MATLAB's ode45 function which employs the fourth order Runge-Kutta-Fehlberg method. 
By specifying absolute tolerance to be $10^{-15}$, we ensure 
that relative rather than absolute tolerance 
determines the mesh setting.\footnote{
The ode45 algorithm requires that relative tolerance {\it or}
absolute tolerance be met.
}

For the fixed $z$ mesh method of Lee and Stewart, we determine a number of mesh points $N$ such that relative error between the fixed mesh inhomogeneous method and the adaptive inhomogeneous method, solved with error tolerance requirements of $10^{-10}$,  is less than $10^{-6}$ but greater than $9*10^{-7}$. For the fixed $z$ mesh method of Lee and Stewart, the reported tolerance is the relative error between the ODE output of these two methods. 

In Tables \ref{E0q0phi1}--\ref{E20q47phi2} we report the computational statistics for the $\mu(x)$ method, the adjoint $\mu(x)$ method, the hybrid method, the method of Lee and Stewart (adaptive mesh in $z$), the method of Lee and Stewart (fixed mesh in $z$), the polar adjoint method, and the polar adjoint radial method, for both ignition functions $\phi(u)=e^{-\mathcal{E}/u}$ and $\phi(u)=e^{-\mathcal{E}/T(u)}$ where $T(u)=1-(u-1.5)^2$. For the first ignition function we record data for $\mathcal{E}=1,\ 10,\ 20$, and $q=0.1,\ 0.3,\ 0.47$, and for the second ignition function we give statistics for $\mathcal{E}=0,\ 10,\ 40$ and $q=0,\ 0.3,\ 0.47$. Although not recorded here, we carried out studies for $q=0.2$ and $q=0.4$ and found the behavior quite similar to that reported.

The outcome, as described in the introduction, is that the Humpherys-Zumbrun
algorithm, implemented as $\mu(x)$, polar, or polar radial method, outperforms
an optimized, adaptive-mesh, version of the Lee-Stewart algorithm by factor
$1$-$10$, and outperforms the actual fixed-mesh version proposed in \cite{LS}
by factor $100$-$1000$, even in the high activation energy/square-wave limit.
The hybrid method did not perform significantly better
in the square-wave limit, and performed significantly worse in other regimes,
hence we recommend that this option be discarded.
Likewise, the Evans method, valid only when $\phi(2)<<\max \phi$,
does not seem to perform as well as the basic 
$\mu(x)$ method (see Table \ref{comp1}),
and so does not seem to be worth the trouble
of implementation in this special regime.
However, this option should perhaps not be discarded
without more systematic investigation than was carried out here.



\begin{table}[!th]
\begin{tiny}
\begin{tabular}{|c|c|c|c|c|}
\hline 
&E=0&E=20&E=30&E=40\\ 
\hline
q=0.01&(0.0001,4,10,0.0086,1) & (0.0001,40,10,0.018,1e+01) & (0.0001,60,10,0.023,1e+01) & (0.0001,80,10,0.026,2e+01)\\ 
\hline
q=0.1&(0.0001,4,10,0.1,1) & (0.0001,40,21,0.2,3) & (0.0001,60,33,0.19,3) & (0.0001,80,40,0.2,4)\\ 
\hline
q=0.2&(0.0001,4,11,0.15,1) & (0.0001,40,39,0.19,3) & (0.0001,60,54,0.2,4) & (0.0001,80,74,0.2,7)\\ 
\hline
q=0.3&(0.0001,4,12,0.17,2) & (0.0001,40,49,0.19,3) & (0.0001,60,75,0.19,6) & (0.0001,256,89,0.19,2e+01)\\ 
\hline
q=0.4&(0.0001,4,17,0.16,2) & (0.0001,40,62,0.2,4) & (0.0001,60,84,0.2,8) & (0.0001,512,107,0.2,4e+01)\\ 
\hline
\end{tabular}
\end{tiny}
\caption{Table of computational statistics, ($M$,$R$,points,relative error,time), where $-M$ is the value of numerical negative infinity, $R$ is the radius of the contour used to compute, points is how many points the adaptive solver used to compute the mesh, err is the maximum relative error between two points of the output, and time is how long it took to compute the contour. Here $\phi(u)=C e^{\mathcal{E}/T(u)}$ where $T(u)=1-(1.5-u)^2$ and $C$ is chosen so that $\bar z$ passes through $0.5$ at $x=-7$. Here we are using the $\mu(x)$ method. We take $R$ to be the maximum of $2 \mathcal{E}$ and the best root fit with $e^{E_0+E_1/\lambda +E_2/\lambda^2}$.   }
\label{comp0}
\end{table}

\begin{table}[!th]
\begin{center}
\begin{tabular}{|c|c|}
\hline
method of Lee and Stewart (adaptive mesh in x) & 37.0 \\ \hline
homogenous method of Lee and Stewart (adaptive mesh in x) &  35.8 \\ \hline
adjoint $\mu(x)$ method & 14.8 \\ \hline
hybrid method &  27.6 \\ \hline
polar adjoint method &  13.5 \\ \hline
inhomogeneous Erpenbeck method &  40.1 \\ \hline
polar adjoint radial method& 20.8 \\ \hline
adjoint $\mu$ method &  18.2 \\ \hline
polar/Drury method &  13.5 \\ \hline
$\mu(x)$ method &  14.5 \\ \hline
 polar radial method & 20.9 \\ \hline 
 $\mu$ method & 17.9 \\ \hline
 Evans function method & 21.2 \\ \hline
homogenous Erpenbeck method & 32.4 \\ \hline
centered inhomogeneous Erpenbeck method & 51.8 \\ \hline
\end{tabular}
\end{center}
\caption{Computational time comparison on a fixed $\lambda$ mesh of 55 points with radius 10 and inner radius $10^{-4}$ for various methods with $E=10$, $q=0.3$, $\phi(u)=e^{-\mathcal{E}/T(u)}$, $T(u)=1-(1.5-u)^2$, $C=10^{21\mathcal{E}/40}$ and relative and absolute tolerance is set at $10^{-12}$ in the ODE solver.}
\label{comp1}
\end{table}

\begin{table}[!th]\begin{tiny}
\begin{tabular}{|c||c|c|c||c|c|c||c|c|c||c|c|c|}
\hline
 &$\lambda=1$& & &$\lambda=i$& & &$\lambda=10$& & &$\lambda=10i$ & &\\
\hline
  & Tol & Time & pnts& Tol & Time & pnts & Tol & Time & pnts & Tol & Time & pnts\\
 \hline
A&-&-&-&-&-&-&-&-&-&-&-&\\ \hline
B&  1.0e-2 & 0.0027 &   2    &  1.0e-2 & 0.0027 &   2    &  1.0e-2 & 0.0026 &   2    &  1.0e-2 & 0.0026 &   2    \\ \hline
C&  1.0e-2 & 0.0077 &   6    &  1.0e-2 & 0.0065 &   5    &  1.0e-2 & 0.014 &   15    &  1.0e-2 & 0.015 &   14    \\ \hline
D&  1.0e-2 & 0.039 &   37    &  1.0e-7 & 0.093 &   101    &  1.0e-8 & 0.88 &   1098    &  1.0e-8 & 1.1 &   1182    \\ \hline    
E&  6.9e-2 & 0.025 &   2    &  1.1e-7 & 26 &   6.55e4 &  8.5e-7 & 3.7e+02 &   1048580&  7.6e-7 & 420 &   1048580\\ \hline
F&  1.0e-02 & 0.0027 &   2    &  1.0e-02 & 0.0028 &   2    &  1.0e-02 & 0.0026 &   2    &  1.0e-02 & 0.0026 &   2    \\ \hline 
G& 1.0e-02& 0.0041&  3 & 1.0e-0.2& 0.0042& 3& 1.0e-02&0.0045&3&1.0e-02&0.0043&3\\ 
\hline
\end{tabular}
\end{tiny}
\captionsetup{font=scriptsize}
\caption{
Efficiency comparison for several methods: A- the $\mu(x)$ method, B- the adjoint $\mu(x)$ method, C- the hybrid method, D- the method of Lee and Stewart (using an adaptive mesh in $z$), E- the method of Lee and Stewart using a 4th order fixed mesh Runge Kutta solver, F- the polar adjoint method, and G- the polar adjoint radial method. Here $\mathcal{E}=0$, $q=0$, $C=10^{21\mathcal{E}/40}$, and $\phi(u)=Ce^{-\mathcal{E}/T(u)}$, $T(u)=1-(u-1.5)^2$.
}
\label{E0q0phi1}
\end{table}
 
%

\clearpage

\begin{table}[!th]\begin{tiny}
\begin{tabular}{|c||c|c|c||c|c|c||c|c|c||c|c|c|}
\hline
 &$\lambda=1$& & &$\lambda=i$& & &$\lambda=10$& & &$\lambda=10i$ & &\\
\hline
  & Tol & Time & pnts& Tol & Time & Mesh & Tol & Time & Mesh & Tol & Time & Mesh\\
 \hline
A&  1.0e-02 & 0.017 &   12    &  1.0e-02 & 0.017 &   11    &  1.0e-02 & 0.067 &   69    &  1.0e-02 & 0.096 &   93    \\ \hline 
B &  1.0e-04 & 0.021 &   16    &  1.0e-04 & 0.023 &   15    &  1.0e-03 & 0.07 &   73    &  1.0e-05 & 0.1 &   104    \\ \hline 
C &  1.0e-02 & 0.02 &   14    &  1.0e-03 & 0.023 &   15    &  1.0e-02 & 0.075 &   81    &  1.0e-04 & 0.12 &   121    \\ \hline 
D&  1.0e-08 & 0.16 &   191    &  1.0e-08 & 0.16 &   177    &  1.0e-09 & 2.1 &   2603    &  1.0e-09 & 2.3 &   2634    \\ \hline 
E&  5.7e-07 & 47 &   131072    &  6.1e-07 & 53 &   131072&&&&&&\\ \hline
F&  1.0e-04 & 0.021 &   16    &  1.0e-05 & 0.025 &   19    &  1.0e-05 & 0.083 &   91    &  1.0e-04 & 0.099 &   99    \\ \hline
G &1.0e-04& 0.0240&24&1.0e-04&0.0271&24& 1.0e-05&0.0931&97& 1.0e-05&0.1180&109 \\ 
\hline
\end{tabular}
\end{tiny}
\captionsetup{font=scriptsize}
\caption{
Efficiency comparison for several methods: A- the $\mu(x)$ method, B- the adjoint $\mu(x)$ method, C- the hybrid method, D- the method of Lee and Stewart (using an adaptive mesh in $z$), E- the method of Lee and Stewart using a 4th order fixed mesh Runge Kutta solver, F- the polar adjoint method, and G- the polar adjoint radial method. Here $\mathcal{E}=0$, $q=0.3$, $C=10^{21\mathcal{E}/40}$, and $\phi(u)=Ce^{-\mathcal{E}/T(u)}$, $T(u)=1-(u-1.5)^2$.
}
\label{E0q3phi1} 
\end{table}

%

\begin{table}[!th]\begin{tiny}
\begin{tabular}{|c||c|c|c||c|c|c||c|c|c||c|c|c|}
\hline
 &$\lambda=1$& & &$\lambda=i$& & &$\lambda=10$& & &$\lambda=10i$ & &\\
\hline
  & Tol & Time & pnts& Tol & Time & Mesh & Tol & Time & Mesh & Tol & Time & Mesh\\
 \hline
A&  1.0e-02 & 0.022 &   19    &  1.0e-02 & 0.024 &   18    &  1.0e-02 & 0.1 &   119    &  1.0e-02 & 0.17 &   176    \\ \hline 
B&  1.0e-04 & 0.029 &   28    &  1.0e-04 & 0.032 &   28    &  1.0e-03 & 0.11 &   126    &  1.0e-05 & 0.2 &   212    \\ \hline 
C &  1.0e-03 & 0.028 &   25    &  1.0e-04 & 0.035 &   30    &  1.0e-02 & 0.12 &   131    &  1.0e-05 & 0.23 &   249    \\ \hline 
D&  1.0e-08 & 0.31 &   379    &  1.0e-08 & 0.32 &   367    &  1.0e-09 & 4.5 &   5641    &  1.0e-09 & 5 &   5701    \\ \hline 
E&4.6e-7&188&524288&4.7e-7&211&524288&-&-&$>10^6$&-&-&\\ \hline
F&  1.0e-04 & 0.027 &   28    &  1.0e-05 & 0.04 &   38    &  1.0e-04 & 0.16 &   185    &  1.0e-05 & 0.2 &   212    \\ \hline
G& 1.0e-05&0.0441&45& 1.0e-06&0.0705&63& 1.0e-06&0.2318&241& 1.0e-05&0.2381&220\\
\hline
\end{tabular}
\end{tiny}
\captionsetup{font=scriptsize}
\caption{
Efficiency comparison for several methods: A- the $\mu(x)$ method, B- the adjoint $\mu(x)$ method, C- the hybrid method, D- the method of Lee and Stewart (using an adaptive mesh in $z$), E- the method of Lee and Stewart using a 4th order fixed mesh Runge Kutta solver, F- the polar adjoint method, and G- the polar adjoint radial method. Here $\mathcal{E}=0$, $q=0.47$, $C=10^{21\mathcal{E}/40}$, and $\phi(u)=Ce^{-\mathcal{E}/T(u)}$, $T(u)=1-(u-1.5)^2$.
\label{E0q47phi1} 
}
\end{table}

\begin{table}[!th]\begin{tiny}
\begin{tabular}{|c||c|c|c||c|c|c||c|c|c||c|c|c|}
\hline
 &$\lambda=1$& & &$\lambda=i$& & &$\lambda=10$& & &$\lambda=10i$ & &\\
\hline
  & Tol & Time & pnts& Tol & Time & Mesh & Tol & Time & Mesh & Tol & Time & Mesh\\
 \hline
A& -&-&-&-&-&-&-&-&-&-&-&-\\ \hline
B&  1.0e-2 & 0.0027 &   2    &  1.0e-2 & 0.0028 &   2    &  1.0e-2 & 0.0026 &   2    &  1.0e-2 & 0.0028 &   2    \\ \hline 
C&  1.0e-2 & 0.01 &   10    &  1.0e-2 & 0.01 &   9    &  1.0e-2 & 0.021 &   23    &  1.0e-2 & 0.037 &   39    \\ \hline 
D &  1.0e-8 & 0.39 &   485    &  1.0e-8 & 0.39 &   438    &  1.0e-9 & 4.9 &   6134    &  1.0e-9 & 6.1 &   6816    \\ \hline 
E&2.7e-6&187&524288&2.7e-6&210&524288&-&-&$>$1e+6&-&-&$>$1e+6\\ \hline
F&  1.0e-02 & 0.0026 &   2    &  1.0e-02 & 0.0026 &   2    &  1.0e-02 & 0.0025 &   2    &  1.0e-02 & 0.0026 &   2    \\ \hline 
G&1.0e-02&0.0040&3&1.0e-02&0.0041&3&1.0e-02&0.0045&3&1.0e-02&0.0050&3\\
\hline
\end{tabular}
\end{tiny}
\captionsetup{font=scriptsize}
\caption{
Efficiency comparison for several methods: A- the $\mu(x)$ method, B- the adjoint $\mu(x)$ method, C- the hybrid method, D- the method of Lee and Stewart (using an adaptive mesh in $z$), E- the method of Lee and Stewart using a 4th order fixed mesh Runge Kutta solver, F- the polar adjoint method, and G- the polar adjoint radial method. Here $\mathcal{E}=10$, $q=0$, $C=10^{21\mathcal{E}/40}$, and $\phi(u)=Ce^{-\mathcal{E}/T(u)}$, $T(u)=1-(u-1.5)^2$.
\label{E10q0phi1} 
}
\end{table}


\begin{table}[!th]\begin{tiny}
\begin{tabular}{|c||c|c|c||c|c|c||c|c|c||c|c|c|}
\hline
 &$\lambda=1$& & &$\lambda=i$& & &$\lambda=10$& & &$\lambda=10i$ & &\\
\hline
  & Tol & Time & pnts& Tol & Time & Mesh & Tol & Time & Mesh & Tol & Time & Mesh\\
 \hline
A&  1.0e-6 & 0.034 &   29    &  1.0e-6 & 0.036 &   29    &  1.0e-3 & 0.035 &   29    &  1.0e-4 & 0.046 &   38    \\ \hline 
B&  1.0e-5 & 0.023 &   19    &  1.0e-6 & 0.031 &   26    &  1.0e-4 & 0.042 &   35    &  1.0e-4 & 0.037 &   32    \\ \hline 
C&  1.0e-6 & 0.039 &   38    &  1.0e-6 & 0.044 &   35    &  1.0e-3 & 0.036 &   34    &  1.0e-4 & 0.056 &   52    \\ \hline 
D&  1.0e-6 & 0.068 &   68    &  1.0e-6 & 0.067 &   60    &  1.0e-8 & 0.38 &   464    &  1.0e-8 & 0.34 &   381    \\ \hline 
E&  2.6e-7 & 23 &   65536    &  2.6e-7 & 26 &   65536    &  3.6e-7 & 93 &   262144    &  3.8e-7 & 1.1e+02 &   262144    \\ \hline 
F&  1.0e-06 & 0.028 &   25    &  1.0e-06 & 0.033 &   29    &  1.0e-03 & 0.028 &   26    &  1.0e-05 & 0.051 &   47    \\ \hline
G& 1.0e-06&0.0422&42&1.0e-06&0.0476&43& 1.0e-04&0.0431&43&1.0e-04&0.0451&41\\ \hline
\end{tabular}
\end{tiny}
\captionsetup{font=scriptsize}
\caption{
Efficiency comparison for several methods: A- the $\mu(x)$ method, B- the adjoint $\mu(x)$ method, C- the hybrid method, D- the method of Lee and Stewart (using an adaptive mesh in $z$), E- the method of Lee and Stewart using a 4th order fixed mesh Runge Kutta solver, F- the polar adjoint method, and G- the polar adjoint radial method. Here $\mathcal{E}=10$, $q=0.3$, $C=10^{21\mathcal{E}/40}$, and $\phi(u)=Ce^{-\mathcal{E}/T(u)}$, $T(u)=1-(u-1.5)^2$.
\label{E10q3phi1} 
}
\end{table}

%

\begin{table}[!th]\begin{tiny}
\begin{tabular}{|c||c|c|c||c|c|c||c|c|c||c|c|c|}
\hline
 &$\lambda=1$& & &$\lambda=i$& & &$\lambda=10$& & &$\lambda=10i$ & &\\
\hline
  & Tol & Time & pnts& Tol & Time & Mesh & Tol & Time & Mesh & Tol & Time & Mesh\\
 \hline
A &  1.0e-5 & 0.032 &   31    &  1.0e-6 & 0.043 &   40    &  1.0e-4 & 0.056 &   57    &  1.0e-4 & 0.068 &   63    \\ \hline 
B&  1.0e-6 & 0.04 &   42    &  1.0e-7 & 0.058 &   58    &  1.0e-4 & 0.052 &   54    &  1.0e-5 & 0.083 &   82    \\ \hline 
C&  1.0e-5 & 0.036 &   37    &  1.0e-6 & 0.047 &   46    &  1.0e-4 & 0.056 &   63    &  1.0e-4 & 0.071 &   71    \\ \hline 
D&  1.0e-6 & 0.074 &   85    &  1.0e-7 & 0.086 &   94    &  1.0e-8 & 0.72 &   895    &  1.0e-8 & 0.79 &   890    \\ \hline 
E&  6.5e-7 & 12 &   32768    &  1.2e-7 & 26 &   65536    &  7.9e-7 & 370 &   $\approx 10^6$   &  6.9e-7 & 420 &   $\approx 10^6$    \\ \hline
F&  1.0e-06 & 0.04 &   39    &  1.0e-07 & 0.06 &   61    &  1.0e-04 & 0.053 &   59    &  1.0e-05 & 0.079 &   82    \\ \hline  
G&1.0e-07&0.0673&69&1.0e-07&0.0809&74&1.0e-05&0.0825&84&1.0e-05&0.0977&90\\ \hline
\end{tabular}
\end{tiny}
\captionsetup{font=scriptsize}
\caption{
Efficiency comparison for several methods: A- the $\mu(x)$ method, B- the adjoint $\mu(x)$ method, C- the hybrid method, D- the method of Lee and Stewart (using an adaptive mesh in $z$), E- the method of Lee and Stewart using a 4th order fixed mesh Runge Kutta solver, F- the polar adjoint method, and G- the polar adjoint radial method. Here $\mathcal{E}=10$, $q=0.47$, $C=10^{21\mathcal{E}/40}$, and $\phi(u)=Ce^{-\mathcal{E}/T(u)}$, $T(u)=1-(u-1.5)^2$.
\label{E10q47phi1} 
}
\end{table}

\begin{table}[!th]\begin{tiny}
\begin{tabular}{|c||c|c|c||c|c|c||c|c|c||c|c|c|}
\hline
 &$\lambda=1$& & &$\lambda=i$& & &$\lambda=10$& & &$\lambda=10i$ & &\\
\hline
  & Tol & Time & pnts& Tol & Time & Mesh & Tol & Time & Mesh & Tol & Time & Mesh\\
 \hline
 A&-&-&-&-&-&-&-&-&-&-&-&-\\ \hline
 B&  1.0e-02 & 0.0027 &   2    &  1.0e-02 & 0.0028 &   2    &  1.0e-02 & 0.0027 &   2    &  1.0e-02 & 0.0027 &   2    \\ \hline 
 C&  1.0e-02 & 0.062 &   70    &  1.0e-02 & 0.14 &   149    &  1.0e-02 & 0.54 &   617    &  1.0e-02 & 1.4 &   1494    \\ \hline 
F&  1.0e-02 & 0.0037 &   2    &  1.0e-02 & 0.0027 &   2    &  1.0e-02 & 0.0025 &   2    &  1.0e-02 & 0.0026 &   2    \\ \hline  
G& 1.0e-02&0.0045&3&1.0e-02&0.0059&3&1.0e-02&0.0061&3&1.0e-02&0.0044&3 \\
\hline
\end{tabular}
\end{tiny}
\captionsetup{font=scriptsize}
\caption{
Efficiency comparison for several methods: A- the $\mu(x)$ method, B- the adjoint $\mu(x)$ method, C- the hybrid method, D- the method of Lee and Stewart (using an adaptive mesh in $z$), E- the method of Lee and Stewart using a 4th order fixed mesh Runge Kutta solver, F- the polar adjoint method, and G- the polar adjoint radial method. Here $\mathcal{E}=40$, $q=0$, $C=10^{21\mathcal{E}/40}$, and $\phi(u)=Ce^{-\mathcal{E}/T(u)}$, $T(u)=1-(u-1.5)^2$.
\label{E40q0phi1}
}
\end{table}

%

\begin{table}[!th]\begin{tiny}
\begin{tabular}{|c||c|c|c||c|c|c||c|c|c||c|c|c|}
\hline
 &$\lambda=1$& & &$\lambda=i$& & &$\lambda=10$& & &$\lambda=10i$ & &\\
\hline
  & Tol & Time & pnts& Tol & Time & Mesh & Tol & Time & Mesh & Tol & Time & Mesh\\
 \hline
 A&  1.0e-07 & 0.096 &   103    &  1.0e-07 & 0.12 &   112    &  1.0e-08 & 0.38 &   450    &  1.0e-07 & 0.37 &   394    \\ \hline 
 B&  1.0e-08 & 0.18 &   203    &  1.0e-07 & 0.13 &   134    &  1.0e-08 & 0.42 &   503    &  1.0e-07 & 0.64 &   687    \\ \hline 
 C&  1.0e-07 & 0.1 &   117    &  1.0e-07 & 0.14 &   145    &  1.0e-08 & 0.4 &   487    &  1.0e-07 & 0.79 &   877    \\ \hline 
 D&  1.0e-08 & 0.2 &   246    &  1.0e-07 & 0.13 &   143    &  1.0e-08 & 0.84 &   1068    &  1.0e-08 & 0.93 &   1063    \\ \hline 
E &  2.17e-06 & 5.9 &   16384    &  2.06e-06 & 6.6 &   16384    &  4.36e-07 & 12 &   32768    &  3.79e-07 & 13 &   32768    \\ \hline 
F &  1.0e-03 & 0.023 &   21    &  1.0e-07 & 0.12 &   118    &  1.0e-03 & 0.062 &   67    &  1.0e-07 & 0.64 &   687    \\ \hline 
G& 1.0e-07&0.1338&137&1.0e-07&0.1720&158&1.0e-07&0.3172&334&1.0e-07&0.7689&718 \\ \hline
\end{tabular}
\end{tiny}
\captionsetup{font=scriptsize}
\caption{
Efficiency comparison for several methods: A- the $\mu(x)$ method, B- the adjoint $\mu(x)$ method, C- the hybrid method, D- the method of Lee and Stewart (using an adaptive mesh in $z$), E- the method of Lee and Stewart using a 4th order fixed mesh Runge Kutta solver, F- the polar adjoint method, and G- the polar adjoint radial method. Here $\mathcal{E}=40$, $q=0.3$, $C=10^{21\mathcal{E}/40}$, and $\phi(u)=Ce^{-\mathcal{E}/T(u)}$, $T(u)=1-(u-1.5)^2$.
\label{E40q3phi1} 
}
\end{table}

%

\begin{table}[!th]\begin{tiny}
\begin{tabular}{|c||c|c|c||c|c|c||c|c|c||c|c|c|}
\hline
 &$\lambda=1$& & &$\lambda=i$& & &$\lambda=10$& & &$\lambda=10i$ & &\\
\hline
  & Tol & Time & pnts& Tol & Time & Mesh & Tol & Time & Mesh & Tol & Time & Mesh\\
 \hline
 A&  1.0e-08 & 0.14 &   159    &  1.0e-08 & 0.17 &   168    &  1.0e-07 & 0.22 &   250    &  1.0e-08 & 0.46 &   493    \\ \hline 
 B &  1.0e-08 & 0.19 &   221    &  1.0e-08 & 0.21 &   226    &  1.0e-08 & 0.37 &   445    &  1.0e-08 & 0.87 &   944    \\ \hline 
 C&  1.0e-08 & 0.15 &   176    &  1.0e-08 & 0.19 &   199    &  1.0e-07 & 0.23 &   273    &  1.0e-08 & 0.85 &   942    \\ \hline 
 D&  1.0e-08 & 0.22 &   273    &  1.0e-07 & 0.15 &   163    &  1.0e-08 & 0.63 &   793    &  1.0e-08 & 0.68 &   765    \\ \hline 
E&  9.82e-07 & 12 &   32768    &  8.75e-07 & 13 &   32768    &  1.78e-07 & 23 &   65536    &  1.62e-07 & 26 &   65536    \\ \hline 
F&  1.0e-03 & 0.026 &   23    &  1.0e-06 & 0.08 &   76    &  1.0e-03 & 0.054 &   57    &  1.0e-07 & 0.55 &   585    \\ \hline
G& 1.0e-08&0.2300&242&1.0e-07&0.1879&174&10e-08&0.4368&458&1.0e-07&0.6668&620 \\ \hline
\end{tabular}
\end{tiny}
\captionsetup{font=scriptsize}
\caption{
Efficiency comparison for several methods: A- the $\mu(x)$ method, B- the adjoint $\mu(x)$ method, C- the hybrid method, D- the method of Lee and Stewart (using an adaptive mesh in $z$), E- the method of Lee and Stewart using a 4th order fixed mesh Runge Kutta solver, F- the polar adjoint method, and G- the polar adjoint radial method. Here $\mathcal{E}=40$, $q=0.47$, $C=10^{21\mathcal{E}/40}$, and $\phi(u)=Ce^{-\mathcal{E}/T(u)}$, $T(u)=1-(u-1.5)^2$.
\label{E40q47phi1} 
}
\end{table}

\clearpage

\begin{table}[!th]\begin{tiny}
\begin{tabular}{|c||c|c|c||c|c|c||c|c|c||c|c|c|}
\hline
 &$\lambda=1$& & &$\lambda=i$& & &$\lambda=10$& & &$\lambda=10i$ & &\\
\hline
  & Tol & Time & pnts& Tol & Time & Mesh & Tol & Time & Mesh & Tol & Time & Mesh\\
 \hline
 A &  1.0e-03 & 0.011 &   9    &  1.0e-03 & 0.011 &   8    &  1.0e-02 & 0.026 &   24    &  1.0e-03 & 0.037 &   33    \\ \hline 
 B  &  1.0e-03 & 0.0082 &   6    &  1.0e-04 & 0.011 &   8    &  1.0e-02 & 0.021 &   21    &  1.0e-04 & 0.038 &   38    \\ \hline 
 C   &  1.0e-03 & 0.011 &   9    &  1.0e-04 & 0.016 &   12    &  1.0e-02 & 0.026 &   26    &  1.0e-04 & 0.059 &   60    \\ \hline 
 D   &  1.0e-04 & 0.029 &   29    &  1.0e-07 & 0.039 &   43    &  1.0e-08 & 0.34 &   423    &  1.0e-08 & 0.37 &   419    \\ \hline 
E &  7.2e-07 & 0.046 &   128     &  8.2e-07 & 0.052 &   128    & 1.8e-07 & 1.5 &   4096     &  1.8e-07 & 1.6 &   4096    \\ \hline
F &  1.0e-03 & 0.0092 &   6    &  1.0e-05 & 0.011 &   9    &  1.0e-02 & 0.02 &   21    &  1.0e-04 & 0.039 &   38    \\ \hline 
G & 1.0e-04&0.0183&18&1.0e-04&0.0206&18&1.0e-04&0.067&70\\ \hline
\end{tabular}
\end{tiny}
\captionsetup{font=scriptsize}
\caption{
	Efficiency comparison for several methods: A- the $\mu(x)$ method, B- the adjoint $\mu(x)$ method, C- the hybrid method, D- the method of Lee and Stewart (using an adaptive mesh in $z$), E- the method of Lee and Stewart using a 4th order fixed mesh Runge Kutta solver, F- the polar adjoint method, and G- the polar adjoint radial method. Here $\E=1$, $q=0.1$, $C=e^{\mathcal{E}/2}$, and $\phi(u)=Ce^{-\mathcal{E}/u}$.
	\label{E1q1phi2} 
	}
\end{table}

\begin{table}[!th]\begin{tiny}
\begin{tabular}{|c||c|c|c||c|c|c||c|c|c||c|c|c|}
\hline
 &$\lambda=1$& & &$\lambda=i$& & &$\lambda=10$& & &$\lambda=10i$ & &\\
\hline
  & Tol & Time & pnts& Tol & Time & Mesh & Tol & Time & Mesh & Tol & Time & Mesh\\
 \hline
 A  &  1.0e-02 & 0.011 &   9    &  1.0e-03 & 0.013 &   10    &  1.0e-02 & 0.028 &   28    &  1.0e-03 & 0.044 &   44    \\ \hline 
 B  &  1.0e-04 & 0.0099 &   10    &  1.0e-04 & 0.011 &   10    &  1.0e-02 & 0.024 &   26    &  1.0e-05 & 0.086 &   91    \\ \hline 
 C   &  1.0e-02 & 0.012 &   10    &  1.0e-04 & 0.016 &   13    &  1.0e-02 & 0.027 &   30    &  1.0e-04 & 0.074 &   79    \\ \hline 
 D &  1.0e-06 & 0.036 &   44    &  1.0e-08 & 0.059 &   66    &  1.0e-08 & 0.43 &   532    &  1.0e-08 & 0.47 &   532    \\ \hline 
E &  6.4e-07 & 0.092 &   256     &  7e-07 & 0.1 &   256&  7.6e-08 & 2.9 &   8192     &  7.8e-08 & 3.3 &   8192 \\ \hline
F&  1.0e-04 & 0.012 &   11    &  1.0e-04 & 0.011 &   10    &  1.0e-04 & 0.028 &   32    &  1.0e-05 & 0.086 &   91    \\ \hline 
G& 1.0e-04&0.0251&25&1.0e-04&0.0280&25&1.0e-05&0.1018&106&1.0e-05&0.1326&122 \\
\hline
\end{tabular}
\end{tiny}
\captionsetup{font=scriptsize}
\caption{
Efficiency comparison for several methods: A- the $\mu(x)$ method, B- the adjoint $\mu(x)$ method, C- the hybrid method, D- the method of Lee and Stewart (using an adaptive mesh in $z$), E- the method of Lee and Stewart using a 4th order fixed mesh Runge Kutta solver, F- the polar adjoint method, and G- the polar adjoint radial method. Here $\E=1$, $q=0.3$, $C=e^{\mathcal{E}/2}$, and $\phi(u)=Ce^{-\mathcal{E}/u}$.
\label{E1q3phi2}
}
\end{table}

\begin{table}[!th]\begin{tiny}
\begin{tabular}{|c||c|c|c||c|c|c||c|c|c||c|c|c|}
\hline
 &$\lambda=1$& & &$\lambda=i$& & &$\lambda=10$& & &$\lambda=10i$ & &\\
\hline
  & Tol & Time & pnts& Tol & Time & Mesh & Tol & Time & Mesh & Tol & Time & Mesh\\
 \hline
 A   &  1.0e-03 & 0.014 &   14    &  1.0e-04 & 0.021 &   19    &  1.0e-02 & 0.035 &   38    &  1.0e-03 & 0.075 &   79    \\ \hline 
 B  &  1.0e-05 & 0.02 &   23    &  1.0e-05 & 0.024 &   25    &  1.0e-03 & 0.034 &   40    &  1.0e-05 & 0.17 &   185    \\ \hline 
 C   &  1.0e-03 & 0.016 &   16    &  1.0e-05 & 0.028 &   28    &  1.0e-02 & 0.034 &   40    &  1.0e-04 & 0.13 &   137    \\ \hline 
 D  &  1.0e-07 & 0.063 &   76    &  1.0e-08 & 0.086 &   96    &  1.0e-08 & 0.67 &   840    &  1.0e-08 & 0.75 &   854    \\ \hline 
E&  1.5e-07 & 0.37 &   1024     &  1.7e-07 & 0.41 &   1024&  2.5e-07 & 5.9 &   16384     &  2.6e-07 & 6.5 &   16384    \\ \hline 
F&  1.0e-05 & 0.02 &   23    &  1.0e-06 & 0.035 &   36    &  1.0e-04 & 0.044 &   52    &  1.0e-05 & 0.17 &   185    \\ \hline 
G& 1.0e-05&0.0493&51&1.0e-06&0.0783&72&1.0e-06&0.2847&297&1.0e-05&0.3128&290\\ \hline
\end{tabular}
\end{tiny}
\captionsetup{font=scriptsize}
\caption{
Efficiency comparison for several methods: A- the $\mu(x)$ method, B- the adjoint $\mu(x)$ method, C- the hybrid method, D- the method of Lee and Stewart (using an adaptive mesh in $z$), E- the method of Lee and Stewart using a 4th order fixed mesh Runge Kutta solver, F- the polar adjoint method, and G- the polar adjoint radial method. Here $\E=1$, $q=0.47$, $C=e^{\mathcal{E}/2}$, and $\phi(u)=Ce^{-\mathcal{E}/u}$.
\label{E1q47phi2} 
}
\end{table}

\begin{table}[!th]\begin{tiny}
\begin{tabular}{|c||c|c|c||c|c|c||c|c|c||c|c|c|}
\hline
 &$\lambda=1$& & &$\lambda=i$& & &$\lambda=10$& & &$\lambda=10i$ & &\\
\hline
  & Tol & Time & pnts& Tol & Time & Mesh & Tol & Time & Mesh & Tol & Time & Mesh\\
 \hline
 A &  1.0e-02 & 0.01 &   8    &  1.0e-04 & 0.016 &   13    &  1.0e-02 & 0.027 &   27    &  1.0e-03 & 0.046 &   45    \\ \hline 
 B  &  1.0e-03 & 0.0089 &   8    &  1.0e-04 & 0.012 &   10    &  1.0e-03 & 0.024 &   26    &  1.0e-04 & 0.052 &   55    \\ \hline 
 C    &  1.0e-02 & 0.0096 &   8    &  1.0e-04 & 0.017 &   14    &  1.0e-02 & 0.027 &   29    &  1.0e-04 & 0.069 &   74    \\ \hline 
 D  &  1.0e-06 & 0.038 &   46    &  1.0e-07 & 0.044 &   49    &  1.0e-08 & 0.4 &   497    &  1.0e-08 & 0.44 &   500    \\ \hline 
E &  3.3e-07 & 0.092 &   256     &  3.3e-07 & 0.1 &   256&  5.6e-07 & 1.5 &   4096     &  5.9e-07 & 1.6 &   4096    \\ \hline 
F&  1.0e-03 & 0.0092 &   8    &  1.0e-04 & 0.012 &   10    &  1.0e-02 & 0.024 &   25    &  1.0e-04 & 0.053 &   55    \\ \hline 
G& 10e-04&0.0202&20&1.0e-04&0.0240&21&1.0e-04&0.0811&84&1.0e-04&01229&113\\ \hline
\end{tabular}
\end{tiny}
\captionsetup{font=scriptsize}
\caption{
Efficiency comparison for several methods: A- the $\mu(x)$ method, B- the adjoint $\mu(x)$ method, C- the hybrid method, D- the method of Lee and Stewart (using an adaptive mesh in $z$), E- the method of Lee and Stewart using a 4th order fixed mesh Runge Kutta solver, F- the polar adjoint method, and G- the polar adjoint radial method. Here $\E=10$, $q=0.1$, $C=e^{\mathcal{E}/2}$, and $\phi(u)=Ce^{-\mathcal{E}/u}$.
\label{E10q1phi2}
}
\end{table}

\begin{table}[!th]\begin{tiny}
\begin{tabular}{|c||c|c|c||c|c|c||c|c|c||c|c|c|}
\hline
 &$\lambda=1$& & &$\lambda=i$& & &$\lambda=10$& & &$\lambda=10i$ & &\\
\hline
  & Tol & Time & pnts& Tol & Time & Mesh & Tol & Time & Mesh & Tol & Time & Mesh\\
 \hline
 A &  1.0e-03 & 0.014 &   13    &  1.0e-05 & 0.035 &   32    &  1.0e-02 & 0.047 &   50    &  1.0e-04 & 0.14 &   155    \\ \hline 
 B &  1.0e-02 & 0.0082 &   8    &  1.0e-05 & 0.027 &   28    &  1.0e-03 & 0.043 &   49    &  1.0e-05 & 0.21 &   228    \\ \hline 
 C   &  1.0e-02 & 0.013 &   12    &  1.0e-05 & 0.033 &   32    &  1.0e-02 & 0.045 &   52    &  1.0e-05 & 0.24 &   260    \\ \hline 
 D &  1.0e-07 & 0.076 &   93    &  1.0e-08 & 0.1 &   116    &  1.0e-08 & 0.81 &   1027    &  1.0e-08 & 0.9 &   1037    \\ \hline 
E &  2.6e-07 & 0.37 &   1024     &  3e-07 & 0.41 &   1024&  4.5e-07 & 5.8 &   16384     &  4.8e-07 & 6.6 &   16384    \\ \hline
F&1.0e-04 & 0.013 &   14    &  1.0e-06 & 0.041 &   42    &  1.0e-05 & 0.062 &   71    &  1.0e-05 & 0.21 &   228    \\ \hline
G&1.0e-03&0.031&32&1.0e-05&0.0528&48&1.0e-06&0.2576&269&1.0e-05&0.3669&339\\ \hline
\end{tabular}
 \end{tiny}
\captionsetup{font=scriptsize}
\caption{
Efficiency comparison for several methods: A- the $\mu(x)$ method, B- the adjoint $\mu(x)$ method, C- the hybrid method, D- the method of Lee and Stewart (using an adaptive mesh in $z$), E- the method of Lee and Stewart using a 4th order fixed mesh Runge Kutta solver, F- the polar adjoint method, and G- the polar adjoint radial method. Here $\E=10$, $q=0.3$, $C=e^{\mathcal{E}/2}$, and $\phi(u)=Ce^{-\mathcal{E}/u}$.
\label{E10q3phi2} 
}
\end{table}

\begin{table}[!th]\begin{tiny}
\begin{tabular}{|c||c|c|c||c|c|c||c|c|c||c|c|c|}
\hline
 &$\lambda=1$& & &$\lambda=i$& & &$\lambda=10$& & &$\lambda=10i$ & &\\
\hline
  & Tol & Time & pnts& Tol & Time & Mesh & Tol & Time & Mesh & Tol & Time & Mesh\\ \hline
 A &  1.0e-03 & 0.029 &   31    &  1.0e-08 & 0.5 &   540    &  1.0e-02 & 0.14 &   174    &  1.0e-09 & 7.3 &   8009    \\ \hline 
 B &  1.0e-02 & 0.019 &   22    &  1.0e-08 & 0.51 &   542    &  1.0e-04 & 0.16 &   186    &  1.0e-08 & 4.7 &   5160    \\ \hline 
 C   &  1.0e-03 & 0.028 &   31    &  1.0e-08 & 0.5 &   542    &  1.0e-02 & 0.14 &   176    &  1.0e-08 & 4.7 &   5169    \\ \hline 
 D &  1.0e-08 & 0.38 &   482    &  1.0e-08 & 0.42 &   481    &  1.0e-09 & 5.9 &   7352    &  1.0e-09 & 6.6 &   7464    \\ \hline 
E &  1.6e-07 & 5.9 &   16384     &  1.4e-07 & 6.6 &   16384&&&&&& \\ \hline 
F&  1.0e-03 & 0.026 &   29    &  1.0e-06 & 0.19 &   201    &  1.0e-04 & 0.21 &   245    &  1.0e-06 & 1.8 &   1933    \\ \hline 
G& 1.0e-07&0.3950&411&1.0e-06&0.6063&565&1.0e-07&3.1675&3218&1.0e-06&6.0556&5575\\ \hline
\end{tabular}
\end{tiny}
\captionsetup{font=scriptsize}
\caption{
Efficiency comparison for several methods: A- the $\mu(x)$ method, B- the adjoint $\mu(x)$ method, C- the hybrid method, D- the method of Lee and Stewart (using an adaptive mesh in $z$), E- the method of Lee and Stewart using a 4th order fixed mesh Runge Kutta solver, F- the polar adjoint method, and G- the polar adjoint radial method. Here $\E=10$, $q=0.47$, $C=e^{\mathcal{E}/2}$, and $\phi(u)=Ce^{-\mathcal{E}/u}$.
\label{E10q47phi2} 
}
\end{table}

\begin{table}[!th]\begin{tiny}
\begin{tabular}{|c||c|c|c||c|c|c||c|c|c||c|c|c|}
\hline
 &$\lambda=1$& & &$\lambda=i$& & &$\lambda=10$& & &$\lambda=10i$ & &\\
\hline
  & Tol & Time & pnts& Tol & Time & Mesh & Tol & Time & Mesh & Tol & Time & Mesh\\
 \hline
 A &  1.0e-02 & 0.011 &   9    &  1.0e-05 & 0.022 &   19    &  1.0e-02 & 0.032 &   33    &  1.0e-04 & 0.08 &   84    \\ \hline 
 B &  1.0e-02 & 0.0081 &   7    &  1.0e-05 & 0.016 &   16    &  1.0e-04 & 0.031 &   35    &  1.0e-05 & 0.11 &   113    \\ \hline 
 C   &  1.0e-02 & 0.011 &   10    &  1.0e-05 & 0.021 &   20    &  1.0e-02 & 0.033 &   36    &  1.0e-05 & 0.14 &   148    \\ \hline 
 D &  1.0e-06 & 0.041 &   50    &  1.0e-07 & 0.048 &   54    &  1.0e-08 & 0.48 &   595    &  1.0e-08 & 0.53 &   597    \\ \hline 
E &  8.1e-08 & 0.18 &   512     &  8e-08 & 0.21 &   512&&&&&&\\ \hline
F&  1.0e-03 & 0.0083 &   8    &  1.0e-03 & 0.0092 &   8    &  1.0e-02 & 0.027 &   30    &  1.0e-04 & 0.07 &   74    \\ \hline 
G&1.0e-04&0.0232&23&1.0e-04&0.0282&25&1.0e-05&0.1116&115&1.0e-05&0.1686&155\\ \hline
\end{tabular}
  \end{tiny}
\captionsetup{font=scriptsize}
\caption{
Efficiency comparison for several methods: A- the $\mu(x)$ method, B- the adjoint $\mu(x)$ method, C- the hybrid method, D- the method of Lee and Stewart (using an adaptive mesh in $z$), E- the method of Lee and Stewart using a 4th order fixed mesh Runge Kutta solver, F- the polar adjoint method, and G- the polar adjoint radial method. Here $\E=20$, $q=0.1$, $C=e^{\mathcal{E}/2}$, and $\phi(u)=Ce^{-\mathcal{E}/u}$.
\label{E20q1phi2} 
}
\end{table}

\begin{table}[!th]\begin{tiny}
\begin{tabular}{|c||c|c|c||c|c|c||c|c|c||c|c|c|}
\hline
 &$\lambda=1$& & &$\lambda=i$& & &$\lambda=10$& & &$\lambda=10i$ & &\\
\hline
  & Tol & Time & pnts& Tol & Time & Mesh & Tol & Time & Mesh & Tol & Time & Mesh\\
 \hline
 A &  1.0e-03 & 0.021 &   22    &  1.0e-09 & 0.4 &   418    &  1.0e-02 & 0.093 &   106    &  1.0e-09 & 3.5 &   3852    \\ \hline 
 B &  1.0e-04 & 0.022 &   26    &  1.0e-09 & 0.41 &   424    &  1.0e-03 & 0.091 &   107    &  1.0e-08 & 2.3 &   2529    \\ \hline 
 C   &  1.0e-04 & 0.026 &   30    &  1.0e-09 & 0.4 &   423    &  1.0e-02 & 0.09 &   108    &  1.0e-08 & 2.3 &   2532    \\ \hline 
 D &  1.0e-08 & 0.21 &   265    &  1.0e-08 & 0.22 &   246    &  1.0e-09 & 2.9 &   3668    &  1.0e-09 & 3.2 &   3705    \\ \hline 
E &  1.8e-07 & 1.5 &   4096     &  1.6e-07 & 1.6 &   4096&&&&&&\\ \hline
F& 1.0e-03 & 0.017 &   19    &  1.0e-06 & 0.094 &   100    &  1.0e-04 & 0.11 &   128    &  1.0e-05 & 0.54 &   583    \\ \hline 
G&1.0e-05&0.0874&90&1.0e-06&0.1626&149&1.0e-06&0.7249&740&1.0e-06&1.3564&1255\\ \hline
\end{tabular}
   \end{tiny}
\captionsetup{font=scriptsize}
\caption{
Efficiency comparison for several methods: A- the $\mu(x)$ method, B- the adjoint $\mu(x)$ method, C- the hybrid method, D- the method of Lee and Stewart (using an adaptive mesh in $z$), E- the method of Lee and Stewart using a 4th order fixed mesh Runge Kutta solver, F- the polar adjoint method, and G- the polar adjoint radial method. Here $\E=20$, $q=0.3$, $C=e^{\mathcal{E}/2}$, and $\phi(u)=Ce^{-\mathcal{E}/u}$.
\label{E20q3phi2} 
}
\end{table}

\clearpage

\begin{table}[!th]\begin{tiny}
\begin{tabular}{|c||c|c|c||c|c|c||c|c|c||c|c|c|}
\hline
 &$\lambda=1$& & &$\lambda=i$& & &$\lambda=10$& & &$\lambda=10i$ & &\\
\hline
  & Tol & Time & pnts& Tol & Time & Mesh & Tol & Time & Mesh & Tol & Time & Mesh\\
 \hline
 A &  1.0e-03 & 0.13 &   156    &  1.0e-09 & 7.8 &   8654    &  1.0e-02 & 1.2 &   1412    &  1.0e-10 & 1.3e+02 &   136649    \\ \hline 
 B &  1.0e-05 & 0.18 &   221    &  1.0e-09 & 8 &   8670    &  1.0e-04 & 1.2 &   1430    &  1.0e-10 & 1.3e+02 &   136946    \\ \hline 
 C   &  1.0e-04 & 0.15 &   175    &  1.0e-09 & 7.9 &   8657    &  1.0e-03 & 1.2 &   1419    &  1.0e-10 & 1.3e+02 &   136948    \\ \hline 
 D &-&-&-&1.0e-10&319.6&258949&-&-&-&1.0e-11&9270.4&3.9e06 \\ \hline
E &  5.9e-07 & 96 &   262144     &  5.2e-08 & 1.1e+02 &   262144&&&&&& \\ \hline
F  &  1.0e-04 & 0.18 &   218    &  1.0e-07 & 2.7 &   2882    &  1.0e-03 & 2.1 &   2023    &  1.0e-07 & 27 &   28684    \\ \hline 
G&1.0e-07&5.8414&5911&1.0e-07&18.2622&16857&-&-&-&-&-&- \\ \hline
\end{tabular}
\end{tiny}
\captionsetup{font=scriptsize}
\caption{
Efficiency comparison for several methods: A- the $\mu(x)$ method, B- the adjoint $\mu(x)$ method, C- the hybrid method, D- the method of Lee and Stewart (using an adaptive mesh in $z$), E- the method of Lee and Stewart using a 4th order fixed mesh Runge Kutta solver, F- the polar adjoint method, and G- the polar adjoint radial method. Here $\E=20$, $q=0.47$, $C=e^{\mathcal{E}/2}$, and $\phi(u)=Ce^{-\mathcal{E}/u}$.
\label{E20q47phi2} 
}
\end{table}


\medskip
{\bf Acknowledgement.}
K.Z. thanks the University of Paris 13 for their hospitality
during a visit in which this work was partly carried out.
The numerical Evans function computations performed
in this paper were carried out
using the STABLAB package developed by Jeffrey Humpherys with
help of the authors.





\end{document}